\documentclass[11pt]{article}
\usepackage{amssymb,amsmath,latexsym}
\usepackage[mathscr]{eucal}
 \usepackage{amsfonts,verbatim,graphicx, subfigure,longtable}

\newtheorem{thm}{Theorem}[section]

\newtheorem{lem}[thm]{Lemma}
\newtheorem{prop}[thm]{Proposition}

\newtheorem{remark}[thm]{Remark}
\numberwithin{equation}{section}

\def\cS{{\cal S}}

\def\E{{\mathbb E}}
\def\P{{\mathbb P}}
\def\HH{{\mathbb H}}
\def\RR{{\mathbb R}}
\def\CC{{\mathbb C}}

\def\QQ{{\mathbb Q}}

\def\D{{\bf D}}
\def\F{{\bf F}}
\def\W{{\bf X}}
\def\W{{\bf W}}

\def\n{{\bf n}}
\def\w{{\bf w}}
\def\dd{{\cal D}}
\def\gg{{\cal G}}

\def\pf{\noindent{\bf Proof.} }
\def\qed{{\hfill $\Box$ \bigskip}}
\def\wh{\widehat}
\def\wt{\widetilde}
\def\dis{\displaystyle}
\def\eps{\varepsilon}
\def\O{\Omega}
\def\c{\check}

\def\b{{\mathscr B}}

\def\1{{\bf 1}}
\def\W{{\bf W}}
\def\x{{\bf x}}
\def\y{{\bf y}}
\def\s{{\bf s}}

\def\H{{\bf H}}

\def\<{\langle}
\def\>{\rangle}
\def\st{\stackrel{\circ}}
\def\hc{{\partial B_\epsilon\cap\mathbb H}}

\begin{document}
\title{\bf  Stochastic Komatu-Loewner evolutions\\
  and  BMD domain constant}
\author{{\bf Zhen-Qing Chen\thanks{Research partially supported
by NSF Grant DMS-1206276}} \quad  and \quad {\bf Masatoshi Fukushima}}
\date{(April 26, 2016)}
\maketitle

\begin{abstract}

 Let $D={\mathbb H} \setminus \cup_{k=1}^N C_k$ be a standard slit domain, where $\mathbb H$ is the upper half plane and $C_k$, $1\leq k\leq N$, are mutually disjoint horizontal line segments in $\HH$.  Given a Jordan arc $\gamma\subset D$ starting at $\partial {\mathbb H},$ let $g_t$ be the unique conformal map from $D\setminus \gamma[0,t]$ onto a standard slit domain $D_t=\HH\setminus \cup_{k=1}^N C_k(t)$ satisfying the hydrodynamic normalization at infinity.
 It has been
  established recently
  that $g_t$ satisfies an ODE called a Komatu-Loewner equation
  in terms of
  the complex Poisson kernel of the Brownian motion with darning (BMD) for $D_t$,   extending the classical
  chordal Loewner equation for the simply connected domain $D={\mathbb H}.$

We randomize
the Jordan  arc  $\gamma$ according to
 a system of probability measures on the family of equivalence classes of Jordan arcs that enjoy
a domain Markov property and a certain conformal invariance property.
We show
that the induced process $(\xi(t), \s(t))$
satisfies a Markov type stochastic differential equation,
where $\xi(t)$ is a motion on $\partial\HH$ and $\s(t)$ represents the motion of
the endpoints of the
slits $\{C_k(t),\; 1\le k\le N \}.$
The diffusion and drift coefficients $\alpha$ and $b$ of $\xi(t)$ are homogeneous functions of degree $0$ and $-1$, respectively, while $\s(t)$ has   drift coefficients only, determined by the BMD complex Poisson kernel that are known
 to be Lipschitz continuous.

Conversely, given
such functions $\alpha$ and $b$ with local Lipschitz continuity,
the corresponding SDE admits
a unique solution $(\xi(t), \s(t))$.
The latter produces random conformal maps $g_t(z)$ via the Komatu-Loewner
equation.
The resulting family of random growing hulls
  $\{F_t\}$ from the conformal mappings is called ${\rm SKLE}_{\alpha,b}.$
We show that it enjoys a certain scaling property and
a domain Markov property.
Among other things,
we further prove that
${\rm SKLE}_{\alpha,-b_{\rm BMD}}$ for a constant
$\alpha >0$ 
has a locality property if and only if $\alpha = 
\sqrt{6}$,
where $b_{\rm BMD}$ is
a BMD-domain constant that describes the
discrepancy of a standard slit domain from $\HH$ relative to BMD.
\end{abstract}

\medskip
\noindent
{\bf AMS 2010 Mathematics Subject Classification}: Primary 60J67, 60J70; Secondary 30C20, 60H10, 60H30

\smallskip\noindent
{\bf Keywords and phrases}: Stochastic Komatu-Loewner evolution, Brownian motion with darning, Komatu-Loewner equation for slits, SDE with homogeneous coefficients,   
generalized Komatu-Loewner equation for image hulls, 
BMD domain constant, locality property

{\small
 \begin{tableofcontents}
 \end{tableofcontents} }

\section{Introduction}
In 2000, Oded Schramm \cite{S}
introduced a
{\it stochastic Loewner evolution }(SLE) on the upper half plane $\HH$ with
driving process
 $\xi(t)=\sqrt{\kappa}B_t$, where $B_t$ is the standard Brownian motion on $\partial \HH$ and $\kappa$ is a positive constant.
 The solution of the SLE is a family of random conformal
 mappings from $\HH \setminus K_t$ to $\HH$ indexed by $t\geq 0$.
 The increasing family of random hulls $\{K_t; t\geq 0\}$
 is nowadays called ${\rm SLE}_\kappa$.
It has a certain conformal invariance and a domain Markov property.
 ${\rm SLE}_\kappa$   is a powerful tool
 in studying two-dimensional critical systems  in statistical physics.
 ${\rm SLE}_\kappa$ has been proved to be the scaling limit
 of various critical two-dimensional lattice models,
 such as loops erased random walk,
 uniform spanning trees, critical percolation, critical Ising model,
 and has been conjectured for a few more including self-avoiding random walks.
 In particular, ${\rm SLE}_6$ was found to have a special property called {\it locality} by Lawler-Schramm-Werner \cite{LSW1, LSW2}.
 Later, it was proved by S. Smirnov that ${\rm SLE}_6$ is the scaling limit of the critical percolation exploration process on two-dimensional triangular lattice.
 In honor of Schramm, SLE is now also called Schramm-Loewner evolution.

In this paper, we extend the SLE theory from the upper half plane $\HH$ to
a standard slit domain
--a specific canonical 
multiply connected planar domain.
Based on recent results of Chen-Fukushima-Rhode  \cite{CFR}  on
chordal Komatu-Loewner (KL) equation and following
the lines briefly laid 
by R. O. Bauer and R. M. Friedrich  
\cite{BF1, BF2, BF3}, 
we show that, for a corresponding evolution in
a standard slit domain $D=\HH\setminus \bigcup_{k=1}^N C_k$, the possible candidates of the driving processes are given by the solution $(\xi(t), \s(t))$ of a special Markov type stochastic differential equation
whose diffusion and drift coefficients are homogeneous function of degree $0$ and $-1$, respectively. Here
$\xi(t)$ is a motion on $\partial\HH$ and $\s(t)$ is a motion of slits $C_k,\;1\le k\le N$.
When no slit is present, $\xi(t)$ becomes $\sqrt{\kappa}B_t$ as in the simply connected domain $\HH$ case. The solution $(\xi(t), \s(t))$ of the SDE then produces a family of random conformal mappings from the multiply connected domains $D\setminus F_t $ to the canonical slit domains  $D(\s (t))$ via KL equations.
This family or its associated increasing family of
random growing $\HH$-hulls
$\{F_t; t\geq 0\}$
is called a  {\it stochastic Komatu-Loewner evolution} (SKLE in abbreviation).
We then study the locality property of SKLE.

We now
recall the setting formulated in \cite{CFR} and some
of its results   that will be utilized in
this paper. They are followed by a detailed description 
of the rest of the paper.

A domain of the form $D=\HH\setminus \bigcup_{k=1}^N C_k$ is called a {\it standard slit domain} where $\{C_k\}$ are mutually disjoint line segments
in $\HH$ parallel to the $x$-axis.
Denote by $\dd$ the collection of 'labeled (ordered)' standard slit domains.  For instance,
$\HH\setminus \{C_1, C_2, C_3,\cdots, C_N\}$ and
$\HH\setminus \{C_2, C_1, C_3, \cdots, C_N\}$ are considered as different
elements of $\dd$ in general although they correspond
to the same subset $\HH\setminus \bigcup_{i=1}^N C_i$ of $\HH.$
For $D$ and  $\wt D$ in $\dd$,  define the distance $d(D,\wt D)$ between them by
\begin{equation}\label{e:1.1}
d(D,\wt D)=\max_{1\le i\le N} \left( |z_k-\wt z_k|+|z_k'-\wt z_k'| \right),
\end{equation}
where, $z_k$ and $z_k'$ (respectively, $\wt z_k$ and $\wt z_k'$)
are the left and right endpoints of the $k$th slit of $D$ (respectively,
$\wt D$).

\begin{figure}[h]
	\begin{center}
    \vspace{-1em}	\includegraphics[scale=.45]{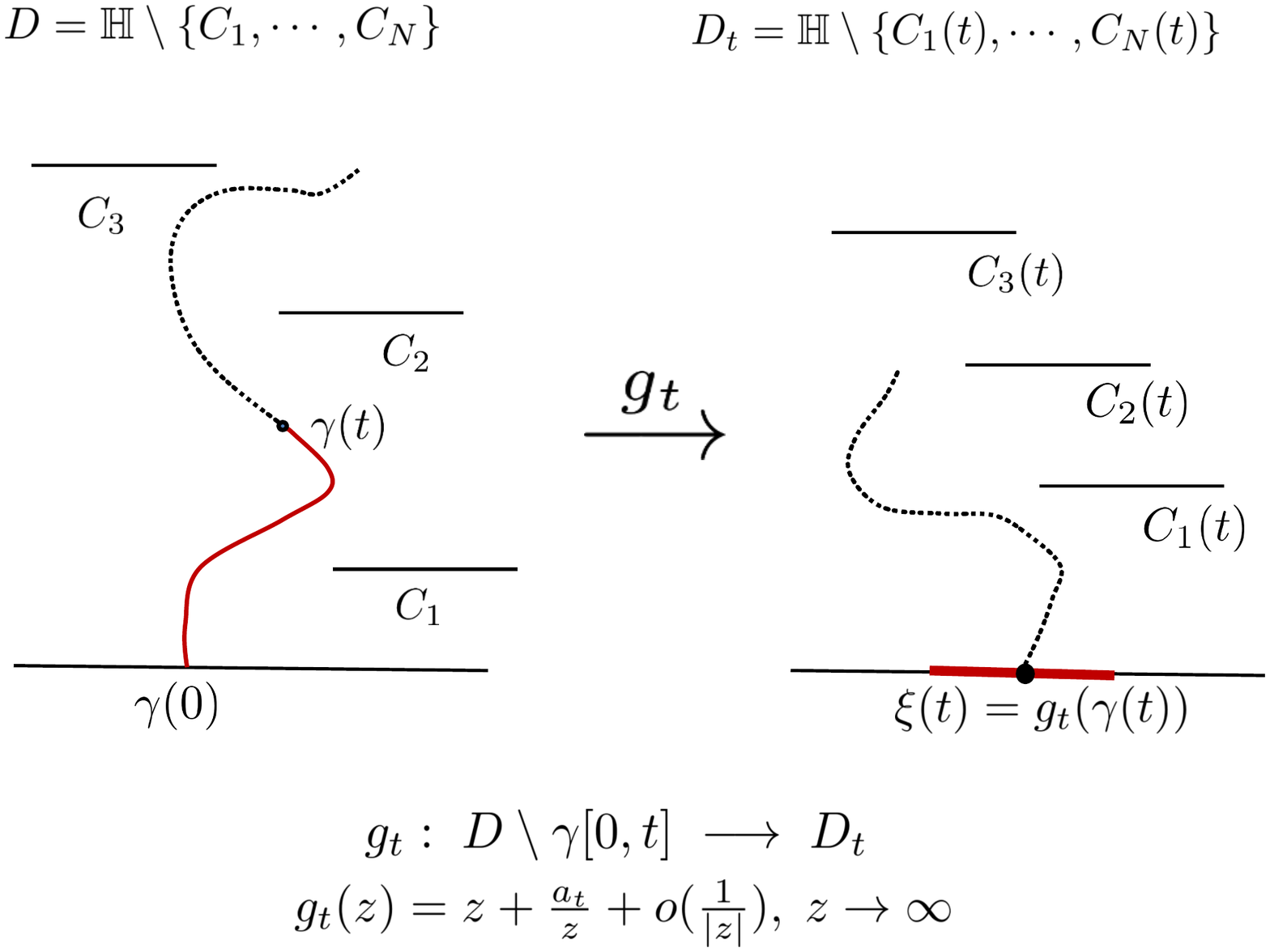}
    \vspace{-4.5em} %adjust the distance between the caption and the graph
    \caption{Conformal mapping $g_t$}\label{fig:1}
   \vspace{-1em}
	\end{center}
	\end{figure}

We fix $D\in \dd$ and consider a Jordan arc
\begin{equation}\label{I.1}
\gamma: [0,t_\gamma) \rightarrow \overline D \quad \hbox{with }\  \gamma(0)\in \partial\HH \ \hbox{ and } \gamma(0,t_\gamma)\subset D
\ \hbox{ for } 0<t_\gamma\le \infty.
\end{equation}
For each $t\in [0,t_\gamma)$, let
\begin{equation}\label{I.2}
g_t: D\setminus \gamma[0,t]\rightarrow D_t
\end{equation}
be the unique conformal map from $D\setminus \gamma[0,t]$ onto some $D_t=\HH\setminus \bigcup_{k=1}^N C_k(t)\in \dd$ satisfying a {\it hydrodynamic normalization}
\begin{equation}\label{I.3}
g_t(z)=z+\frac{a_t}{z}+o(1/|z|),\qquad z\to \infty.
\end{equation}
The coefficient $a_t$ is called the {\it half-plane capacity} of $g_t$.
The slits $C_k(t)$, $1\leq k \leq N$,
 are uniquely determined by $D$ and $\gamma [0, t]$.
 See Figure \ref{fig:1}.
Let $\s (t)$ denote the endpoints of these slits $C_k(t)$
(see \eqref{e:3.1a} below for a precise definition) and
denote $D_t$ by $D(\s (t))$.
We also define
\begin{equation}\label{I.4}
\xi(t)=g_t(\gamma(t))  \in \partial \HH ,\quad 0\le t < t_\gamma.
\end{equation}

For a Borel set $A\subset \overline \HH$,
we use $\partial_p A$ to denote the boundary of $A$
with respect to the topology induced by the path distance
in $\HH\setminus A$.
For instance, when $A\subset \HH$ is a horizontal line segment,
then $\partial _p A$ consists of the upper part $A^+$
and the lower part $A^-$ of the line segment $A$.

In \cite[Section 8]{CFR}, 
the following
results
are established:
\begin{description}
\item[(P.1)] For every $0<s< t_\gamma$,
$g_t(z)$ is jointly continuous in
$(t, z)\in [0, s]\times
(( D\cup \partial_p K\cup \partial \HH)
\setminus \gamma [0, s])$, where $K=\bigcup_{k=1}^N C_k.$
\item[(P.2)]\ $a_t$ is strictly increasing and continuous in $t\in [0,t_\gamma)$ with $a_0=0$ so that the arc $\gamma$ can be reparametrized
    in such a way
    that $a_t=2t,\ 0\le t < t_\gamma,$ which is called the {\it half-plane capacity parametrization}.

\item[(P.3)] $\xi(t)\in \partial \HH$ is continuous in $t\in [0,t_\gamma).$

\item[(P.4)] $D_t\in \dd$ is continuous in $t\in [0,t_\gamma)$ with respect to
 the metric \eqref{e:1.1} on $\dd$.
\end{description}

Historically $g_t(z)$ has been obtained by solving the extremal problem to maximize the coefficient $a_t$ among all univalent functions on $D\setminus \gamma[0,t]$ with the hydrodynamic normalization at infinity.
But, in order to prove the above continuity properties, we
used the following  probabilistic representation of $g_t(z)$ given in
\cite[\S 7]{CFR}:

\medskip
\noindent
Let $Z^{\HH,*}=\{Z_t^{\HH,*}, \ \P_z^{\HH,*}, \ z\in D^*\}$ be the {\it Brownian motion with darning} ({\it BMD})
on $D^*:=D\cup \{ c_1^*, \dots, c_N^*\}$ obtained from
the absorbing Brownian motion in $\HH$ by rendering (or shorting)
each slit $C_k$ into one single point $c_k^*$.
That is, $Z^{\HH,*}$ is an $m$-symmetric diffusion process on $D^*$ whose
 subprocess
killed upon leaving $D$ is the absorbing Brownian motion in $D$. Here $m$ is the measure on
$D^*$ that does not charge on $\{c^*_1,  \dots, c_N^*\}$ and its restriction to $D$ is the Lebesgue measure
in $D$. BMD $Z^{\HH,*}$ is unique in law 
and spends zero Lebesgue amount of time on  $\{c^*_1,  \dots, c_N^*\}$;
see  \cite{CFR} for details.
Set
$F_t:=\gamma[0,t]$, $\Gamma_r:=\{z=x+iy: y=r\},\ r>0.$
For a set $A\subset D^*$,  define $\sigma_{A}= \inf\{t>0:  Z_t^{\HH,*} \in A\}$. Then (cf.
\cite[Theorem 7.2]{CFR}
\begin{equation}\label{I.0}
\Im\; g_t(z)=\lim_{r\to\infty} r\cdot \P_z^{\HH,*}(\sigma_{\Gamma_r}<\sigma_{F_t}).
\end{equation}
Here $\Im g_t (z)$ stands for the imaginary part of the conformal map $g_t(z)$. The above formula
was first obtained by Lawler \cite{L2}
with {\it Excursion reflected Brownian motion} ({\it ERBM})
formulated there in place of BMD.
See \cite[Remark 2.2]{CFR}, \cite[\S 6]{CF2} and Remark 6.13 of the present paper  for the relationship between these two processes.

It is proved in \cite[Theorem 9.9]{CFR} that the family
of conformal maps $\{g_t (z); t\geq 0\}$
satisfies the {\it Komatu-Loewner (KL) equation} under the half-plane capacity parametrization of $\gamma$:
\begin{equation}\label{I.5}
\partial_t g_t (z)
 = -2\pi \Psi_t(g_t(z), \xi(t)) \quad \hbox{for }
0\le t< t_\gamma \  \hbox{ with } \
 g_0(z)=z\in (D\cup \partial_p K)\setminus \gamma[0,t_\gamma),
\end{equation}
where $\Psi_t(z,\xi)$,  $z\in D_t$,  $\xi\in \partial\HH$,
 is the {\it BMD-complex Poisson kernel} for $D_t$, namely,
the unique analytic function in $z\in D_t$ vanishing at $\infty$ whose imaginary part is the Poisson kernel of the BMD for the standard slit domain $D_t$
with pole $\xi \in \partial \HH$. Here $\partial_t := \frac{\partial}{\partial t}$ stands for the partial derivative in $t$.

The ODE \eqref{I.5} was derived in \cite{BF3} as well as in its original
form by Y. Komatu in \cite{K},
but only in the sense of left-derivative
 in $t$ on its left hand side.
 In \cite[\S 9]{CFR}, a Lipschitz continuity of the complex Poisson kernel $\Psi(z,\xi)$ of the BMD for $D\in \dd$ under the perturbation of $D\in \dd$ is
 established, which together with {\bf(P.4)} yields
the following property by taking $K(t)=\bigcup_{k=1}^N C_k(t)$:

\medskip
\noindent
{\bf(P.5)}  $\Psi_t(z,\xi)$ is jointly continuous in $ (t,z,\xi)\in \bigcup_{t\in [0,t_\gamma)}\{t\}\times
\left( D_t\cup\partial_pK(t)\cup(\partial\HH\setminus \{\xi\})\right)$.

\medskip
\noindent
{\bf (P.1), (P.3), (P.5)} imply that the righthand side of the equation \eqref{I.5} is continuous in $t$ and consequently it becomes a genuine ODE.

The rest of this paper is organized as follows.
In Section \ref{S:2}, we show under the above mentioned setting of \cite{CFR} that the
endpoints $\s(t)$
of the slits $C_j(t),\ 1\le j\le N,$ satisfy an
ODE analogous to the KL equation,
in terms of the boundary trace of the BMD-complex Poisson kernel $\Psi_t(z,\xi).$

In Section \ref{S:3}, we introduce a probability measure on a collection of Jordan arcs $\gamma$ using the half-plane capacity parametrization that satisfies a domain Markov property and an invariance property
under linear conformal map.
We then study the basic properties of the induced process $\W_t=(\xi(t), \s(t))$.
In particular, under certain regularity assumption (conditions {\bf (C.1)} and {\bf (C.2)} in subsection \ref{S:3.5}),
$\W_t$ is shown to satisfy
an SDE
for which the diffusion and drift coefficients for $\xi(t)$ are
homogeneous functions $\alpha(\s)$ and $b(\s)$ of degree $0$ and $-1$, respectively,
 and the endpoints $\s (t)$ of the slit motion component satisfy
the KL equation.

Conversely, given locally
Lipschitz continuous homogeneous functions $\alpha$ and $b$
of degree 0 and -1, respectively,
 we establish in Section \ref{S:4} that the corresponding SDE for $\W_t$ has a unique strong solution.  We show that the solution $ (\xi(t), \s(t))$ has
 a Brownian scaling property and is homogeneous in $x$-direction.

The solution $(\xi(t), \s(t))$ obtained above 
generates 
 a family of (random) conformal mappings $\{g_t(z)\}$ via the Komatu-Loewner equation \eqref{e:5.25}.   
The associated random growing hulls $\{F_t\}$ is called the SKLE drivn by $(\xi(t), \s(t))$ determined by coefficients $\alpha, b$ and is denoted   
by ${\rm SKLE}_{\alpha, b}.$
Its basic properties including pathwise properties as a solution of an ODE
as well as a certain scaling property and a domain Markov property
of its distribution
are studied in Section \ref{S:5}.
The induced random measures on $\{F_t; t\geq 0\}$ are shown to have the
domain Markov property and a dilation and translation invariance property.

In Section \ref{S:6},
we introduce a constant $b_{\rm BMD}$
measuring a discrepancy of a standard slit domain from $\HH$ relative to the BMD.   We call this constant the {\it BMD domain constant}.
This section concerns the locality of ${\rm SKLE}_{\alpha,b}$-hulls $\{F_t\}$, 
which is a property that $\{\Phi_A (F_t) \}$, after a suitable time change,
has the same distribution as $\{F_t\}$ for any hull $A\subset D\in \dd $ and 
 its associated canonical map 
$\Phi_A: D\setminus A\mapsto  \wt D \in  \dd$. 
We do not know if $\{F_t\}$ is generated even by a continuous curve.  Nevertheless,   
a generalized Komatu-Loewner equation \eqref{e:6.25} for the map $\wt g_s(z)$ 
associated with the image hulls  
 $\{\Phi_A(F_t)\}$ can be derived by first establishing 
the joint continuity of $\wt g_t(z)$ using BMD and the absorbing Brownian motion. This equation and a generalized It\^o formula will lead us to Theorem  \ref{T:6.12}
 asserting that 
 ${\rm SKLE}_{\alpha,-b_{\rm BMD}}$ for a constant 
$\alpha >0$ enjoys a locality property if and only if $\alpha = \sqrt{6}$.

To establish the equation \eqref{e:6.25},  
 we need a comparison of half-plane capacities
obtained by S. Drenning \cite{D} using ERBM.
A full proof of this comparison theorem
using BMD instead of ERBM will be supplied 
in the Appendix, Section \ref{S:7}, of this paper. 

An SKLE is produced by a pair $(\xi(t),\s(t))$ of a motion $\xi(t)$ on $\partial\HH$ and a motion $\s(t)$ of slits via Komatu-Loewner equation, while an SLE is produced by a motion on $\partial\HH$ alone via Loewner equation. 
They are subject to different mechanisms.
Nevertheless. as a family of random growing hulls,  it is demonstrated
 in \cite{CFS} that,  when $\alpha$ is constant,   ${\rm SKLE}_{\alpha,b}$ is, up to
some random hitting time and modulo a time change, equivalent in distribution to the chordal ${\rm SLE}_{\alpha^2}$.
Moreover, it is shown in \cite{CFS} that, after a reparametrization in time,
 ${\rm SKLE}_{ \sqrt{6},-b_{\rm BMD}}$ has 
the same distribution as chordal ${\rm SLE}_6$
in upper half space $\HH$. 
  In relation to Theorem \ref{T:6.12} of the present paper, the locality of 
${\rm SLE}_6$ will be revisited and examined in \cite{CFS}. 

The present paper only treats chordal SKLEs.  The study of K-L equations and SKLEs for other canonical multiply connected planar domains as annulus, circularly slit annulus and circularly slit disk will be recalled and examined in \cite{CFS}.

Throughout this paper, we use ``:=" as a way of definition.
For $a, b\in \RR$, $a\vee b:=\max \{a, b\}$ and $a\wedge b := \min \{a, b\}$.

\section{Komatu-Loewner equation for slits}\label{S:2}

We keep the setting and the notations of \cite{CFR} that are described in the preceding section.

In this and the next sections, we consider simple curves only.  
We use them to find out what kind of driving processes should be
 for general stochastic Komatu-Loewner equation.
 We parameterize the Jordan arc $\gamma$ by its half-plane capacity,
 which is always possible in view of  {\bf (P.2)}.
  For $t\in [0,t_\gamma)$, the conformal map $g_t$ from $D\setminus \gamma[0,t]$ onto $D_t$ can be extended analytically to $\partial_p K$ in the following manner.

Fix $j\in \{1, \dots,  N\}$, and denote the left and right endpoints of $C_j$ by $z_j=a+ic$ and $ z_j^r=b+ic$,  respectively.
Denote the open slit $C_j\setminus \{z_j,\;z_j^r\}$ by  $C_j^0$.
  Consider the open rectangles
\[R_+ :=\{z: a<x<b,\; c<y<c+\delta\},\qquad  R_- :=\{z: a<x<b,\; c-\delta<y<c\},\]
and $ R:=R_+\cup C_j^0\cup R_-$,
where $\delta>0$ is sufficiently small so that
$R_+\cup R_-\subset D\setminus \gamma[0,t_\gamma).$  Since $\Im g_t(z)$ takes a constant value at  the slit $C_j$,
 $g_t$  can be extended to an analytic function $g_t^+$ (resp. $g_t^-$) from $R_+$ (resp. $R_-$) to $R$ across $C_j^0$
 by the Schwarz reflection.

We next take $\eps>0$ with $\eps<(b-a)/2$ so that $B(z_j,\eps)\setminus C_j\subset D\setminus \gamma[0,t_\gamma].$ Then $\psi(z)=(z-z_j)^{1/2}$ maps $B(z_j,\eps)\setminus C_j$ conformally onto $B(0,\sqrt{\eps})\cap \HH.$  As in the proof of \cite[Theorem 7.4]{CFR},
$$
f_t^\ell(z)=g_t\circ \psi^{-1}(z)=g_t(z^2+z_j)
$$
 can be extended
analytically to $B(0,\sqrt{\eps})$
by the Schwarz reflection and by noting that the origin $0$ is a removable singularity for $f_t^\ell.$
Similarly,
we can induce an analytic function $f_t^r$ on $B(0,\sqrt{\eps})$ from $g_t$ on
 $B(z_j^r,\eps)\setminus C_j.$

  For an analytic function $u(z),$ its derivatives in $z$ will be denoted by $u'(z)$, $ u''(z)$ and so on.

\begin{lem}\label{L:2.1}\begin{description}
\item {\rm(i)}\ $\partial_t   g_t^\pm  (z)$, $(g_t^\pm)'(z)$ and
$(g_t^\pm)''(z)$ are continuous in $(t,z)\in [0,t_\gamma)\times R$.

\item{\rm (ii)}\ $\eta_t(z,\zeta):=\Psi_t(g_t(z),\zeta)$ can be extended to an analytic function $\eta_t^+(z,\zeta)$ {\rm(}resp. $\eta_t^-(z,\zeta)${\rm)} from $R_+$ {\rm(}resp. $R_-${\rm)} to $R$
by the Schwarz reflection, and
\begin{equation}\label{e:2.1}
( \eta_t^\pm(z,\zeta) )'\ \hbox{\rm are continuous in}\ (t,z,\zeta) \in [0,t_\gamma)\times R\times \partial\HH.
\end{equation}

\item{\rm(iii)}\ $(g_t^\pm)'(z)$ are differentiable in $t\in (0,t_\gamma)$
and
\begin{equation}\label{e:2.2}
\partial_t (g_t^\pm)'(z)
\ \hbox{\rm are continuous in}\ (t,z)\in [0,t_\gamma)\times R.
\end{equation}

\item{\rm (iv)}\ $\partial_t f_t^\ell (z)$,
$(f_t^\ell)'(z)$ and $(f_t^\ell)''(z)$ are continuous
in $(t,z)\in [0,t_\gamma]\times B(0,\sqrt{\eps}).$

\item{\rm (v)} \
 $\wt \eta_t(z,\zeta) :=
    \Psi_t(f_t^\ell(z),\zeta)=\Psi_t(g_t(\psi^{-1}(z)),\zeta)$ can be extended to an analytic function from
$B({\bf 0},\sqrt{\eps})\cap\HH$ to $B({\bf 0},\sqrt{\eps})$ and
\begin{equation}\label{e:2.3}
( \wt \eta_t(z,\zeta))'\ \hbox{\rm is continuous in}\ (t,z,\zeta)\in [0,t_\gamma)\times B({\bf 0},\sqrt{\eps})\times \partial \HH.
\end{equation}

\item{\rm (vi)}\ $(f_t^\ell)'(z)$ is differentiable in $t\in (0,t_\gamma)$ and
\begin{equation}\label{e:2.4}
\partial_t (f_t^\ell)'(z)
\ \hbox{\rm is continuous in}\ (t,z)\in [0,t_\gamma)\times B({\bf 0},\sqrt{\eps}).
\end{equation}

\item {\rm (vii)}\ The statements {\rm (iv),\;(v),\;(vi)} in the above remain valid with $f_t^r$ in place of $f_t^\ell.$
\end{description}
\end{lem}
\pf\ (i)\ This follows from the Cauchy integral formulae of derivatives of
$g_t^\pm$ combined with the property {\bf(P.1)}
and \eqref{I.5}.

\medskip

(ii)\ This can be proved in the same way as (i) using
{\bf (P.1)} and {\bf (P.5)}.

\begin{figure}[h]
	\begin{center}
    \vspace{-1em}
	\includegraphics[scale=.45]{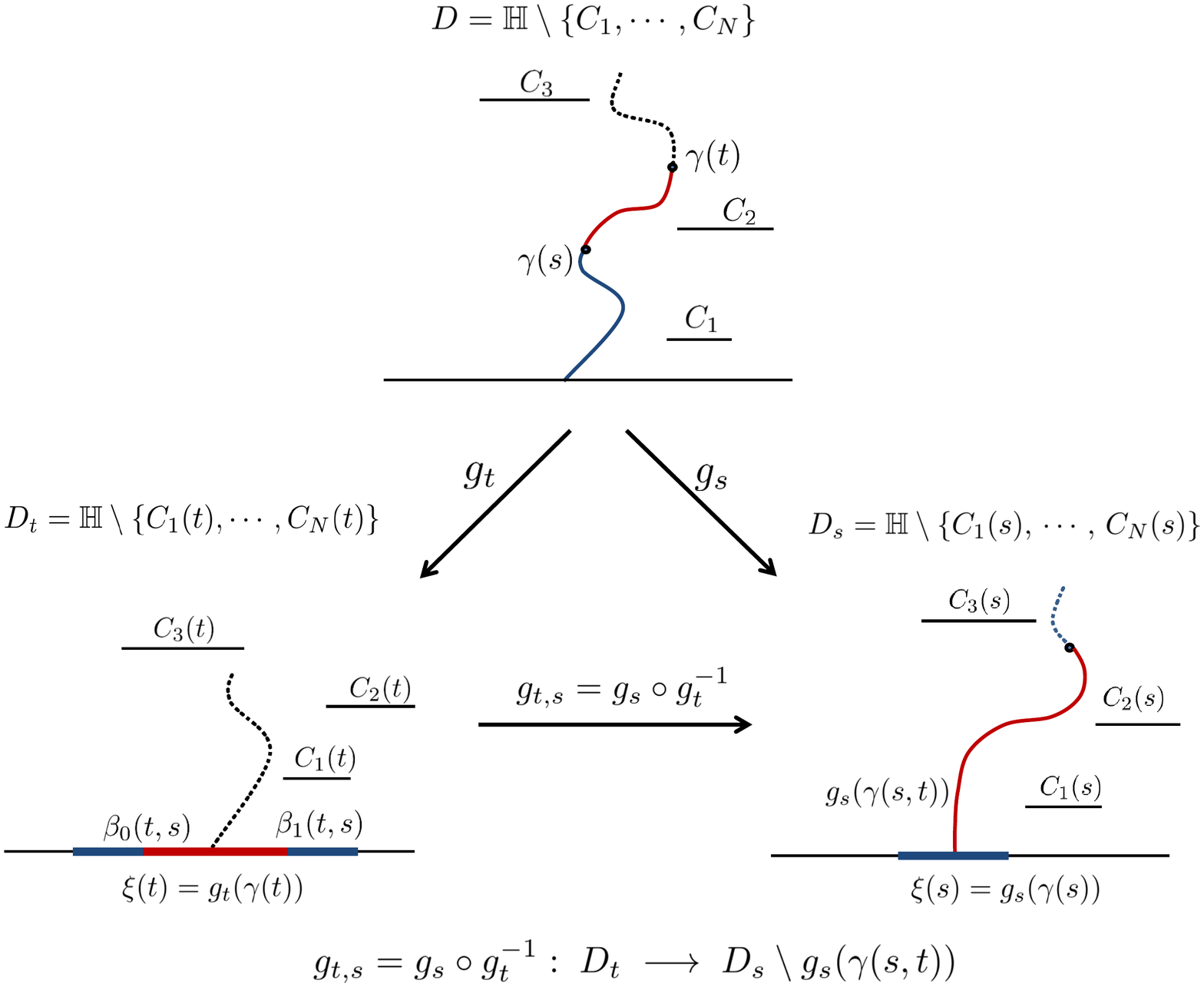}
    \vspace{-1em} %adjust the distance between the caption and the graph
    \caption{Conformal mapping $g_{t,s}$}\label{fig:2}
    \vspace{-1em}
	\end{center}
	\end{figure}

\medskip
(iii)\ For $0<s<t<t_\gamma,$ define
$ g_{t,s}=g_{s}\circ g_t^{-1}$,
which is a conformal map from $D_t$ onto $D_{s}\setminus g_{s}(\gamma[s,t])$.
Define
$
\xi(t)=g_t(\gamma(t)) =\lim_{z\to\gamma(t), z\in D\setminus \gamma[0,t]}g_t(z).
$
It is easy to see that
$\xi(t)\in \partial\HH.$  Furthermore, there exist unique points $\beta_0(t,s)<\beta_1(t,s)$ from $\partial\HH$ such that
$
\beta_0(t,s)<\xi(t)<\beta_1(t,s),\quad
g_{t,s}(\beta_0(t,s))=g_{t,s}(\beta_1(t,s))=\xi(s),
$
and
$$
\Im g_{t,s}(x+i0+)
\begin{cases}
= 0 \quad &\hbox{for } x\in \partial\HH\setminus(\beta_0(t,s),\beta_1(t,s)),\\
> 0 &\hbox{for } x\in (\beta_0(t,s),\beta_1(t,s)).
\end{cases}
$$
See Figure \ref{fig:2}.
We know from \cite[(6.22)]{CFR} that
\begin{equation}\label{e:2.5a}
g_{s}(z)-g_t(z)=\int_{\beta_0(t,s)}^{\beta_1(t,s)} \Psi_t(g_t(z),x)\Im g_{t,s}(x)dx.
\end{equation}
Taking derivative in $z$ yields
\[
(g_s^\pm)'(z)-(g_t^\pm)'(z)=\int_{\beta_0(t,s)}^{\beta_1(t,s)}
( \eta^\pm_t(z,x) )' \, \Im g_{t,s}(x)dx.
\]
On the other hand, it is established in \cite[Lemma 6.2]{CFR} that
 for $0\le s<t< t_\gamma$ that
\begin{equation}\label{e:2.6a}
2(s-t)=a_t-a_{s}=\frac{1}{\pi}\int_{\beta_0(t,s)}^{\beta_1(t,s)} \Im g_{t,s}(x+i0+)dx.
\end{equation}
Taking quotient of the last two displays and then passing $s\uparrow t$ yields
\begin{equation}\label{e:2.5}
\partial_t^- (g_t^\pm)'(z)
=-2\pi (\eta^\pm_t(z, \xi(t)))' ,
\end{equation}
where $\partial_t^-$ denotes the left-derivative in $t$.
In view of \eqref{e:2.1} and the property {\bf(P.3)},
the right hand side of \eqref{e:2.5} is continuous in
$t\in [0, t_\gamma )$.
Thus,  as in the proof of \cite[Thoerem 9.9]{CFR}, $(g_t^\pm)'(z)$ is differentiable in $t$. Consequently, \eqref{e:2.2} follows in view of  \eqref{e:2.5}.

\medskip
(iv)\ Let $\wt\gamma$ be a closed
smooth Jordan curve in $B({\bf 0},\sqrt{\eps}).$  By Cauchy's integral formula
\[(f_t^\ell)'(z)=\frac{1}{2\pi i}\int_\eta \frac{f_t^\ell(\zeta)}{(\zeta-z)^2}d\zeta,\quad z\in {\rm ins}\ \wt\gamma.\]
Since $f_t^\ell(\zeta)=g_t(\zeta^2+z_j)$ is continuous in $t$ uniformly in $\zeta\in \wt\gamma$ by {\bf(P.1)}, we get the desired continuity.  The same is true for $(f^\ell_t)''(z).$

\medskip

(v)\ Since $\Im \wt \eta_t(z,\zeta)$ is constant in $z$ on
$B({\bf 0},\sqrt{\eps})\cap\partial\HH\setminus \{{\bf 0}\}$, it extends
analytically to $B({\bf 0},\sqrt{\eps})\setminus \{{\bf 0}\}$ by the
Schwarz reflection.
Note that ${\bf 0}$ is a removable singularity because
$\Im \eta_t(z,\zeta)$ is bounded near $\{{\bf 0}\}$.
The second assertion can be shown as the proof of (ii) using {\bf (P.1)} and {\bf(P.5)}.

\medskip
(vi)\ Taking $z$ to be $\psi^{-1}(z)$ in \eqref{e:2.5a},
we have
$$
f_s^\ell (z)- f_t^\ell (z)=\int_{\beta_0(t,s)}^{\beta_1(t,s)}
 \wt \eta_t(z,x)
 \Im g_{t,s}(x)dx \quad \hbox{for }  z\in B({\bf 0},\sqrt{\eps})
  \hbox{ and }  s<t.
$$
Differentiating in $z$ gives
\[
(f_s^\ell)'(z)-(f_t^\ell)'(z)=\int_{\beta_0(t,s)}^{\beta_1(t,s)}
 ( \wt \eta_t(z,x))'\Im g_{t,s}(x)dx
 \quad \hbox{for }  z\in B({\bf 0},\sqrt{\eps})
  \hbox{ and }  s<t.
\]
Taking quotient of the above with \eqref{e:2.6a} and   passing
 $s\uparrow t$ yields
$
\partial_t^- (f_t^\ell)'(z) = 2\pi ( \wt \eta_t(z,\xi(t)))'.
$
Since the right hand side is continuous in $t$ by \eqref{e:2.3} and {\bf (P.3)}, we arrive at the conclusion (vi).
\qed

Denote by $z_j(t)$ and $ z^r_j(t)$ the left and right endpoints of the slit $C_j(t)$,
where $1\le j\le N$,  for $D_t\in \dd$ and $t\in [0,t_\gamma)$.
Since $g_t$ is a homeomorphism between $\partial_pC_j$ and $\partial_pC_j(t)$,
for each $t\in [0,t_\gamma)$ and  $1\le j\le N$,
there exist unique
$$
\wt z_j(t) = \wt x_j(t)+iy_j\in \partial_pC_j  \quad
\hbox{and} \quad \wt z_j^r(t) = \wt x_j^r (t)+iy_j\in \partial_pC_j
$$
so that  $g_t(\wt z_j(t))=z_j(t)$ and $g_t(\wt z_j^r (t))=z_j^r (t)$.

\begin{lem}\label{L:2.2}
\ {\rm(i)}\ If $\wt z_j(t)\in C_j^+\setminus \{z_j,z_j^r \},$ then
\begin{equation}\label{e:2.6}
(g_t^+)'(\wt z_j(t))=0,\qquad (g_t^+)''(\wt z_j(t))\neq 0.
\end{equation}
{\rm(ii)}\ If $\wt z_j(t)\in C_j^-\setminus \{z_j,z_j^r \},$ then
{\rm\eqref{e:2.6}} holds with $g_t^-$ in place of $g_t^+.$

\noindent
{\rm(iii)}\ If $\wt z_j(t)\in \partial_p C_j\cap B(z_j,\eps),$ then,
for $\psi(z)=(z-z_j)^{1/2},$
\begin{equation}\label{e:2.7}
(f_t^\ell)'(\psi(\wt z_j(t)))=0,\qquad (f_t^\ell)''(\psi(\wt z_j(t)))\neq 0.
\end{equation}

\noindent
{\rm(iv)}\ If $\wt z_j(t)\in \partial_p C_j\cap B(z_j^r ,\eps),$ then,
for $\psi(z)=(z-z_j^r )^{1/2},$
\begin{equation}\label{e:2.8}
(f_t^r)'(\psi(\wt z_j(t)))=0,\qquad (f_t^r)''(\psi(\wt z_j(t)))\neq 0.
\end{equation}

\noindent
{\rm (v)}\ The above four statements also hold for $\wt z_j^r (t)$ in place of $\wt z_j(t).$
\end{lem}

\pf It suffices to prove (i) and (iii).
 $g_t^+$ is analytic on $R$ and $\wt z_j(t)\in R.$
 Suppose $g_t^+(z)-z_j(t)$ has a zero of order $m$ at $\wt z_j(t)$: for some analytic function $h$ with $h(\wt z_j(t))\neq 0,$
\[g_t^+(z)-z_j(t)=g_t^+(z)-g_t^+(\wt z_j(t))=(z-\wt z_j(t))^m h(z).\]
Then, in view of \cite[p131,\;Th.11]{A}, there exists $\eps_0>0$ with $B(\wt z_j(t), \eps_0)\subset R$ and $\delta_0>0,$ such that, for any $w\in B(z_j(t),\delta_0)$, $(g_t^+)^{-1}(w)\cap B(\wt z_j(t),\eps_0)$ consists of $m$ distinct points.
Since $g_t$ is homeomorphic between $\partial_pC_j$ and $\partial_pC_j(t),$ there exists $\delta_{00}>0$ such that, for any $\delta\in (0, \delta_{00})$ and for any $w\in B(z_j(t),\delta)\cap C_j(t)$ with $w\neq z_j(t),$ $(g_t^+)^{-1}(w)\subset B(\wt z_j(t),\eps_0)\cap C_j^+$ consists of two points because $z_j(t)$ is an endpoint of $C_j(t)$ and so $w$ corresponds to two distinct points of $\partial_pC_j(t).$  Hence $m=2.$

\noindent
(iii)\ Except for the last part, the following proof is similar to that of (i).

$f_t^\ell$ is analytic on
$B({\bf 0},\sqrt{\eps})$
and $\psi(\wt z_j(t))\in B({\bf 0},\sqrt{\eps}).$
 Suppose $f_t^\ell(z)-z_j(t)$ has a zero of order $m$ at $\psi(\wt z_j(t))$: for some analytic function $h$ with $h(\psi(\wt z_j(t)))\neq 0,$
\[f_t^\ell(z)-z_j(t)=f_t^\ell(z)-f_t^\ell(\psi(\wt z_j(t)))=(z-\psi(\wt z_j(t)))^m h(z), \ z \in B({\bf 0},\sqrt{\eps}).
\]
Then, as in the proof of (i), there exists $\eps_0>0$ with $B(\psi(\wt z_j(t)), \eps_0)\subset B({\bf 0},\sqrt{\eps})$
 and $\delta_0>0,$ such that, for any $w\in B(z_j(t),\delta_0)$, $(f_t^\ell)^{-1}(w)\cap B(\psi(\wt z_j(t)),\eps_0)$ consists of $m$ distinct points.
Since $z_j(t)$ is the endpoint of $C_j(t)$ and $g_t$ is homeomorphic between $\partial_pC_j$ and $\partial_pC_j(t),$ there exists $\delta_{00}>0$ such that, for any $\delta\in (0, \delta_{00})$ and for any $w\in B(z_j(t),\delta)\cap C_j(t)$ with $w\neq z_j(t),$ $(g_t)^{-1}(w)\subset \partial_pC_j\cap B(z_j,\eps)$ consists of two points.
In fact, $w$ corresponds to two distinct points $w_+\in C_j^+(t),\; w_-\in C_j^-(t)$ so that $g_t^{-1}(w)=\{\wt w_+,\wt w_-\}$ with $g_t(\wt w_\pm)=w_\pm.$
Then $(f_t^\ell)^{-1}(w)=\psi(g_t^{-1}(w))=\{\psi(\wt w_+),\psi(\wt w_-)\}$ consists of two distinct points of $B({\bf 0},\sqrt{\eps}).$  Therefore $m=2.$
 \qed

We let $h(t,z)=(g_t^+)'(z).$  Then $h(t,z)$ is a $C^1$-function in $(t,z)\in (0,t_\gamma)\times R$ by virtue of Lemma \ref{L:2.1}.

Assume that $\wt z_j(t_0)\in C_j^+\setminus \{z_j, z_j^r \}.$
By \eqref{e:2.6},
\begin{equation}\label{e:2.9}
h(t, \wt z_j(t) )=0 \qquad \hbox{for } t\in (t_0-\delta_1, t_0+\delta_1)
\end{equation}
for some $\delta_1>0$.
On the other hand, $|h'(t, z)|=|(g_t^+)''(z)| >0$  by Lemma \ref{L:2.2}.
So by the implicit function theorem , there is some
 $\delta_2\in (0, \delta_1)$ so that $t\mapsto \wt z_j(t)$
is  $C^1$  in $t\in (t_0-\delta_2, t_0+\delta_2)$.
 Differentiating \eqref{e:2.9} in $t$ yields
\begin{equation}\label{e:2.10}
\frac{d}{dt} \wt z_j(t)= -\frac{\partial_t h(t,z)}{h'(t,z)}\Big|_{z=\wt z_j(t)}
\qquad \hbox{for } t\in (t_0-\delta_2, t_0+\delta_2) .
 \end{equation}

The same assertions hold for $\wt z_j(t)$ when $\wt z_j(t)\in C_j^-\setminus \{z_j, z_j^r \}$.   A similar argument shows, by using (iii) and (iv) of Lemma \ref{L:2.2}, that $\psi(\wt z_j(t))$ is a $C^1$ function of $t$ in a neighborhood of $t_0$
 when $\wt z_j(t_0)\in \partial_pC_j\cap\left(B(z_j,\eps)\cup B(z_j^r ,\eps)\right)$.

\begin{thm}\label{T:2.3}\ The endpoints $z_j(t)=x_j(t)+iy_j(t),\ z_j^r (t)=x_j^r (t)+iy_j(t),$ of $C_j(t)$ satisfy the following equations for $1\le j\le N$:
\begin{equation}\label{e:2.11}
\frac{d}{dt} y_j(t)=-2\pi \Im \Psi_t(z_j(t), \xi(t)),
\end{equation}
\begin{equation}\label{e:2.12}
\frac{d}{dt} x_j(t)=-2\pi \Re \Psi_t(z_j(t), \xi(t)),
\end{equation}
\begin{equation}\label{e:2.13}
\frac{d}{dt} x_j^r (t)=-2\pi \Re \Psi_t(z_j^r (t), \xi(t)) .
\end{equation}
\end{thm}
\pf
It suffices to prove \eqref{e:2.11}-\eqref{e:2.12}.
 It follows from \eqref{I.5} and (i), (iv) of Lemma \ref{L:2.1} that
\begin{equation}\label{e:2.16}
 \partial_t g^\pm_t (z)  =  -2\pi \Psi_t ( g^\pm_t (z), \xi (t))
 \quad \hbox{ for } z\in \partial_p C_j \setminus \{z_j, z_j^r \}
\end{equation}
and
\begin{equation}\label{e:2.1y}
 \partial_t f^\ell_t (z)  =  -2\pi \Psi_t ( f^\ell_t (z), \xi (t))
\quad \hbox{ for }  z\in \partial_p C_j \cap B(z_j, \eps).
\end{equation}
Note that $z_j (t)= g^\pm_t (\wt z_j(t))$ when
 $\wt z_j (t)\in \partial_p C_j \setminus \{z_j, z_j^r \}$ and
$z_j (t)= f^\ell_t (\wt z_j(t))$
when $ \wt z_j (t) \in \partial_p C_j \cap B(z_j, \eps)$.
Since $\wt z_j(t)$ is $C^1$ in $t$, we have by
 Lemma \ref{L:2.2}
$$
\frac{d}{dt} z_j(t)= \frac{d}{dt}   (g^\pm_t ( \wt z_j(t))
=  \partial_t g^\pm_t (\wt z_j(t)) + (g^\pm_t)'
(\wt z_j(t)) \frac{d}{dt} \wt z_j(t)
 = -2\pi \Psi_t (z_j (t), \xi (t))
$$
when $z_j (t)\in \partial_p C_j \setminus \{z_j, z_j^r \}$,
and
$$
\frac{d}{dt} z_j(t)= \frac{d}{dt} f^\ell_t  ( \wt z_j(t))
=  \partial_t f^\ell_t (\wt z_j(t)) + (f^\ell_t)'
(\wt z_j(t)) \frac{d}{dt} \wt z_j(t)
= -2\pi \Psi_t (z_j (t), \psi (t))
$$
when $z_j (t)\in \partial_p C_j \cap B(z_j, \eps)$.
This proves   \eqref{e:2.11}-\eqref{e:2.12}.
\qed

\begin{remark}\rm\ (i)\ The equation \eqref{e:2.11}-\eqref{e:2.13} was derived in \cite{BF3} by assuming that $\wt z_j(t)\in \partial_pC_j\setminus \{z_j, z_j^r \}$ and also by taking for granted the smoothness of
\ \lq$\frac{d}{dz}g_t(z)$\rq\ in two variables
$(t, z)$, which is now established by Lemma \ref{L:2.1}.

(ii)\ If
\begin{equation}\label{e:2.14}
g_t(z_j)= z_j(t) \ \hbox{ and } \  g_t(z_j^r )=z_j^r (t)\qquad
 \hbox{for } t\in (0, t_\gamma) \hbox{ and } 1\le j\le N,
\end{equation}
then Theorem \ref{T:2.3} is merely a special case of the Komatu-Loewner equation \eqref{I.5} with $z=z_j$ and $z=z_j^r $, $1\le j\le N$, respectively.
But in general \eqref{e:2.14} is not true.  \qed
\end{remark}

We call \eqref{e:2.11}-\eqref{e:2.13} the
{\it Komatu-Loewner equation} for the slits.

\section{Randomized curve $\gamma$ and induced process $\W$}\label{S:3}

\subsection{Random curve with domain Markov property and a conformal invariance}\label{S:3.1}

As in the previous sections, for
 a standard slit domain $D=\HH\setminus \bigcup_{k=1}^N C_k$,
the left and right endpoints of the $k$th slit $C_k$ are denoted by $z_k=x_k+iy_k$ and $z_k^r =x_k^r +iy_k^r $, respectively.
Recall that $\dd$ is the collection of all  labeled (or, ordered)
standard slits domains equipped with metric $d$ of \eqref{e:1.1}.
We define an open subset
$\cS$
of the Euclidean space $\RR^{3N}$ by
\begin{eqnarray}\label{e:3.1a}
\cS&=& \Big\{\s:=(\y,\x,\x^r )\in \RR^{3N}:
 \ \y,\;\x,\:\x^r \in \RR^N,\ \y>{\bf 0}, \ \x<\x^r ,
 \nonumber \\
&& \hskip 0.3truein {\rm either}\ x_j^r <x_k\ {\rm or}\ x_k^r <x_j\ {\rm whenever}\ y_j=y_k,\ j\neq k \Big\}.
\end{eqnarray}
The Borel $\sigma$-field on $\cS$ will be denoted as $\b (\cS)$.
The space $\dd$ can be identified with $\cS$ as a topological space.
We write $\s(D)$ (resp. $D(\s)$) the element in $\cS$ (resp. $\dd$) corresponding to $D\in \dd$ (resp. $\s\in \cS$).

A set $F\subset \CC$ is called an $\HH$-{\it hull} if $\overline F$ is compact,
 $F=\overline{F}\cap \HH$ and $\HH\setminus F$ is simply connected.
For $D\in \dd$ and an $\HH$-hull $F\subset D$, there exists a unique conformal map $g$ from $D\setminus F$ onto some $\wt D\in \dd$ satisfying the hydrodynamic normalization
$g(z)=z+\frac{a}{z}+o\left(\frac{1}{|z|}\right)$ as $z\to \infty.$
In what follows, such a map $g$ will be called a {\it canonical map from} $D\setminus F$.  The associated constant $a$ (which is real and non-negative) will be called the {\it half-plane capacity} of $g$
and can be evaluated as
\begin{equation}\label{e:3.1}
a=\lim_{z\to\infty} z(g(z)-z).
\end{equation}

Set
\[\wh\dd=\{\wh D=D\setminus F: \, D\in \dd
\hbox{ and } F\subset D \hbox{ is an } \HH \hbox{-hull} \}.
\]
 For $\wh D=D\setminus F\in \wh\dd$, let
$$
\O(\wh D)=\left\{\gamma=\{\gamma(t): 0\le t< t_\gamma\}: \hbox{\rm Jordan arc}, \
\gamma(0, t_\gamma)\subset \wh D,\; \gamma(0)
\in \partial(\HH\setminus F),\ 0< t_\gamma\le \infty\right\}.
$$
Two curves $\gamma,\;\wt \gamma\in \O(\wh D)$ are regarded
equivalent if $\wt \gamma$ can be
 obtained from $\gamma$ by a reparametrization.
Denote by  $\dot \O(\wh D)$   the equivalence classes of $\O(\wh D)$.

Given $\gamma\in \O(\wh D)$ for $\wh D=D\setminus F\in \wh\dd,$ let $g_t$ be the canonical map from $\wh D\setminus \gamma[0,t]=D\setminus (\gamma[0,t]\cup F)$  with the half-plane capacity $a_t$, $t\in [0,t_\gamma)$.
Note that $g_t=\wh g_t\circ g$,  where $g$ is the canonical map from $D\setminus F$ onto some $\wh D\in \dd$ and
$\wh g_t$
is the canonical map from $\wh D\setminus g(\gamma[0,t])$ onto some $D_t\in \dd$. It then follows from \eqref{e:3.1} that $a_t=a+\wh a_t$, where $a$ and
$\wh a_t$
are the half-plane capacity of $g$ and
 $\wh g_t$ respectively.

Since $\wh a_t=a_t-a_0$ is strictly increasing and continuous in $t\in [0,t_\gamma)$ with $\wh a_0=0$ by {\bf(P.2)}, the curve $\gamma$ can be reparametrized as $\wt\gamma(t)=\gamma(\wh a_{2t}^{-1})$ for $0\leq t < t_{\wt \gamma}:= \frac12 \wh a_{t_\gamma}$ so that the corresponding half-plane capacity becomes $a_0+2t$.
The curve $\wt \gamma$ is called the {\it half-plane capacity renormalization} of $\gamma.$

Throughout the rest of this paper, each
$\dot\gamma\in \dot \O(\wh D)$ will be represented by a curve (denoted by $\dot \gamma$ again) belonging to this class parametrized by the half-plane capacity. We conventionally adjoin an extra point $\Delta$ to $\overline\HH$ and define $\dot\gamma(t)=\Delta$ for $t\ge t_{\dot\gamma}$ so that $\dot\gamma$ can be regarded as a map from $[0,\infty]$ to $\overline\HH\cup\{\Delta\}.$
 We then introduce
a filtration $\{\dot\gg_t(\wh D); t\geq 0\}$  on $\dot \O(\wh D)$ by
\[
\dot\gg_t(\wh D):=\left(\sigma\{\dot\gamma(s):\;0\le s\le t\}\right)\cap \{t<t_{\dot\gamma}\}, \qquad  \dot\gg(\wh D):=\sigma\{\dot\gamma(s): s\ge 0\}.
\]

For each $D\in  \dd$, we consider a family of probability measures 
$\{\P_{D,z}; \, z\in \partial \HH \}$
on $(\dot \O( D), \dot\gg( D))$ that satisfies the property
\begin{equation}\label{e:3.2}
\P_{ D,z}(\{\dot\gamma(0)=z\})=1, \qquad z\in \partial\HH,
\end{equation}
as well as the following {\bf (DMP)} and
{\bf (IL)}.

\medskip

For each $\wh D\in \wh \dd$ and $t\ge 0,$ define the shift operator
 $\dot\theta_t: \dot \O( \wh D)\cap \{t<t_{\dot\gamma}\} \mapsto
  \dot \O( \wh D\setminus \dot\gamma[0,t])$ by
\begin{equation}\label{e:3.4a}
(\dot\theta_t\dot\gamma)(s)=\dot\gamma(t+s) \qquad \hbox{for }
 s \in [0,t_{\dot\gamma}-t).
\end{equation}

\bigskip
\noindent
{\bf (DMP)} (Domain Markov property): for  any
 $D \in \dd$, $t\ge 0$ and $z\in \partial \HH$,
\begin{equation}\label{e:3.3}
\P_{D,z}\left(\dot\theta_t^{-1}\Lambda\big|\dot\gg_t(D)\right)
=\P_{g_t (D\setminus\dot\gamma[0,t]), g_t (\dot\gamma(t))}(g_t (\Lambda )) \quad
\hbox{for every }
 \Lambda\in \dot\gg( D\setminus\dot\gamma[0,t]) .
\end{equation}
Here $g_t (z)$ is the canonical map from 
$D\setminus\dot\gamma[0,t]$.
Note that $g_t (D\setminus \dot \gamma [0, t])\in \dd$ and
$g_t (\dot\gamma(t))\in \partial \HH$ is well-defined
since $g_t (z)$ can be extended continuously
to  $\partial_p (D\setminus\dot\gamma[0,t])$.

\medskip
\noindent
{\bf (IL)} (Invariance under linear conformal map):
for any $D\in \dd$ and any linear map $f$ from $D$ onto $f(D)\in \dd$,
\begin{equation}\label{e:3.4}
\P_{f(D),f(z)}=
\P_{D,z}  \circ f^{-1}
\qquad  \hbox{for every } z\in \partial \HH .
\end{equation}

\medskip

\begin{remark}\label{R:3.1} \rm
For $\wh D =D\setminus F\in \wh \dd$, let $\Phi$ be the canonical map that maps
$\wh D$ onto $\Phi (\wh D)\in \dd$.
Suppose that $\Phi$ can be extended continuously
to $\partial_p (\HH \setminus F)$.
Then for each $z\in \partial_p (\HH \setminus F)$,
one can define
$\P_{\wh D, z} = \P_{\Phi (\wh D), \Phi (z)} \circ \Phi^{-1}$.
We can therefore restate \eqref{e:3.3}  as
\begin{equation}\label{e:3.3c}
\P_{D,z}\left(\dot\theta_t^{-1}\Lambda\big|\dot\gg_t(D)\right)
=\P_{D\setminus \dot\gamma[0,t], \dot\gamma(t)}(\Lambda) \quad
\hbox{for every }
 \Lambda\in \dot\gg(D\setminus\dot\gamma[0,t]) .
\end{equation}
This explains why we call \eqref{e:3.3} the domain Markov property.
The formulation \eqref{e:3.3} avoids the technical issue whether $\Phi$ can be extended
continuously to  $\partial_p (\HH\setminus F)$
for general $\wh D=D\setminus F\in \wh \dd$.
See Proposition \ref{P:DMP} and Theorem \ref{T:DMP} in Section \ref{S:5}. \qed
\end{remark}

\subsection{Markov property of $\W$}\label{S:3.2}

For each $D\in \dd,\ \dot\gamma\in \dot \O(D)$ and $t\in [0,t_{\dot\gamma}),$ $\dot\gamma$ induces the conformal map $g_t$ from $D\setminus \dot\gamma[0,t]$ onto $D_t=g_t(D)\in \dd.$
The conformal map $g_t (z)$
can be extended to a continuous map from $D\cup\partial_pK\cup \partial_p\gamma[0,t]\cup\partial\HH$ onto $\overline \HH$.
We occasionally write $g_t$ as $g_t^D$ or
$g_{D\setminus \dot\gamma[0,t]}$ to indicate
its dependence on $\dot\gamma\in \dot \O(D)$. Note that
$g_t$ sends $\dot\gamma(t)$ to $\xi(t)\in \partial \HH$.

Let $\{\s(t)=\s(D_t), \; t\in [0,t_{\dot\gamma})\}$ be the induced slit motion
with $D_0:=D$.
We will consider the joint process
\[\W_t=\left\{
\begin{array}{ll}
(\xi(t), \s(t))\in \RR \times S \subset \RR^{3N+1},&\quad 0\le t<t_{\dot\gamma},\\
\delta, &\quad t\ge t_{\dot\gamma}.
\end{array}
\right.
\]
Here the real part of $\xi(t)\in \partial\HH$ is designated by $\xi(t)$ again and $\delta$ is an extra point conventionally adjoined to $\RR\times S.$
We shall occasionally write $\s(t)$ as $g_t^D(\s)$ with $\s=\s(D).$

To establish the Markov property of $\W_t$, we need the following measurability results.

\begin{lem}\label{L:3.1}\ Fix $D\in \dd$ and $t\ge 0.$

\begin{description}
\item{\rm(i)}\ For each $z\in D$, $\Im g_t(z)$ is a $[0,\infty)$-valued $\dot\gg_t(D)$-measurable function on $\dot\Omega(D).$

\item{\rm(ii)}\ $g_t(z)$ is an $\overline\HH$-valued
$\b(D\cup\partial_pK\cup\partial_p\dot\gamma[0,t]\cup\partial\HH)\times \dot\gg_t(D)$-measurable function on
$ ( D\cup\partial_pK\cup\partial_p\dot\gamma[0,t]\cup\partial\HH  ) \times \dot\O(D).$

\item{\rm(iii)}\ $\W_t$ is an $\RR^{3N+1}$-valued $\dot\gg_t(D)$-measurable function on $\dot \O(D).$
\end{description}
\end{lem}

\pf\ (i)\ We make use of the probabilistic representation \eqref{I.0} of $\Im g_t(z).$ \quad Take $r>0$ such that the set $\HH_r=\{z\in \HH: \Im z<r\}$ contains $\dot\gamma(0,t]\cup K.$ \quad It suffices to show $F_{z,r}(\dot\gamma)=\P_z^{\HH,*}(\sigma_{\Gamma_r}<\sigma_{\dot\gamma(0,t]})$ is a $\dot\gg_t(D)$-measurable function on $\dot\O(D)$ for each fixed $z\in D.$

Let $Z^{\HH_r,*}=(Z_t^{\HH_r,*}, \zeta, \P_z^{\HH_r,*})$ be the BMD on $\HH_r^*=
(D\cap \HH_r)
\cup \{c_1^*,\cdots, c_N^*\}$ obtained from the absorbing Brownian motion on $\HH_r$ by rendering each hole $C_k$ into a single point $c_k^*$,
with life time $\zeta$. Then
\begin{equation}\label{e:3.5}
F_{z,r}(\dot\gamma)
=\P_z^{\HH_r,*}(\sigma_{\dot\gamma(0,t]}=\infty)
=\P_z^{\HH_r,*}(\dot\gamma(0,t]\cap Z_{[0,\zeta)}^{\HH_r,*}=\emptyset).
\end{equation}

Let $\HH_r^*\cup\{\Delta\}$ be the one-point compactification of $\HH_r^*.$.
 As the sample space $(\Xi, \b(\Xi))$ of $Z^{\HH_r,*}$, we take
\[\Xi=\{Z\in C([0,\infty)\mapsto \HH_r^*\cup\{\Delta\}): Z_t=\Delta,\ t\ge \zeta(=\sigma_\Delta)\}\]
and $\b(\Xi)=\sigma\{Z_t, t\ge 0\}.$\quad We consider the direct product $\dot\O(D)\times \Xi$ of the measurable space $(\dot\O(D), \dot\gg_t(D))$ and $(\Xi, \b(\Xi)).$  \quad Then the set
$\Lambda=\{(\dot\gamma,Z)\in \dot\O(D)\times \Xi: \dot\gamma(0,t]\cap Z_{[0,\infty)}=\emptyset\}$ is $\dot\gg_t(D)\times \b(\Xi)$-measurable because
$$
\Lambda=\bigcup_{n=1}^\infty\bigcap_{u\in [0,t]\bigcap \QQ_+}\bigcap_{v\in \QQ_+}\{|\dot\gamma(u)-Z_v|>1/n\},
$$
where $\QQ_+$ denotes the set of positive rational numbers.

In view of \eqref{e:3.5}, $F_{z,r}(\dot\gamma)=\P_z^{\HH_r,*}(\Lambda_{\dot\gamma})$ for the $\dot\gamma$-section $\Lambda_{\dot\gamma}=\{Z\in \Xi:( \dot\gamma,Z)\in \Lambda\}$ of $\Lambda$ and so $F_{z,\gamma}(\dot\gamma)$ is $\dot\gg_t(D)$-measurable by the Fubini Theorem.

\medskip

(ii)\ By (i) and \eqref{I.0}, $\Im g_t(z)=\lim_{r\to\infty}rF_{z,r}(\dot\gamma)$ is $\dot\gg_t(D)$-measurable in $\dot\gamma$ for each $z\in D$. On the other hand,
 it is continuous in $z\in D\cup\partial_pK\cup \partial_p\dot\gamma[0,t]\cup\partial\HH$ for each $\dot\gamma \in \dot\O(D).$\quad Therefore $\Im g_t$ is
$\b(D\cup\partial_pK\cup\partial_p\dot\gamma[0,t]\cup\partial\HH)\times
\dot\gg_t(D)$-measurable in $(z,\dot\gamma).$

Since $g_t$ is obtained from $\Im g_t$
explicitly via \cite[(10.17)]{CFR},
$g_t$ enjoys the same joint measurability.

\medskip

(iii)\ This follows from (ii).   \qed

For
$\xi\in \RR$ and $\s\in S$,
we denote the probability measure
$\P_{D(\s), \xi+i0}$ on $(\dot \O(D(\s)), \dot\gg(D(\s)))$ by
$\P_{(\xi, \s)}$.

\begin{thm}\label{T:3.2}
{\rm (Time homogeneous Markov property)}
The process $ \{\W_t, t\geq 0; \P_{(\xi, \s)}, \, \xi \in \RR, \,
\s \in \cS\}$
 is $\{\dot\gg_t(D(\s(0)); t\geq 0\}$-adapted, and
\begin{equation}\label{e:3.6}
\P_{(\xi, \s)}(\W_0=(\xi, \s))=1,
\end{equation}
\begin{equation}\label{e:3.7}
\P_{(\xi, \s)}\left(\W_{t+s}\in B\;\big|\;\dot\gg_t(D(\s))\right)
=\P_{\W_t}(\W_s\in B)
\quad \hbox{for } t, s\ge 0,\ B\in \mathscr B(\RR \times \cS).
\end{equation}
\end{thm}

\pf\ $\W_t$ is $\dot\gg_t(D(\s(0))$-measurable by Lemma \ref{L:3.1}.
\eqref{e:3.6} follows from \eqref{e:3.2}.

For $D=D(\s)\in \dd$,
$\dot\gamma\in \dot \O(D)$ and $t\in [0,t_{\dot\gamma})$, $g_t^D$ is a conformal map from $D\setminus \dot\gamma[0,t]\in \wh \dd$ onto $D_t\in \dd$ sending $\dot\gamma(t)$ to $\xi(t)\in \partial\HH$ and so, by
\eqref{e:3.3}, for $\Lambda\in \dot\gg(D\setminus \dot\gamma[0,t])$ and
$z\in \partial\HH$
\begin{equation}\label{e:3.9}
\P_{D,z}(\dot\theta_t^{-1}\Lambda\; \big|\; \dot\gg_t(D))
=\P_{D_t,\xi(t)}(g_t^D(\Lambda)).
\end{equation}

Set, for $t, s\ge 0$ and $B\in \b(\RR\times \cS),$
$$
\Lambda_{t,s}= \left\{ \dot \eta \in \dot \O (D\setminus \dot \gamma [0, t]):
\dot \eta (0)=\dot \gamma (t), \ (\wh \xi (s) , \wh \s (s)) \in B  \right\}.
$$
Here,
by means of the canonical (conformal) map
$g^{D\setminus \dot \gamma [0, t]}_s$
from $(D\setminus \dot \gamma [0, t])\setminus \dot\eta [0, s]$ onto
$\wh D_s \in \dd$, we define
$\wh \xi (s) = g^{D\setminus \dot \gamma [0, t]}_s (\dot\eta (s))$ and
$\wh \s (s)= \s ( \wh D_s)$.  Then we have
$$
\theta^{-1}_t \Lambda_{t,s} = \{ \dot \gamma \in \dot \O (D): \W_{t+s} \in B\}
\quad \hbox{and} \quad
 g^D_t (\Lambda_{t,s})= \left\{ \dot \gamma \in \dot \O (D_t): \dot\gamma(0)=\xi(t),\ \W_s \in B\right\}.
$$
In fact, the first identity is due to the relation
$$ g^D_{t+s} (z) = g_{D_t \setminus g^D_t \dot \gamma [t, t+s]} \circ
g^D_t (z),$$
while the second one is obtained by the observation that $g^D_t$ induces a one-to-one map between $\dot\O(D\setminus \dot\gamma[0,t])$ and $\dot\O(D_t).$

The conclusion of the theorem now follows from \eqref{e:3.9}. \qed

\begin{remark}\rm The filtration $\{\dot\gg_t(D(\s)),\; t\ge 0\}$ in the identity \eqref{e:3.7}
depends on the second component $\s$
of the initial state $(\xi,\s).$  Nevertheless we can regard the process $(\W_t, \P_{(\xi,\s)})$ as a Markov process on $\RR\times \cS$ in a usual sense.
If
we write $\w=(\xi,\s)$ and introduce a transition function $P_t$ on $\RR\times \cS$ by
$$
P_tf(\w)=\E_\w[f(\W_t)],\quad  f\in \b_b(\RR\times \cS),
 $$
 then \eqref{e:3.7} implies that, for any $0\le t_1<t_2<\cdots<t_n$,
 $f_1,f_2,\cdots,f_n\in \b_b(\RR\times \cS)$, and $\w\in \RR\times \cS$,
$$
\E_\w \left[ \prod_{k=1}^n f_k(\W_{t_k}) \right]
= \int_{(\RR\times \cS)^n}  \prod_{k=1}^{n} f_k(\w_k )  P_{t_k-{t_{k-1}}}(\w_{k-1}, d\w_k)
$$
with $t_0:=0$ and $\w_0:=w$.  \qed
\end{remark}

\subsection{Brownian scaling for $\W$}\label{S:3.3}

\begin{lem}
For $D\in \dd, \gamma\in \O(D),$ let $a_t=a_t(\gamma,D)$
be the associated half-plane capacity.  Then for any $c>0$
\begin{equation}\label{e:3.16}
a_t(c\gamma, cD)= c^2a_t(\gamma,D),\quad t\in [0,t_\gamma).
\end{equation}
In particular, if $\gamma$ is parameterized by the half-plane the half-plane capacity,
then
\begin{equation}\label{e:3.17}
\left(\dot{c\gamma}\right)(t)= c\;\dot\gamma(c^{-2}t),\quad 0\le t < c^2 t_{\dot\gamma} =: t_{\dot{(c\gamma)}} .
\end{equation}
is the half-plane capacity  parametrization of the curve $c\gamma $ in $cD$.

\end{lem}
\pf\ Let $g_t(z)$ be the canonical map from $D\setminus \gamma.$
Then $g_t^c(z)=cg_t(z/c)$ is the canonical map from $cD\setminus c\gamma.$  \eqref{e:3.16} follows from \eqref{e:3.1} and
\[z(g_t^c(z)-z)=c^2\frac{z}{c}\left(g_t(\frac{z}{c})-\frac{z}{c}\right).\]
\eqref{e:3.17} follows from $a_t(\dot\gamma,D)=2t$ and \eqref{e:3.16}.   \qed

We make a convention that $c\;\Delta=\Delta$ for any constant $c>0.$  Then the identity \eqref{e:3.17} holds for any $t\ge 0$;\ for $t\ge c^2t_{\dot\gamma},$ the both hand sides of \eqref{e:3.17} equal $\Delta.$   Keeping this in mind, we show the following:

\begin{prop}\label{P:3.5}\ For $D\in \dd,\; z\in \partial\HH$ and any $c>0$
\begin{equation}\label{e:3.18}
\{c^{-1} \dot\gamma(c^2t),\; t\ge 0\}\ {\rm under}\ \P_{cD,cz}\ \hbox{ has the same
distribution as } \
\{\dot\gamma(t),\; t\ge 0\}\ {\rm under}\ \P_{D,z}.
\end{equation}
\end{prop}

\pf\ For a fixed $c>0$,  $f(z)=cz$ is a
conformal map from $D$ onto $cD\in \dd$.
By the
invariance under linear conformal map
\eqref{e:3.4}, we have for $D\in \dd,\ z\in \partial\HH$
\begin{equation}\label{e:3.19}
\P_{D,z}(\Lambda)=f^{-1}_*\P_{cD,cz}(\Lambda)=\P_{cD,cz}(f(\Lambda)),\quad \Lambda \in \dot\gg(D).
\end{equation}
For $\Lambda=\{\dot\gamma\in \dot \O(D): \dot\gamma\in B\}\in
\dot\gg(D)$ with $B\in \b\left((\overline \HH\cup\{\Delta\})^{[0,\infty)}\right)$,  $f(\Lambda)=\{\dot\gamma\in \dot \O(cD):
\dot{(\gamma/c)}\in B\}$.
By \eqref{e:3.17},
\begin{equation}\label{e:3.20}
  \dot{(\gamma/c)}(t)= c^{-1}  \dot \gamma(c^2t),\quad t\ge 0,
\end{equation}
and so \eqref{e:3.18} follows from \eqref{e:3.19}.  \qed

\begin{thm}\label{T:3.6}
{\rm (Brownian scaling property of $\W$)} \
For $\xi\in \RR$, $\s\in \cS$ and $c>0$
\begin{equation}\label{e:3.21}
\{c^{-1} \W_{c^2t},\; t\ge 0\}\ {\rm under}\ \P_{(c\xi, c\s)}
\  \hbox{ has the same distribution as } \
\{\W_t,\; t\ge 0\}\ {\rm under}\ \P_{(\xi, \s)}.
\end{equation}
\end{thm}
\pf\
For fixed $D\in \dd$ and $c>0,$ consider the canonical map $g_t^{cD}$ associated with $cD$ and a curve $\{\dot\gamma(t),\; t\ge 0\}\subset \dot \O(cD).$  The induced process
 $\W_t=(\xi (t), \s(t))$, $t\in [0,t_{\dot\gamma})$,
is given by $\s(t)=g_t^{cD}(c\s)$ for $\s=\s(D)$ and $\xi(t)=g_t^{cD}(\dot\gamma(t)).$

Now the curve on the left hand side of \eqref{e:3.18} defined for $\dot\gamma\in \dot\O(cD)$ belongs to $\dot \O(D)$ in view of \eqref{e:3.20} and the associated canonical map $\wt g_t^D$ from $D$ is given by
\begin{equation}\label{e:3.22}
\wt g_t^D(z)= c^{-1} g_{c^2t}^{cD}(cz),\quad z\in D,\quad{\rm for}\ t\in [0,t_{\dot\gamma/c^2}),
\end{equation}
which induces the motion $\{c^{-1}\W_{c^2t}:\
t\in [0,t_{\dot\gamma/c^2})\},$
because 
$\wt g_t^D(\s)=c^{-1}\s(c^2t)$ and $\wt g_t^D(c^{-1}\dot\gamma(c^2t))= c^{-1}\xi(c^2t)$
for  $t\in [0,t_{\dot\gamma/c^2})$.

Let $\{\W_t,\; t\ge 0\}$ be the $(\RR\times S)$-valued motion produced by $D\in \dd$ and $\dot\gamma\in \dot\O(D).$  Then, for $0\le t_1<t_2<\cdots <t_n$, $(\W_{t_1}, \W_{t_2},\cdots,\W_{t_n})$ equals an $(\RR\times S)^n$-valued $\dot\gg(D)$-measurable function $F(\dot\gamma)$ of $\dot\gamma\in \dot\O(D)$ by virtue of Lemma \ref{L:3.1}.
Therefore we can conclude from \eqref{e:3.18} and the above observation that \eqref{e:3.21} holds. \qed

\subsection{Homogeneity of $\W$ in $x$-direction}\label{S:3.4}

\begin{lem}\label{L:3.7}
For $D\in \dd, \gamma\in \O(D),$ let $a_t=a_t(\gamma,D)$ be the associated half-plane capacity.
Then for any ${r}\in \RR$,
\begin{equation}\label{e:3.23}
a_t(\gamma+{r}, D+ {r})= a_t(\gamma,D),\quad t\in [0,t_\gamma).
\end{equation}
In particular, the half-plane capacity parametrization of the curve $\gamma+{r}$ in
$D+{r}$ is given by
$\dot\gamma+{r}$; in other words,
\begin{equation}\label{e:3.24}
\dot{(\gamma+{r})}(t)= \dot\gamma(t)+{r},\quad 0\le t < t_{\dot\gamma}.
\end{equation}
\end{lem}
\pf\ Let $g_t(z)$ be the canonical map associated with $(\gamma, D).$
Then $g_t^{r}(z)=g_t(z-{r})+{r}$ is the canonical map associated with $(\gamma+{r}, D+{r}).$  \eqref{e:3.23} follows from \eqref{e:3.1} and
\[z(g_t^{r}(z)-z)=\frac{z}{z-{r}}\cdot (z-{r})\left(g_t(z-{r})-(z-{r})\right).\]
\eqref{e:3.24} follows from $a_t(\dot\gamma,D)=2t$ and \eqref{e:3.23}.   \qed

The identity \eqref{e:3.24} holds for any $t\ge 0$ because both hand sides of \eqref{e:3.24} equal $\Delta$  when  $t\ge t_{\dot\gamma}$.

\begin{prop}\label{P:3.8}\ For $D\in \dd,\; z\in \partial\HH$ and any ${r}\in \RR$
\begin{equation}\label{e:3.25}
\{ \dot\gamma(t)-{r},\; t\ge 0\}\ {\rm under}\ \P_{D+{r},z+{r}}\ \hbox{ has the same
distribution as } \
\{\dot\gamma(t),\; t\ge 0\}\ {\rm under}\ \P_{D,z}.
\end{equation}
\end{prop}
\pf\ For a fixed ${r}\in \RR$, consider the sift
$f(z)=z+{r},\quad z\in D$.
By the invariance under linear conformal map \eqref{e:3.4},
we have for $D\in \dd,\ z\in \partial\HH$
\begin{equation}\label{e:3.26}
\P_{D,z}(\Lambda)=f^{-1}_*\P_{D+{r},z+{r}}(\Lambda)=\P_{D+{r},z+{r}}(f(\Lambda)),\quad \Lambda \in \dot\gg(D).
\end{equation}
$\Lambda\in \dot\gg(D)$ can be expressed as $\Lambda=\{\dot\gamma\in \dot \O(D): \dot\gamma\in B\}$ for $B\in \b\left((\overline \HH\cup\{\Delta\})^{[0,\infty)}\right).$  Then
\[f(\Lambda)=\{\dot\gamma\in \dot \O(D+{r}):
\dot{(\gamma-{r})}\in B\}.\]
This combined with \eqref{e:3.24} and \eqref{e:3.26} leads us to \eqref{e:3.25}.  \qed

For ${r}\in \RR,$ denote by $\wh{r}$ the vector in $\RR^{3N}$ whose first $N$ entries are $0$ and the last $2N$ entries are ${r}$.   Note that
$\s(D+{r})=\s(D)+\wh{r} \quad \hbox{ for }  D\in \dd,\ {r}\in \RR.$

\begin{thm}\label{T:3.9} {\rm (Homogeneity of $(\W_t,\;\P_{(\xi, \s)})$ in $x$-direction)}
For $\xi\in \RR$, $\s\in S$ and ${r}\in \RR$,  
$\{(\xi(t)-{r}, \s(t)-\wh{r}),\; t\ge 0\}$  under $\P_{(\xi+{r}, \s+\wh{r})}$
  has the same distribution as  
$\{(\xi(t),\s(t)),\; t\ge 0\}$  under $ \P_{(\xi, \s)}$.
\end{thm}

\pf\
Fixed $D\in \dd,\ z=\xi+i0\in \partial\HH$,  ${r}\in \RR$ and put $\s=\s(D).$
Consider the canonical map $g_t^{D+{r}}$ associated with $D+{r}$ and a curve $\{\dot\gamma(t),\; t\ge 0\}\subset \dot \O(D+{r}).$  The process $\W_t=(\s(t),\xi(t)),\ t\in [0,t_{\dot\gamma}),$ being considered under $\P_{D+{r},z+{r}}=\P_{(\xi+{r},\s+\wh{r})}$ is induced from $g_t^{D+{r}}$ by
\[\xi(t)=g_t^{D+{r}}(\dot\gamma(t)),\quad
\s(t)=g_t^{D+{r}}(\s+\wh{r}).\]

Now the curve on the left hand side of \eqref{e:3.25} belongs to $\dot \O(D)$ in view of \eqref{e:3.24} and the associated canonical map $\wt g_t^D$ is given by
$$
\wt g_t^D(z)=  g_t^{D+{r}}(z+{r})-{r},\quad z\in D,\quad{\rm for}\ t\in [0,t_{\dot\gamma}).
$$
The induced motion is
\[
\begin{cases}
\wt g_t^D(\dot\gamma(t)-{r})=g_t^{D+{r}}(\dot\gamma(t))-{r}=\xi(t)-{r},\\
\wt g_t^D(\s)= g_t^{D+{r}}(\s+\wh{r})-\wh{r}=\s(t)-\wh{r}.
\end{cases}
\] 
The theorem 
now follows from \eqref{e:3.25} by the same reason as in the last paragraph of the proof of Theorem \ref{T:3.6}. \qed

\subsection{Stochastic differential equation for $\W$}\label{S:3.5}

We write $\w=(\xi, \s)\in \RR\times \cS .$
We know from
Theorem \ref{T:3.2} that $\W=(\W_t, \P_\w)$ is a time homogeneous Markov process taking values in $\RR\times \cS \subset \RR^{3N+1}.$  The sample path of $\W$ is continuous up to its
lifetime $t_{\dot\gamma}\le \infty$ owing to {\bf(P.3)} and {\bf(P.4)}. Let  $P_t$ be its transition semigroup defined as
$$
P_tf(\w)=\E_\w[f(\W_t)],  \qquad t\ge 0, \ \w\in \RR\times \cS.
$$
Denote by $C_\infty(\RR\times \cS )$ the space of all continuous functions on $\RR\times \cS $ vanishing at infinity.

In this section, we  assume that
the Markov process $\W$ satisfies properties {\bf(C.1)} and {\bf(C.2)} stated below.

\medskip
\noindent {\bf (C.1)}
\ $P_t(C_\infty(\RR\times \cS ))\subset C_\infty(\RR\times \cS ),\ t>0,\quad C_c^\infty(\RR\times \cS )\subset \dd(L),$

\medskip
\noindent
where $L$ is the infinitesimal generator of $\{P_t,\;t>0\}$ defined by
\begin{eqnarray}
Lf(\w)&=&\lim_{t\downarrow 0}\frac{1}{t}(P_tf(\w)-f(\w)),\quad \w\in \RR\times \cS ,\nonumber\\
\dd(L)&=&\{f\in C_\infty(\RR\times \cS ):\ \hbox{\rm the right hand side above}\nonumber\\
&& \hskip 0.8truein \hbox{\rm converges uniformly in}\ \w\in \RR\times \cS  \}. \label{e:3.29}
\end{eqnarray}
Under condition {\bf (C.1)},
 $\W=\{\W_t, \P_\w\}$ is a {\it Feller-Dynkin diffusion} in the sense of \cite{RW}.  In view of \cite[III, (13.3)]{RW}, the restriction $\cal L$ of $L$ to $C_c^\infty(\RR\times \cS )$ is a second order elliptic partial differential operator expressed as
\begin{equation}\label{e:3.30}
{\cal L}f(\w)=\frac12\sum_{i,j=0}^{3N} a_{ij}(\w)f_{w_iw_j}(\w)+\sum_{i=0}^{3N}b_i(\w)f_{w_i}(\w)+c(\w)f(\w),\quad \w\in \RR\times \cS ,
\end{equation}
where $a$ is a non-negative definite symmetric matrix-valued continuous function, $b$ is a vector-valued continuous function and $c$ is a non-positive continuous function.

The second assumption on $\W$ is

\medskip\noindent
{\bf(C.2)}\ $c(\w)=0,\quad \w\in \RR\times \cS$.

\medskip\noindent
This property is fulfilled if $\W$ is conservative: $\P_{D,z}(t_{\dot\gamma}=\infty)=1$ for any $D\in \dd$ and $z\in \partial\HH,$ or equivalently,
\begin{equation}\label{e:3.31}
P_t1(\w)=1\quad \hbox{\rm for any}\ t\ge 0\ {\rm and}\ \w\in \RR\times \cS .
\end{equation}
In fact, $c(\w)$ can be evaluated as
\[c(\w)=\lim_{t\downarrow 0}\frac1t (P_t1(\w)-1),\quad \w\in \RR\times \cS ,\]
according to Theorem 5.8 and its Remark in \cite{Dy}.  Hence \eqref{e:3.31} implies {\bf(C.2)}.  Condition {\bf(C.2)} means that $\W$ admits no killing inside $\RR\times \cS ,$ and so it is much weaker than the conservativeness of $\W.$

We take this opportunity to point out that the exit time $V_{\eta,x}$ employed in \cite[III, Lemma (12.1)]{RW} and in the formula following it should be corrected to be $V_{\eta,x}\wedge \zeta$, where $\zeta$ is  the lifetime, as this lemma was taken from \cite[V, Lemma 5.5]{Dy} where an exit time had been defined in the latter
form.

Recall that, for $\s=(\y,\x,\x^r )$,   $z_j=x_j+iy_j$, $z_j^r =x_j^r +iy_j$
denote the endpoints of the $j$th slit $C_j$ in $D(\s)\in \dd$.
  For $\s\in \cS $,  denote by  $\Psi_\s(z,\xi)$  the complex Poisson kernel of the Brownian motion with darning (BMD) on $D(\s).$
  The  KL equations
\eqref{e:2.11}-\eqref{e:2.13} established in \S 2 for slits can be stated as
\begin{equation}\label{e:3.37}
\s_j(t)-\s_j(0)=\int_0^t b_j(\W(s))ds,\quad t\ge 0,\quad 1\le j\le 3N,
\end{equation}
where
\begin{equation}\label{e:3.38}
b_j(\w)=
\begin{cases}
-2\pi\Im \Psi_\s(z_j, \xi),\quad & 1\le j\le N,\\
-2\pi\Re \Psi_\s(z_j, \xi),\quad & N+1\le j\le 2N,\\
-2\pi\Re \Psi_\s(z_j^r , \xi),\quad & 2N+1\le j\le 3N.
\end{cases}
\end{equation}
It follows that $b_j(\w)$ in \eqref{e:3.30} is given by the above expression \eqref{e:3.38} for $j\geq 1$ and
$a_{ij}(\w)=0$ for $i+j\geq 1$.
Thus under the condition of {\bf (C.1)} (in fact, \eqref{e:3.30}) and {\bf (C.2)},
it is known (see for example,  \cite[VII,(2.4)]{RY}) that $\W_t=(\xi (t), \s(t))$ satisfies
\begin{equation}\label{e:3.28}
\begin{cases}
d\xi (t)=\sqrt{a_{00} (\W_t)} dB_t + b_0(\W_t) dt , \\
d \s_j(t)= b_j (\W_t)) dt, \qquad j=1, \dots, 3N,
\end{cases}
\end{equation}
where $B$ is a one-dimensional Brownian motion.

\medskip
A real-valued function $u(\w)=u(\xi,\s)$ on $\RR\times \cS $ is called {\it homogeneous with degree} $0$\ (resp. $-1$)\ if
\[u(c\w)=u(\w)\quad {\rm (\ resp.}\ u(c\w)=c^{-1} u(\w)\ \hbox{\rm)\quad for any}\ c>0\ {\rm and}\ \w\in \RR\times \cS .\]
The same definition of the homogeneity is in force for a real-valued
 function $u(\s)$ on $\cS$.

\begin{lem}\label{L:3.10} Assume conditions {\bf (C.1)} and {\bf (C.2)} hold.
\begin{description}
\item{\rm(i)}\ $a_{00}(\w)$ is a homogenous function of degree $0$,
while $b_i(\w)$ is a homogenous function of degree $-1$ for every $0\le i\le 3N$.

\item{\rm(ii)} \ For every $0\le j\le 3N$,  $\xi\in \RR$, $\s\in \cS$  and ${r}\in \RR.$
\begin{equation}\label{e:3.32}
 a_{00}(\xi+{r}, \s+\wh{r})=a_{00}(\xi,\s),\qquad  b_j(\xi+{r}, \s+\wh{r})=b_j(\xi,\s).
\end{equation}
\end{description}
\end{lem}

\pf\ (i)\ By virtue of the Brownian scaling property \eqref{e:3.21}, we have $P_t(\w,E)=P_{c^2t}(c\w,cE)$. Consequently,  $P_tf(\w)=P_{c^2t}f^{(c)}(c\w)$ and ${\cal L}f(\w)=c^2 {\cal L}f^{(c)}(c\w)$,  where $f^{(c)}(\w)=f(\w/c).$  Hence we get the stated properties of the coefficients $a_{ij}$ and $b_i$ of $\cal L.$

\noindent
(ii)\ By virtue of the homogeneity 
in $x$-direction from Theorem \ref{T:3.9},  
we have $P_tf(\w)=P_tf^{r}(\w+({r},\wh{r}))$ so that $Lf(\w)=Lf^{r}(\w+({r},\wh{r}))$ where $f^{r}(\w)=f(\w-({r},\wh{r})).$  Hence we get \eqref{e:3.32}.
\qed

 \begin{remark}\label{R:3.12} \rm\ The properties of $b_j$ for $1\leq j\leq 3N$ stated in  the above lemma
can be derived without using conditions {\bf (C.1)}-{\bf (C.2)}.
In fact they follow directly from their definition \eqref{e:3.38} combined with the conformal invariance of the BMD on $D\in \dd$
established in Theorem 7.8.1 and Remark 7.8.2 of \cite{CF1}.
Indeed, let $K^*_\s(z,\xi)$, $z\in D(\s)$, $\xi\in \partial D(\s)$,
 be the Poisson kernel of the BMD on $D(\s)\in \dd$ for $\s\in \cS .$  Then, by the stated invariance of the BMD under the dilation $\phi(z)=cz$ for  $c>0$
 that maps $D(\s)$ to $D(c\s)$, we have
\[\int_{-c\eps }^{c\eps } K^*_{c\s}(cz, \xi)d\xi
=\int_{-\eps }^{\eps } K^*_{\s}(z, \xi)d\xi,\quad \eps >0.\]
Dividing the both hand side by $2c\eps $ and letting $\eps \downarrow 0,$
we get $K^*_{c\s}(cz,0)=c^{-1} K^*_\s(z,0).$  Since the complex Poisson kernel $\Psi_\s(z,\xi)$, $z\in D(\s)$, $\xi\in \partial \HH$, is the unique analytic function in $z$ with the imaginary part $K_\s^*(z,\xi)$ satisfying $\lim_{z\to\infty}\Psi_\s(z,\xi)=0,$ we obtain
\begin{equation}\label{e:3.42}
\Psi_{c\s}(cz,0)=c^{-1} \Psi_\s(z,0),\quad z\in D(\s).
\end{equation}
Therefore $b_j(0,\s)$ is homogeneous in $\s\in \cS $ with degree $-1$ for $b_j$ defined by \eqref{e:3.38}, $1\le j\le 3N.$  A similar consideration for the shift $\phi(z)=z+{r},\ {r}\in \RR,$
leads us to
\begin{equation}\label{e:3.43}
K_\s^*(z,\xi)=K^*_{\s+\wh {r}}(z+{r},\xi+{r}) \quad
\hbox{ and } \quad
\Psi_\s(z,\xi)=\Psi_{\s+\wh {r}}(z+{r},\xi+{r})
\end{equation}
for $\s\in \cS$ ,  $z\in D(\s)$ and  $\xi, {r}\in \RR$, and
so  the second property in (ii) holds.   \qed
\end{remark}

Let
$$
\alpha(\s)=\sqrt{  a_{00}(0, \s)},\qquad b(\s)= b_0(0, \s),\quad \s\in \cS .
$$
It follows from Lemma \ref{L:3.10} that $\alpha (s)$ and $ b $  are homogeneous functions on $\cS $ with degree $0$ and $-1$, respectively. Moreover,
$$ \sqrt{a_{00} (\xi, \s)}=\alpha (\s -\wh \xi) \quad \hbox{and} \quad
   b_0( \xi, \s)=b(\s -\wh \xi).
$$
Thus we have the following from \eqref{e:3.28} and Lemma \ref{L:3.10}

\begin{thm}\label{T:3.11}
 Assume conditions {\bf (C.1)} and {\bf (C.2)} hold. 
 \begin{description}
\item{\rm(i)}\ The diffusion $\W_t=(\xi(t),\s(t))$ satisfies under $\P_{(\xi,\s)}$ the following stochastic differential equation:
\begin{eqnarray}
\xi(t)&=&\xi+\int_0^t\alpha(\s(s)-\wh\xi(s))dB_s+\int_0^t b(\s(s)-\wh\xi(s))ds , \label{e:3.39}\\
\s_j(t)&=&\s_j+\int_0^t b_j(\xi(s),\s(s))ds,\quad t\ge 0,\quad 1\le j\le 3N. \label{e:3.40}
\end{eqnarray}
\item{\rm (ii)} For each $1\leq j \leq 3N$,
 $b_j (\xi, \s)$ is given by \eqref{e:3.38}, which has the properties that
 $b_j (\xi , \s)=b_j (0, \s -\wh \xi)$ and that $b_j (0, \s)$
 is a   homogeneous function
  on $\cS$ of degree $-1$.
\end{description}
\end{thm}

\section{Solution of SDE having homogeneous coefficients}\label{S:4}

We consider the following local Lipschitz condition
for a real-valued function $f=f(\s)$ on $\cS $:

\medskip
\noindent
{\bf (L)}\ For any $\s^{(0)}\in \cS $ and any finite open interval $J\subset \RR,$ there exist a neighborhood $U(\s^{(0)})$ of $\s^{(0)}$ in $\cS $ and a constant $L>0$ such that
\begin{equation}\label{e:4.1}
|f(\s^{(1)}-\wh \xi)-f(\s^{(2)}-\wh \xi)|\le L\;|\s^{(1)}-\s^{(2)}|
\quad \hbox{for } \s^{(1)},  \s^{(2)}\in U(\s^{(0)}) \hbox{ and } \xi \in J.
\end{equation}
  Recall that $\wh \xi$ denotes the vector in $\RR^{3N}$ whose first $N$-entries are $0$ and the last $2N$ entries are $\xi$.

\medskip
Recall that the coefficient $b_j(\xi,\s)$ in the equation \eqref{e:3.40}
 is defined by \eqref{e:3.38} and satisfies
 \eqref{e:3.32}.

\begin{lem}\label{L:4.1}\begin{description}
 \item{\rm(i)}\
 The function $\wt b_j(\s):=b_j(0, \s)$
 enjoys property {\rm\bf(L)} for every $1\le j\le 3N$.

 \item{\rm(ii)}\ If a function $f$ on $\cS $ satisfies the condition {\rm\bf(L)}, then
it holds for any $\s^{(1)},\ \s^{(2)}\in U(\s^{(0)})$ and for any $\xi_1,\ \xi_2 \in J$
that
\begin{equation}\label{e:4.2}
|f(\s^{(1)}-\wh \xi_1)-f(\s^{(2)}-\wh \xi_2)|\le L\left(|\s^{(1)}-\s^{(2)}|+\sqrt{2N}|\xi_1-\xi_2|\right).
\end{equation}
\end{description}
\end{lem}

\pf (i) This follows immediately from \cite[Theorem 9.1]{CFR}.

\noindent
(ii)\ Suppose a function $f$ on $\cS $ satisfies the condition {\rm\bf(L)}.
For any $\s^{(1)},\ \s^{(2)}\in U(\s^{(0)})$ and for any $\xi_1,\ \xi_2 \in J$ with $\xi_1<\xi_2,$ we have
\[
|f(\s^{(1)}-\wh \xi_1)-f(\s^{(2)}-\wh \xi_2)|\le
|f(\s^{(1)}-\wh \xi_1)-f(\s^{(2)}-\wh \xi_1)|+
|f(\s^{(2)}-\wh \xi_1)-f(\s^{(2)}-\wh \xi_2)|.
\]
 Since $\s^{(2)}\in U(\s^{(0)}),$ there exists $\delta>0$ such that $\s^{(2)}-\wh \xi\in U(\s^{(0)})$ for any $\xi\in \RR$ with $|\xi|<\delta.$  Choose points ${r}_i,\;0\le i\le \ell,$ with ${r}_0=\xi_1,\;0<{r}_i-{r}_{i-1}<\delta,\;1\le i\le \ell, \; {r}_\ell=\xi_2.$
The first term of the righthand side of the above inequality is dominated by $L|\s^{(1)}-\s^{(2)}|.$
The second term is dominated by
$ \sum_{i=1}^\ell|f(\s^{(2)}-\wh {r}_i)-f(\s^{(2)}-\wh {r}_{i-1})|
= \sum_{i=1}^\ell|f\left((\s^{(2)}-(\wh {r}_i-\wh {r}_{i-1}))-\wh {r}_{i-1})-f((\s^{(2)}-\wh {r}_{i-1}\right)|
\le \sum_{i=1}^\ell L|\wh {r}_i-\wh {r}_{i-1}|=L\sqrt{2N}(\xi_2-\xi_1).
$
\qed

In the rest of this section and throughout the next section,
we assume that we are given a non-negative homogeneous function $\alpha(\s)$ of $\s\in \cS $ with degree $0$ and a homogeneous function $b(\s)$ of $\s\in \cS $ with degree $-1$ both satisfying the condition {\bf (L)}.

\begin{thm}\label{T:4.2}\ The SDE {\rm\eqref{e:3.39}, \eqref{e:3.40}} admits a unique strong solution $\W_t=(\xi(t), \s(t)),\ t\in [0, \zeta)$, where $\zeta$ is the time when $\W_t$ approaches the point at infinity of $\RR\times \cS .$
\end{thm}
\pf\  In view of Lemma \ref{L:4.1}, every coefficient, say, $f(\xi,\s),\; \xi\in \RR,\:\s\in \cS ,$ in \eqref{e:3.39} and \eqref{e:3.40} is locally Lipschitz continuous on $\RR\times \cS \;(\subset \RR^{3N+1})$ in the following sense:
 for any $\s^{(0)}\in \cS $ and for any finite open interval $J\subset \RR$, there exists a ball $U(\s^{(0)})\subset \cS $ centered at $\s^{(0)}$ and a constant $L_0$ such that
\[|f(\xi_1,\s^{(1)}-f(\xi_2,\s^{(2)})|\le L_0(|\s^{(1)}-\s^{(2)}|+|\xi_1-\xi_2|),\quad \s^{(1)},\;\s^{(2)}\in U(\s^{(0)}),\ \xi_1,\;\xi_2\in J.\]
Thus
\eqref{e:3.39} and \eqref{e:3.40} admit a unique local solution.  It then suffices to patch together those local solutions just as in \cite[Chapter V, \S 1]{IW}.
\qed

\begin{prop}\label{P:4.3}\ The solution $\W_t=(\xi(t), \s(t)),\;t\in [0, \zeta),$ of the SDE {\rm\eqref{e:3.39}, \eqref{e:3.40}} enjoys the following properties:

\noindent
{\rm(i)}\ {\rm(Brownian scaling property)}\ For $\s\in \cS ,\;\xi\in \RR$ and for any $c>0,$
$$
 \{c^{-1} \W_{c^2t},\;t\ge 0\} \hbox{ under } \P_{(c\xi,c\s)}
 \hbox{ has the same distribution as } \
\{\W_t,\;t\ge 0\} \hbox { under } \P_{(\xi,\s)}.
$$

\noindent
{\rm(ii)}\ {\rm(homogeneity in}\ $x$-{\rm direction)}\ For $s\in \cS , \xi\in \RR$ and for any ${r}\in \RR$,
\[ \{(\xi(t)-{r}, \s(t)-\wh{r}), t\ge 0\}
\hbox { under }  \P_{(\xi+{r},\s+\wh{r})}
 \hbox{ has the same distribution as } \
\{(\xi(t),\s(t))\;t\ge 0\} \hbox{ under }  \P_{(\xi,\s)}.
\]
\end{prop}

\pf \ (i)\ We put $\W_c(t)=c^{-1}\W(c^2t)=(\xi_c(t), \s_c(t))$ with
$\xi_c(t)=c^{-1}\xi(c^2t),\ \s_c(t)=c^{-1}\s(c^2t).$  $\W(t)=(\xi(t),\s(t))$ under $\P_{(c\xi,c\s)}$ satisfies the equation \eqref{e:3.39} with $c\xi$ in place $\xi.$  Hence, by taking the homogeneity of $\alpha, b$ into account, we get
\begin{eqnarray*}
\xi_c(t)&=&\xi+c^{-1}\int_0^{c^2t}\alpha(\s(s)-\wh\xi(s))dB_s+c^{-1}\int_0^{c^2t}b(\s(s)-\wh\xi(s))ds\\
&=&\xi+c^{-1}\int_0^t\alpha(c(\s_c(s)-\wh\xi_c(s)))dB_{c^2s}+c\int_0^tb(c(\s_c(s)-\wh\xi_c(s)))ds\\
&=&\xi+\int_0^t\alpha(\s_c(s)-\wh\xi_c(s))d\tilde B_s+\int_0^tb(\s_c(s)-\wh\xi_c(s))ds,
\end{eqnarray*}
where $\tilde B_s=c^{-1}B_{c^2s}.$  Therefore the equation \eqref{e:3.39} with a new Brownian motion $\tilde B_s$
is satisfied by $\W_c(t)$ under $\P_{(c\xi,c\s)}.$
Similarly, \eqref{e:3.40}
is also satisfied by $\W_c(t)$ under $\P_{(c\xi,c\s)}.$

\noindent
(ii)\ This is immediate from the expressions \eqref{e:3.39} and \eqref{e:3.40} of the SDE and the property
\eqref{e:3.32}. \qed

\section{Stochastic Komatu-Loewner evolutions}\label{S:5}

\subsection{Stochastic Komatu-Loewner evolutions}\label{S:5.1}

Let us fix a pair of functions $(\xi(t), \s(t)),\ t\in [0,\zeta),$ taking values in $\RR\times \cS $ satisfying the two following properties {\bf(I)} and  {\bf (II)}:

\medskip\noindent
{\bf(I)}\ \hskip 0.2truein  $\xi(t)$ is a real-valued continuous function of $t\in [0,\zeta)$.

\medskip \noindent
{\bf(II)} \hskip 0.1truein $(\xi(t), \s(t)),\ t\in [0,\zeta),$ satisfies the equation \eqref{e:3.40}  with
 $b_j$, $1\le j\le 3N$, given by \eqref{e:3.38}.

\medskip

We have freedom of choices of such a pair in two ways.  \\
The first way is to take any deterministic real continuous function $\xi(t),\ t\in [0,\infty)$, substitute it into the right hand side of \eqref{e:3.40} and get the unique solution $\s(t)$ on a maximal time interval $[0,\zeta)$ of the resulting ODE
by using Lemma \ref{L:4.1}. \\

The second way is to choose any solution path $\W_t=(\xi(t), \s(t)),\ t\in [0,\zeta),$
of the SDE \eqref{e:3.39} and \eqref{e:3.40} obtained in Theorem 4.2 for a given homogeneous functions $\alpha$ and $b$ on $\cS $ with degree $0$ and $-1$, respectively, both satisfying condition {\bf (L)}.

\smallskip
We write $D_t=D(\s(t))\in \dd$, $t\in [0,\zeta),$ and define
\begin{eqnarray*}
 G&=&\bigcup_{t\in [0,\zeta)}\; \{t\}\times D_t, \\
 \wh G &=& \bigcup_{t\in [0,\zeta)}\; \{t\}\times \left( D_t\cup \partial_pK(t)\cup
(\partial\HH\setminus \{\xi(t)\})\right),
\end{eqnarray*}
where  $K(t)=\cup_{j=1}^N C_j(t)$ and $D_t=\HH\setminus K(t)$.
For each $1\le j\le N$, let  $\partial_p C_j^0(t)=C_j^{0,+}(t)\cup C_j^{0,-}(t)$
 denote the set $\partial_p C_j(t)$ with its two endpoints being removed,
 and $\partial_p K^0(t):=\cup_{j=1}^N \partial_pC_j^0(t).$
Note that
$G$ is a domain of $[0,\zeta)\times \HH$ in $\RR^3$ because $t\mapsto D_t$ is continuous.

We first study  the unique existence of local solutions $z(t)$ of the following equation
\begin{equation}\label{e:5.1}
\frac{d}{dt}z(t)=-2\pi \Psi_{\s(t)}(z(t), \xi(t))
\end{equation}
with initial condition
 \begin{equation}\label{e:5.2}
z(\tau)=z_0\in D_{\tau}\cup\partial_pK^0(\tau)\cup(\partial\HH\setminus \xi(\tau))
\end{equation}
for $\tau\in [0,\zeta)$.

\begin{prop}\label{P:5.1}\begin{description}
 \item{\rm(i)}\ $\Psi_{\s(t)}(z,\xi(t))$ is jointly continuous in $(t,z)\in \wh G.$

\item{\rm(ii)}\ $\lim_{z\to\infty}\Psi_{\s(t)}(z,\xi(t))=0$
uniformly in $t $ in every  finite time interval $I \subset [0,\zeta).$

\item{\rm(iii)}\ $\Psi_{\s(t)}(z,\xi(t))$ is locally Lipschitz continuous
in $z$ in the following sense: for any $(\tau,z_0)\in G$, there exist $t_0>0$, $\rho>0$ and $L>0$ such that
$$
V=[(\tau-t_0)^+,\tau+t_0]\times \{z:|z-z_0|\le \rho\}\subset G
$$
 and
\begin{equation}\label{e:5.4}
|\Psi_{\s(t)}(z_1,\xi(t))-\Psi_{\s(t)}(z_2,\xi(t))|\le L\, |z_1-z_2|
\quad \hbox {\rm for any}\ (t,z_1),\;(t,z_2) \in V.
\end{equation}

\item{\rm(iv)}\  Fix $1\le j\le N.$  For any $\tau\in [0,\zeta)$ and
$z_0\in C_{j}^{0,+}(\tau)$,
there exist $t_0 >0$, $L>0$ and an open rectangle $R\subset \HH$ with sides parallel to
the axes
centered at $z_0$ such that
$$
R\cap C_j^+(t)\neq \emptyset\ {\rm and}\ R\cap C_j^+(t) \subset C_j^{0,+}(t)\ \hbox{\rm for every}\ t\in [(\tau-t_0)^+,\tau+t_0],
$$
 and the function $\Psi_{\s(t)}^+(z,\xi(t))$ satisfies {\rm \eqref{e:5.4}} for any $(t,z_1),\;(t,z_2)\in V_j,$ where $V_j=[(\tau-t_0)^+,\tau+t_0]\times R$ and
$\Psi_{\s(t)}^+(z,\xi(t))$ is the extension of $\Psi_{\s(t)}(z,\xi(t))$ from the upper side of $R\setminus C_j^+(t)$ to $R$ by the Schwarz reflection for each $t\in [(\tau-t_0)^+, \tau+t_0],$

An analogous statement holds for $z_0\in C_{j}^{0,-}(\tau)$.

\item{\rm(v)}\   For any $\tau\in [0,\zeta)$ and $z_0\in \partial\HH\setminus \{\xi(\tau)\},$
there exist $t_0>0$, $\rho>0$ and $L>0$
 such that
$$
V_0=[(\tau-t_0)^+,\tau+t_0]\times \{z\in \overline \HH: |z-z_0|\le \rho\} \subset \bigcup_{t\in[(\tau-t_0)^+, \tau+t_0]} \{t\}\times (D_t\cup(\partial \HH\setminus \{\xi(t)\})
$$
and {\rm \eqref{e:5.4}} holds for any $(t,z_1),\;(t,z_2)\in V_0.$

\item{\rm(vi)}\  For every $\tau\in [0,\zeta)$ and
$z_0\in D_\tau\cup (\partial\HH\setminus \xi(\tau))$,
 there exists a unique local solution
  $\{z(t); t\in (\tau -t_0, \tau +t_0)\cap [0, \zeta) \}$
 of  \eqref{e:5.1}  and \eqref{e:5.2}
 satisfying $z(\tau)=z_0$.

\item{\rm(vii)}\ Fix $1\le j\le N.$
For each initial time $\tau\in [0,\zeta)$ and initial position
$z_0\in C_j^{0,+}(\tau),$
there exists a unique local solution
 $\{z(t); t\in (\tau -t_0, \tau +t_0)\cap [0, \infty ) \}$
of the equation {\rm\eqref{e:5.1}}
 with
 $\Psi^+_{\s(t)}(z,\xi(t))$ in place of $\Psi_{\s(t)}(z,\xi(t))$
 and
  $z(\tau)=z_0.$

An analogous statement holds for $z_0\in C_j^{0,-}(\tau).$
\end{description}
\end{prop}

\pf \ (i)\ This can be shown in the same way as that for  {\bf(P.5)} in \cite[\S 9]{CFR} using the continuity of $t\mapsto D_t=D(\s(t))$.

\noindent

(ii)\ Take $R>0$ sufficiently large a large so that the closure of the set $\cup_{t\in I}(\cup_{j=1}^N C_j(t))\cup \xi(t)$ is contained in $B({\bf 0},R)=\{z\in \CC:|z|<R\}.$
Extend the analytic function $h(z,t)=\Psi_{\s(t)}(z,\xi(t))$ from $\HH\setminus B({\bf 0},R)$ to $\CC\setminus B({\bf 0},R)$ by the Schwarz reflection.  By (i), $M=\sup_{z\in \partial B({\bf 0},R), t\in I}|h(z,t)|$ is finite.
Define $\wh h(z,t)=h(1/z,t),\ |z|>R.$  Since $h(z,t)$ tends to zero as $z\to\infty$, $\wh h(z,t)$ is analytic on $B({\bf 0},1/R)$ and, by \cite[(28)-(29) in Chapter 4]{A},
\[
\frac1z\wh h(z,t)=\frac{1}{2\pi i}\int_{|\zeta|=1/R} \frac{\wh h(\zeta,t)}
{\zeta(\zeta-z)}d\zeta=\frac{R}{2\pi}\int_0^{2\pi}\frac{h(Re^{i\theta},t)}{(e^{i\theta}-R^2 z)}d\theta .
\]
Consequently,
\begin{equation}\label{e:5.7}
\sup_{t\in I}\;|z\Psi_{\s(t)}(z,\xi(t))|\le 2RM\quad{\rm if}\quad |z|\ge 2R^2.
\end{equation}

\noindent
(iii)\ $\Psi_{\s(t)}(z,\xi(t))$
is jointly continuous by virtue of (i) and
analytic in $z\in D_t,$
Therefore we readily get (iii) from the Taylor expansion \cite[(28)-(29) of Chapter 4]{A} for $n=1$ again.

\noindent
For (iv) and (v), we extend
$\Psi_{\s(t)}(z,\xi(t))$
using Schwarz reflections.

\noindent
  (vi) and (vii) follow from (iii), (iv) and (v).
\qed

\begin{lem}\label{L:5.2}\  {\rm(i)}\ Fix $1\le j\le N.$   For any $\tau\in [0,\zeta)$ and $z_0\in C_j^{0,+}(\tau)$, there exists
a unique
 solution $z(t),\ t\in [(\tau-t_0)^+, \tau+t_0],$ of {\rm\eqref{e:5.1}} and {\rm\eqref{e:5.2}} for some $t_0>0$ such that
\begin{equation}\label{e:5.8}
 z(\tau)=z_0,\quad z(t)\in C_j^{0,+}(t)\ \hbox{\rm for every}\ t\in [(\tau-t_0)^+,\tau+t_0].
\end{equation}
An analogous statement holds for $z_0\in C_j^{0,-}(\tau).$

\noindent
{\rm(ii)}\ For any $\tau\in [0,\zeta)$ and $z_0\in \partial\HH\setminus \{\xi(\tau)\}$, there exists
a unique
solution $z(t),\ t\in [(\tau-t_0)^+, \tau+t_0],$ of {\rm\eqref{e:5.1}} and {\rm\eqref{e:5.2}} for some $t_0>0$ such that
\begin{equation}\label{e:5.9}
 z(\tau)=z_0,\quad z(t)\in \partial\HH\setminus \{\xi(t)\}\ \hbox{\rm for every}\ t\in [(\tau-t_0)^+,\tau+t_0].
\end{equation}
\end{lem}

\pf\ (i)\ In view of the explicit expression (5.2) in {[CFR]},
when $z\in \partial_p C_j (t)$,   $\Im \Psi_{\s(t)}(z,\xi(t))$ is a bounded function $\eta (t)$ of $t$ independent of $z$.
Thus \eqref{e:5.1} under the requirement \eqref{e:5.8} becomes
 $\Im z(t)= z_0\exp\left(-\int_{t_0}^t \eta(s)ds\right)$ and
\begin{equation}\label{e:5.10}
 \frac{d}{dt} \Re z(t)= -2\pi \Re \Psi_{\s(t)}(\Re z(t)+i\Im z(t), \xi(t)).
\end{equation}
 Equation \eqref{e:5.10} has  a unique solution for $\Re z(t)$ in view of   Proposition \ref{P:5.1}.

\noindent
(ii)\ By (5.2) in {[CFR]}, we have $\Im\Psi_{\s(t)}(z,\xi(t))=0$ for $z\in \partial\HH\setminus \{\xi(t)\}.$
Hence the equation \eqref{e:5.1}
under the requirement \eqref{e:5.9} implies that $\Im z (t)=0$ and
\begin{equation}\label{e:5.11}
\frac{d}{dt} \Re z(t)= -2\pi \Re \Psi_{\s(t)}(\Re z(t) , \xi(t)).
\end{equation}
The above equation uniquely determines $\Re z(t)$ in view of Proposition \ref{P:5.1}. \qed

Denote by $z_j(t)$ and $z_j^r (t)$ the left and right endpoints of the $j$th slit $C_j(t)$ of $\s(t)$.
We know from \eqref{e:3.38} and \eqref{e:3.40}
\begin{equation}\label{e:5.12}
\frac{dz_j(t)}{dt}=-2\pi \Psi_{\s(t)}(z_j(t), \xi(t)),\quad t\in [0,\zeta).
\end{equation}

A solution $\{z(t),\; t\in I\}$ of the equation \eqref{e:5.1} for a time interval $I\subset [0,\zeta)$ is said to {\it pass through} $G \subset
\RR^3$ if $(t,z(t))\in G$ for every $t\in I.$

\begin{lem}\label{L:5.3}\  Fix $ 1\le j\le N$ and let $I=(\alpha,\beta)$ be a finite open subinterval of $[0,\zeta).$ Let

\begin{description}
\item{\rm(i)} Suppose that
$\{z(t);  t\in I\}$
is a solution of \eqref{e:5.1} passing through $\wh G$
with $z(\beta)= z_j (\beta)$ but $z(t)\not= z_j (t)$ for $t\in (\alpha, \beta)$.
Then there exists $t_0 \in (0, \beta-\alpha)$ so that $z(t) \in \partial_p  C^0_j (t) $ for $t\in [\beta-t_0,\beta)$.
The same conclusion holds if $z_j (\beta)$ and $z_j(t)$ are replaced
by $z_j^r  (\beta)$ and $z_j^r (t)$.

\item{\rm(ii)}\  Suppose that
$\{z(t);  t\in I\}$
is a solution of \eqref{e:5.1} passing through $\wh G$
with $z(\alpha )= z_j (\alpha )$ but $z(t)\not= z_j (t)$ for $t\in (\alpha, \beta)$.
 Then there exists
$t_0 \in (0, \beta-\alpha)$ so that $z(t) \in \partial_p  C^0_j (t) $ for $t\in [\alpha, \alpha + t_0)$.
The same conclusion holds if $z_j (\alpha)$ and $z_j(t)$ are replaced
by $z_j^r  (\alpha)$ and $z_j^r (t)$.
\end{description}
\end{lem}

\pf\  We only prove (i) as the proof for (ii) is analogous.
For $\zeta\in \CC$ and $\eps >0$, we use $B(\zeta , \eps)$ to denote
the ball   $\{z\in \CC:|z-\zeta|<\eps \}$ centered at $\zeta$ with radius $\eps$.

Suppose that  $\{z(t),\; t\in [\beta-t_1,\beta) \}$ a solution  of \eqref{e:5.1} passing through $G$ and that $z(\beta)=z_j(\beta)$.
Taking $t_1$ smaller if needed, we may assume that there is
$\eps>0$  so that
\begin{equation}\label{e:5.13}
B(z_j(t), \eps)\subset \HH \quad \hbox{ and } \quad  z_j^r (t)\notin B(z_j(t),\eps)
\quad \hbox{for every }  t\in [\beta-t_1, \beta].
\end{equation}
We can further choose $t_0\in (0,t_1]$ so that
\begin{equation}\label{e:5.14}
z(t)\in B(z_j(t),  \eps/2) \cap (D_t \cup \partial_p C_j (t) )
\quad\hbox{ for every }  t\in [\beta-t_0,\beta ].
\end{equation}
For each $t\in (\beta-t_1,\beta]$, let
\[\psi_t(z)=\sqrt{z-z_j(t)}\;:\;  B(z_j(t),\eps)\setminus C_j(t)   \to B({\bf 0},\sqrt{\eps})\cap   \HH,\]
and
\begin{eqnarray*}
f_t(z)
=\Psi_{\s(t)}(\psi_t^{-1}(z),\xi(t))
=\Psi_{\s(t)}(z^2+z_j(t),\xi(t))\, :\; B({\bf 0},\sqrt{\eps})\cap\HH \to \CC.
\end{eqnarray*}
Then $f_t$ is an analytic function on $B({\bf 0},\sqrt{\eps})\cap\HH$, which can be extended to be
an analytic function on $B({\bf 0},\sqrt{\eps})\setminus \{{\bf 0}\}$ by the Schwarz reflection because $\Im f_t(z)$ is constant on $B({\bf 0},\sqrt{\eps})\cap \partial \HH$.
On account of \cite[Chap. 4 (28),(29)]{A},
it holds
for every $a\in (\sqrt{\eps}/2,\sqrt{\eps})$ and $z\in B({\bf 0}, a)$,
\begin{equation}\label{e:5.15}
f_t(z)-f_t({\bf 0})=z\; h_t(z) \quad \hbox{with }\
h_t(z)=  \frac{1}{2\pi i}\;\int_{\partial B({\bf 0}, a)}  \frac{f_t(\zeta)}{\zeta(\zeta-z)} d\zeta,\quad z\in B({\bf 0}, a),
\end{equation}
In particular, $|h_t'(z)|$ is uniformly bounded in $(z, t)\in B({\bf 0},\sqrt{\eps}/2) \times [\beta-t_0,\beta]$
in view of Proposition 6.1(i). Accordingly $h_t(z)$ is Lipschitz continuous on $B({\bf 0},\sqrt{\eps}/2)$ uniform in $t\in [\beta-t_0,\beta]$:
\begin{equation}\label{e:5.16}
|h_t(z_1)-h_t(z_2)|\le L|z_1-z_2|,\quad z_1,  z_2 \in B({\bf 0},\sqrt{\eps}/2),
\end{equation}
for a constant $L>0$ independent of $t\in [\beta-t_0,\beta].$

We now let
$\wh z(t)=\psi_t(z(t))  =\sqrt{z(t)-z_j(t)}$ for $t\in [\beta-t_0, \beta)$.
On account of \eqref{e:5.14}, $\wh z(\beta)=0$,
\begin{equation}\label{e:5.17}
\wh z(t)\in B({\bf 0},\sqrt{\eps}/2)\cap \overline \HH
 \quad\hbox{\rm for every}\ t\in [\beta-t_0, \beta),
\end{equation}
and
\begin{equation}\label{e:5.18}
\frac{dz(t)}{dt}=-2\pi \Psi_{\s(t)}(z(t),\xi(t))=-2\pi f_t(\wh z(t)),\quad
t\in [\beta-t_0, \beta).
\end{equation}
By \eqref{e:5.12}, $\frac{d z_j(t)}{dt}=-2\pi f_t({\bf 0}).$  Therefore we have by \eqref{e:5.15}, \eqref{e:5.17} and \eqref{e:5.18} that  for any $t\in [\beta-t_0, \beta)$
$$
\frac{d\wh z(t)}{dt}=
\frac{1}{2 \wh z(t)}\left(\frac{d z(t)}{dt} - \frac{d z_j(t)}{dt}\right)
= -\frac{\pi}{\wh z(t)}\left(f_t(\wh z(t)) - f_t({\bf 0})\right)
=-\pi\; h_t(\wh z(t)).
$$
Since $h_t(z)$ is Liptschitz on $B({\bf 0},\sqrt{\eps}/2)$ uniform in $t\in [\beta-t_0,\beta]$, the solution $\wh z (t)$ to the above equation with
$\wh z(\beta)={\bf 0}$ exists and is unique. On the other hand,  note that $\Im ( f_t (z)- f_t ({\bf 0}))=0$ on $B({\bf 0}, \sqrt{ \eps})\cap \partial \HH$.
Thus by \eqref{e:5.15}
\begin{equation}\label{e:5.16b}
\Im h_t (z)=0 \qquad \hbox{on } B({\bf 0}, \sqrt{ \eps})\cap \partial \HH .
\end{equation}
It follows that the unique solution $\wh z$
to $\frac{d\wh z(t)}{dt}=-\pi\; h_t(\wh z(t))$ with $\wh z (\beta ) = {\bf 0}$
is real-valued.  It follows then $z(t) \in \partial_p C_j (t)$.
A similar argument shows that the second part of (i) holds as well.
\qed

\medskip

Due to (i) and (iii) of Proposition \ref{P:5.1}, and a general theorem in ODE (see e.g. \cite{H}),
there exists, for each $(\tau,z_0)\in G,$ a unique solution $z(t)$ of the equation \eqref{e:5.1} satisfying the initial condition $z(\tau)=z_0$ and passing through $G$ with a maximal time interval
$I_{\tau,z_0} (\subset [0,\zeta)) $ of existence.  Such a solution of \eqref{e:5.1} will be designated by
$\varphi(t; \tau,z_0)$, $t\in I_{\tau,z_0}$.
Let $\alpha$ and $\beta$ be the left and right endpoints of $I_{\tau,z_0},$ respectively, both depending on $(\tau, z_0).$  Then
$(t, \varphi (t;\tau,z_0))\in G$ for any $t\in I_{\tau,z_0}\setminus \{\alpha, \beta\}.$

\begin{prop}\label{P:5.4}\ For any $(\tau, z_0)\in G,$ the maximal time interval $I_{\tau,z_0}$ of existence of the unique solution $\varphi (t;\tau,z_0)$ of {\rm\eqref{e:5.1}} with $\varphi(\tau; \tau, z_0)=z_0$ passing through $G$ is
$ [0, \beta)$ for some $\beta > \tau$
and
\begin{equation}\label{e:5.21}
\lim_{t\uparrow \beta} \Im\varphi (t;\tau,z_0)=0,
\quad  
\lim_{t\uparrow \beta}
|\varphi (t;\tau, z_0)-\xi(\beta )|=0\ \hbox{ whenever }  \beta <\zeta.
\end{equation}
\end{prop}

\pf\ Fix $\beta_0\in (0,\zeta)$ and $z_0\in D_{\beta_0}$. Let $(\alpha, \beta)$ be
the largest subinterval of $(0, \zeta)$ so that the equation \eqref{e:5.1} has a unique solution $z(t)=\varphi (t;\beta_0, z_0)$ in $t\in (\alpha,\beta)$ satisfying $z(\beta_0)=z_0$
and passing through $G$.
By (i) and (iii) of  Proposition \ref{P:5.1}, such an interval $(\alpha, \beta)$ exists with $0\leq \alpha < \beta_0 < \beta \leq \zeta$.
For simplicity, we write $\varphi (t;\tau,z_0)$ as $\varphi(t)$.
We claim that
\begin{equation} \label{e:5.19}
\alpha=0 \quad \hbox{and} \quad  \varphi(0+ ):=\lim_{t\downarrow 0}
\varphi (t )\in D .
\end{equation}

Since the right hand side of the equation \eqref{e:5.1} is negative, $\Im \varphi (t)$ is decreasing in $t$. By (i) and (ii) of Proposition \ref{P:5.1},
$\varphi(\alpha+ ):=\lim_{t\downarrow \alpha}
\varphi (t )$ exists with $\Im \varphi(\alpha+ )>0$.
Set $\varphi(\alpha )=\varphi(\alpha+ )$, which takes value
in $D_{\alpha}\cup \cup_{j=1}^N \partial_p  C_j (\alpha)$.
By Proposition \ref{P:5.1}(vi), Lemma \ref{L:5.2}(i) and Lemma \ref{L:5.3},
$\varphi(\alpha ) \notin \cup_{j=1}^N \partial_p C_j (\alpha)$
as $\varphi(t )\in D_t$ for $t\in (\alpha, \beta_0)$.
Thus $\varphi(\alpha )\in D_\alpha$. If $\alpha >0$, then
the solution $\varphi(t)$ of \eqref{e:5.1} can be extended
to $(\alpha -\eps, \beta_0]$ for some $\eps \in (0, \alpha)$. This contradicts
to the maximality of $(\alpha, \beta)$.  Thus $\alpha=0$ and the claim \eqref{e:5.19} is  proved.

Since $\Im \varphi (t )$ is decreasing in $t$,
$\lim_{t\uparrow \beta} \Im \varphi (t)$ exists.
Assume $\beta <\zeta$. Were $\lim_{t\uparrow \beta} \Im \varphi (t)
>0$,  it follows from
(i) and (ii) of Proposition \ref{P:5.1} that
$\varphi(\beta- ):=\lim_{t\uparrow \beta} \varphi (t)$
exists and takes value in $D_\beta \cup \cup_{j=1}^N \partial_p C_j (\beta)$.
By Proposition \ref{P:5.1}(vi), Lemma \ref{L:5.2}(i) and Lemma \ref{L:5.3}
again,
$\varphi(\beta- ) \notin \cup_{j=1}^N \partial_p C_j (\beta)$
as $\varphi(t  )\in D_t$ for $t\in (\beta_0, \beta)$.
Hence $\varphi(\beta- )\in D_\beta$ and thus the solution
$\varphi (t)$ of \eqref{e:5.1} can be extended
to $[\beta_0, \beta+\eps)$ for some $\eps \in (0, \zeta-\beta )$.
This contradicts to the maximality of $(\alpha, \beta)$ and so
$ \lim_{t\uparrow \beta} \Im \varphi (t )=0$.

We now proceed to prove the second claim in \eqref{e:5.21}.
 Suppose $\limsup_{t\uparrow \beta}|\varphi(t)-\xi(\beta)|>0$.
Then by the continuity of $\xi$,  
 $\limsup_{t\uparrow \beta} |\varphi(t)-\xi(t)|>0$.
 Thus there is an $\eps >0$ and a sequence $\{t_n; n\geq 1\}\subset (\beta-\eps, \beta)$ increasing to $\beta$ so that 
 $\inf_{s\in [\beta-\eps, \beta]}|\varphi(t_n)-\xi(s)|>\eps$ for every $n\geq 1$.
 By  (i) and (ii) of Proposition \ref{P:5.1}, $\Psi_{\s (t)}(z. \xi (t))$ is bounded on 
 $$
 \wh G_0:=\Big\{(s, z)\in \wh G: s\in [\beta-\eps, \beta], 
 \inf_{s\in [\beta-\eps, \beta]}|z-\xi(s)|\geq \eps/2 \Big\},
 $$
 say, by $M>0$.  So as long as $(t, \varphi (t))\in \wh G_0$, $|\frac{d}{dt} \varphi  (t)|\leq 2\pi M$.
 Let $\delta = \eps /(4\pi M)$.
 This observation implies that $|\varphi (t_n)-\varphi (t)|\leq 2\pi M (t-t_n)\leq \eps/2$ 
 for every $t\in [t_n, t_n+\delta]\cap [t_0, \beta)$.
 Consequently,  $\varphi(\beta-)=\lim_{t\uparrow \beta}\varphi(t)$ exists and
takes value in $\partial\HH \setminus \{\xi(\beta)\}.$
But this contradicts to Proposition \ref{P:5.1}(vi) and Lemma \ref{L:5.2}(ii)
as $\varphi (t)\in D_t$ for $t\in [t_1, \beta)$. This implies that 
$\lim_{t\uparrow \beta}|\varphi(t)-\xi(\beta)|=0$.
 \qed

We write $D_0=D(\s(0))\in \dd$ as $D.$

\begin{thm}\label{T:5.5}\begin{description}
 \item{\rm(i)}\ For each $z\in D,$ there exists a unique solution $g_t(z),\; t\in [0, t_z),$ of the equation
\begin{equation}\label{e:5.25}
 \partial_t g_t(z)=-2\pi \Psi_{\s(t)}(g_t(z), \xi(t)) \quad \hbox{with }  g_0(z)=z\in D
\end{equation}
passing through $G$, where $[0,t_z),\ t_z>0,$ is the maximal time interval of its existence.  It further holds that
\begin{equation}\label{e:5.26}
\lim_{t\uparrow t_z} \Im g_t(z)=0,
\quad \lim_{t\uparrow t_z}
|g_t(z)-\xi(t_z)|=0\ {\rm whenever}\ t_z<\zeta.
\end{equation}

\item{\rm(ii)}\ Define
\begin{equation}\label{e:5.27}
F_t=\{z\in D: t_z\le t\},\quad t>0.
\end{equation}
Then $D\setminus F_t$ is open and $g_t$ is a conformal map from $D\setminus F_t$ onto $D_t$ for each $t>0.$
\end{description}
\end{thm}

\pf (i)\ This just follows from Proposition \ref{P:5.4}
with $(\tau, z_0)=(0, z)$.

\noindent
(ii)\   Since $\Psi_{\s(t)}( z, \xi (t))$ is analytic in $z$ and jointly continuous in $(t, x)$ by Proposition \ref{P:5.1}, by a general theorem on ODE (see e.g \cite{CL}), $g_t(z)$ is continuous in $(t, z)\in [0,t_z) \times D$ (and so
$D\setminus F_t=g_t^{-1}(D_t)$ is open) and $g_t (z)$ is
analytic in $D\setminus F_t$.
It follows from Proposition \ref{P:5.4} that $g_t$ is a one-to-one   map from   $D\setminus F_t$ onto $D_t$.
\qed

Note that  the complex Poisson kernel of the absorbing Brownian motion (ABM) in $\HH$ is
\begin{equation}\label{e:5.28}
\Psi^\HH(z,\xi)=-\frac{1}{\pi} \frac{1}{z-\xi},\quad z\in \HH,\ \xi \in\partial \HH,
\end{equation}
whose imaginary part $P(z,\xi):= \Im\Psi^\HH(z,\xi)=\frac{1}{\pi}\frac{y}{(x-\xi)^2+y^2}$  is the Poisson kernel of ABM in $\HH$.

Let $I$ be a finite subinterval of $[0,\zeta)$, and  
$R, M$ be the  positive constants in  the proof of Proposition \ref{P:5.1}(ii).

\begin{lem}\label{L:5.6}\ {\rm(i)}\ Let  $\wt M:=\sup_{t\in I}|\xi(t)|$.
  Then
 $$ \sup_{t\in I}\left| \frac{\Psi_{\s(t)}(z,\xi(t))}{\Psi^\HH(z,\xi(t))}\right|\le 4\pi RM \quad \hbox{for }  |z|\ge 2R^2\vee \wt M.
 $$

\noindent
{\rm(ii)}\ For any $R_1\ge R,$
\[ \sup_{t\in I}\sup_{z\in D_t,\; |z|\le R_1}|\Psi_{\s(t)}(z,\xi(t))-\Psi^\HH(z,\xi(t))|<\infty.\]
\end{lem}

\pf \ (i)\ This follows from \eqref{e:5.7} as
$$
\left| \frac{\Psi_{\s(t)}(z,\xi(t))}{\Psi^\HH(z,\xi(t))}\right|
= \pi |z-\xi(t)| \, |\Psi_{\s(t)}(z,\xi(t))|
\le 2 \pi |z\Psi_{\s(t)}(z,\xi(t)|
\quad \hbox{for } t\in I \hbox{ and }  |z|\ge \wt M.
$$

\noindent
(ii)\ For $z\in D_t=D(\s(t))$ and $\xi\in \partial\HH$, let
\[
\H_t(z,\xi)=\Psi_{\s(t)}(z,\xi)-\Psi^\HH(z,\xi), \quad v_t(z,\xi)=K_t^*(z,\xi)-P(z,\xi),
\]
where $K_t^*(z,\xi)=\Im \Psi_{\s(t)}(z,\xi)$, which
is the BMD-Poisson kernel on $D_t.$
Since $\Im \H_t(z,\xi)=v_t(z,\xi)$ vanishes for $z\in \partial\HH\setminus \{\xi\}$, by the Schwarz reflection for each $\xi \in \partial \HH$,
we extend $z\mapsto \H_t(z,\xi)$ analytically to
$D_t\cup \Pi D_t\cup (\partial \HH\setminus \{\xi\})$ which is still denoted
as $\H_t (z, \xi)$.
Here $\Pi$ denotes the mirror reflection with respect to the $x$-axis in the plane.
On the other hand, it follows from the explicit expression of
$v_t(z,\xi)$ given  by (5.2) and (12.24) from \cite{CFR} that $z\mapsto v_t(z,\xi)=\Im \H_t(z,\xi)$ is bounded in a neighborhood of $\xi$.
Hence $\xi$ is a removable singularity of $\H_t(z,\xi)$ and so $\H_t(z,\xi)$ is analytic for $z\in D_t\cup \Pi D_t\cup \partial \HH.$

Choose $\eps>0$ and $\ell>0$ so that the set $\Lambda=\{w=u+iv: |u|<\ell,\ 0\le v<\eps\}$ contains $J=\overline{\{\xi(t):t\in I\}}$ but does not intersect with the slits of $D_t$ for any $t\in I.$  On account of Proposition \ref{P:5.1}(i), we see
that, for any $R_1>\ell$,  $\sup_{t\in I}\sup_{z\in D_t\setminus \Lambda,\; |z|\le R_1}|\H_t(z,\xi(t))|=M_1<\infty$.  Due to the maximum principle for an analytic function, $\H_t(z,\xi(t))$ has the same bound for $z\in \Lambda.$
\qed

We fix $T\in (0,\zeta)$ and set $I=[0,T]$.  By Lemma \ref{L:5.6},
\begin{equation}\label{e:5.29}
M_1:=\sup_{t\in I} \sup_{z\in D_t} 2\pi |z-\xi(t)||\Psi_{\s(t)}(z,\xi(t))|<\infty.
\end{equation}
The next lemma
extends \cite[Lemma 4.13]{L1} from the simply connected domain $\HH$ to multiply connected domains.

\begin{lem}\label{L:5.7}\ For every $t\in I$,  $F_t\subset B(\xi(0),4 R_t)$,
where
$R_t= \sup_{0\le s\le t} |\xi(s)-\xi(0)| \vee \sqrt{ M_1t/2}$.
\end{lem}

\pf\
Fix $t\in I$. For $z\in D$ with $|z-\xi(0)|\geq 4R_t$,
 define $\sigma=\inf\{s: |g_s(z)-z|\ge R_t\}.$  If $s\le t\wedge \sigma,$ then $|g_s(z)-z|<R_t$ and
$$
|\xi(s)-g_s(z)| \ge |(\xi(s)-\xi(0))-(z-\xi(0))|-|g_s(z)-z|
>3R_t-R_t=2R_t.
$$
Hence we have by \eqref{e:5.29}
$$
|\partial_s g_s(z)|=
\left|2\pi \Psi_{\s(s)}(g_s(z),\xi(s))\right|\le \frac{M_1}{|g_s(z)-\xi(s)|
} \le \frac{M_1}{2R_t}.
$$
Consequently,
$ |z- g_s(z)|= | \int_0^s \partial_r g_r (z) dr |\le \frac{M_1}{2R_t}s$
for $s\in [0,  t\wedge \sigma]$.
We claim that $\sigma\geq t$. Suppose otherwise, then by the definition of
$\sigma$, we would have
$R_t=|z- g_\sigma(z)|\le \frac{M_1}{2R_t}\sigma$
and so $\sigma\ge \frac{2}{M_1}R_t^2\geq t$. This contradiction establishes that
 $\sigma\geq t$. So for all $s\in [0,t]$,
 we have $|g_s(z)-z|\leq R_t$ and $|\xi(s)-g_s(z)|\ge 2R_t$.
Thus we have by \eqref{e:5.26} that  $t<t_z$ and $z\in \HH\setminus F_t.$
\qed

\begin{thm}\label{T:5.8}\begin{description}  \item{\rm(i)} The conformal map $g_t(z)$ in {\rm Theorem \ref{T:5.5}} satisfies the
hydrodynamic normalization \eqref{I.3} at infinity.

\item{\rm(ii)}\
The set $F_t$ defined by  \eqref{e:5.27} is an $\HH$-hull; that is,
$F_t$ is relatively closed in $\HH$ and bounded,
and moreover $\HH\setminus F_t$ is simply connected.

\item{\rm (iii)}\ $\{F_t\}$ is strictly increasing in $t$.  It
has the property
\begin{equation}\label{e:5.31}
\bigcap_{\delta>0} \overline{g_t(F_{t+\delta}\setminus F_t)}=\{\xi(t)\} \quad
\hbox{for } t\in [0, \zeta ).
\end{equation}
\end{description}
\end{thm}

\pf\ (i)\ From \eqref{e:5.25}, we have
\[g_t(z)-z=-2\pi \int_0^t\Psi_{\s(s)}(g_s(z), \xi(s))ds.\]
We let $z\to \infty.$  Since the right hand side remains bounded by \eqref{e:5.29},
$g_t(z)\to\infty$ as $z\to\infty.$
Then we can use \eqref{e:5.29} again to see that right hand side converges to $0$ as $z\to\infty,$ yielding the desired conclusion.

\medskip \noindent
(ii)\ It follows from Theorem \ref{T:5.5} and Lemma \ref{L:5.7}
that $F_t$ is relatively closed and bounded.
Were $\HH\setminus F_t$ not simply connected,  $D\setminus F_t$ would be multiply connected of degree at least $N+2$, which is absurd as the conformal image of $D\setminus F_t$ under $g_t$ is the $(N+1)$-ply connected slit domain $D_t$.

\medskip \noindent
(iii)\ Suppose $F_t=F_{t'}$ for some $t'>t\geq 0$.
Then both $g_t$ and $g_{t'}$ are conformal maps from $D\setminus F_t$ onto standard slit domains satisfying the hydrodynamic normalization.
By the uniqueness, we get $g_t(z)=g_{t'}(z),\ z\in D\setminus F_t,$ which is absurd because $\Im g_t(z)$ is strictly decreasing as $t$ increases.

By Lemma \ref{L:5.7} and the fact that $\lim_{t\to 0} R_t=0$,
we have $\cap_{\delta >0} \overline F_\delta = \{ \xi (0)\}$.
So \eqref{e:5.31} holds for $t=0$.
For every $t_0\in (0, \zeta)$, 
 $\{\wh F_{t_0} := g_t(F_{t_0+t}\setminus F_{t_0}); t\in [0, \zeta -t_0)\}$
 is the family of increasing closed sets associated with
  associated with KL-equation \eqref{e:5.25} in Theorem \ref{T:5.5} but with
  $\s(t)$, $\xi (t)$ and $D$ being replaced by $\wh \s (t):= \s (t_0+t)$,
   $\wh \xi (t) := \xi (t_0+t)$ and $\wh D:=D(t_0)$, respectively.
Thus the same argument for $t=0$ above applied to $\{\wh F_\delta; \delta>0 \}$ yields that \eqref{e:5.31}
holds for $t=t_0$.   \qed

In accordance with \cite[p 96]{L1}, we call the property \eqref{e:5.31} the {\it right continuity at} $t$ {\it with limit} $\xi(t).$

\medskip

 We started
this subsection by fixing a pair of functions $(\xi(t), \s(t))$
 satisfying properties {\bf(I), (II)}.
In the rest of this subsection, we shall make a special choice of it, namely,
we fix a solution path $\W_t=(\xi(t),\s(t)), \; t\in [0,\zeta),$ of the SDE \eqref{e:3.39}, \eqref{e:3.40} in Theorem \ref{T:4.2}
 for a given non-negative homogeneous function $\alpha(\s)$ of $\s\in \cS $ with degree $0$ and a given homogeneous function $b(\s)$ of $\s\in \cS $ with degree $-1$ both satisfying the condition {\bf (L)}.

We can now view the associated family $\{g_t(z), t\in [0,t_z)\}$ of conformal maps and the associated growing
$\HH$-hulls
$\{F_t,t\ge 0\}$ constructed in Theorem \ref{T:5.5} and studied in Theorem \ref{T:5.8} as random processes.  Indeed, Proposition \ref{P:4.3} combined with Remark \ref{R:3.12} implies the following scaling properties.

\begin{prop}\label{P:5.9}\
Let $\s\in \cS$, $\xi\in \RR$, $r>0$ and $c\in \RR$.

\begin{description}
\item{\rm(i)}  $\{rg_{t/r^2}(z/r), \;t\ge 0\}$  under $\P_{(\xi/r,\s/r)}$
 has the same distribution as $\{g_t(z),\;t\ge 0\}$ under $\P_{(\xi,\s)}$.

\item{\rm(ii)}
$\{r F_{t/r^2}, \;t\ge 0\}$ under $\P_{(\xi/r,\s/r)}$
 has the same distribution as  $\{F_t,\;t\ge 0\}$  under  $\P_{(\xi,\s)}$.

\item{\rm (iii)}
$\{g_{t}(z-c)+c, \;t\ge 0\}$ and $\{F_{t}-c, \;t\ge 0\}$
under $\P_{(\xi+c, \s +c)}$ have the same distribution as
$\{g_{t}(z), \;t\ge 0\}$ and $\{F_{t}, \;t\ge 0\}$
under $\P_{(\xi, \s)}$, respectively.
\end{description}
\end{prop}

\pf\ (i)\ Let $\W(s)=(\xi(s),\s(s))$ be the solution of the SDE
\eqref{e:3.39}-\eqref{e:3.40} with initial value $(\xi, \s)$.
Note that by Brownian scaling, $\wt \W(s):=r^{-1}\W(r^2s)$ is a solution
to SDE \eqref{e:3.39}-\eqref{e:3.40} driven by Brownian motion
$\wt B_s=r^{-1}B_{r^2 s}$ with initial value $(\xi/r, \s/r)$.
Let $g_t(z)$ be the unique solution of the Komatu-Loewner equation \eqref{e:5.1} driven by $\wt \W$:
 $$
 g_t(z)-z=-2\pi \int_0^t \Psi_{r^{-1}\s(r^2s)}(g_s(z),r^{-1}\xi(r^2s))ds,\ z\in D.
 $$
By Theorem \ref{T:5.5}(i) and Proposition \ref{P:4.3}(i), it suffices to show that $h_t(z):=rg_{t/r^2}(z/r)$ solves the equation \eqref{e:5.1}.

By the homogeneity \eqref{e:3.42} and \eqref{e:3.43},
\begin{eqnarray*}
 \Psi_{r^{-1}\s(r^2s)}(g_s(z), r^{-1}\xi(r^2s))
&=&\Psi_{r^{-1}(\s(r^2s)-\wh\xi(r^2s))}(g_s(z)- r^{-1}\xi(r^2s),0) \\
&=& r\Psi_{\s(r^2s)-\wh\xi(r^2s)}(rg_s(z)- \xi(r^2s),0)
= r\Psi_{\s(r^2s)}(rg_s(z), \xi(r^2s)).
\end{eqnarray*}
and so
$
g_t(z)-z=-2\pi r\int_0^t \Psi_{\s(r^2s)}(rg_s(z),\xi(r^2s))ds
= -\frac{2\pi}{r}\int_0^{r^2t} \Psi_{\s(s)}(rg_{s/r^2}(z),\xi(s))ds.
$\\
Consequently,
$h_t(z)-z=-2\pi \int_0^t \Psi_{\s(s)}(h_s(z),\xi(s))ds.$
That is, $\{h_t (z); t\geq 0\}$ under $\P_{(\xi/r, \s/r)}$ has the same
distribution as $\{g_t (z); t\geq 0\}$ under $\P_{(\xi, \s)}$.

(ii)\ By Theorem \ref{T:5.5}, we have
$F_t=\{z\in D:t_z\le t\}=\{z\in D: \Im g_{s-}(z)=0,\ \hbox{\rm for some}\ s\le t\}.$ \\
Hence the hulls $\{\wh F_t; \geq 0\}$ associated with $\{h_t (z); t\geq 0\}$ is given by
\begin{eqnarray*}
\wh F_t&=&\{z\in D: \Im h_{s-}(z)=0,\ \hbox{\rm for some}\ s\le t\}\\
&=&\{z\in D: \Im g_{(s/r^2)-}(z/r)=0,\ \hbox{\rm for some}\ s\le t\}=rF_{t/r^2}.\end{eqnarray*}
(ii) now follows from (i).

(iii) \ Let $\W(s)=(\xi(s),\s(s))$ be the unique solution of the SDE
\eqref{e:3.39}-\eqref{e:3.40} with initial value $(\xi, \s)$,
and $g_t(z)$ be the unique solution of the Komatu-Loewner equation \eqref{e:5.25} driven by $\W (t)$.
As $b_j (\xi, \s)=b_j (0, \s -\wh \xi)$, $\W(t)+c= (\xi (t)+c, \s (t)+c)$
is the unique solution of the SDE
\eqref{e:3.39}-\eqref{e:3.40} with initial value $(\xi + c, \s +c)$.
In view of second identity in \eqref{e:3.43},
$h_t (z):=g_t (z-c)+c$
is the
unique solution of the Komatu-Loewner equation \eqref{e:5.25} driven by $\W (s)+c$
with $h_0 (z)=z$ for $x\in D+c:=\{w\in \HH: w-c\in D\}$.
This implies the conclusion of (iii).
\qed

See \cite[Proposition 2.1]{RS} for corresponding statements for the case of the simply connected domain $\HH.$

We call the family of random growing hulls $\{F_t; t\ge 0\}$ in Theorem \ref{T:5.8} the {\it stochastic Komatu-Loewner evolution} (SKLE) driven by the solution of the SDE \eqref{e:3.39}-\eqref{e:3.40} with coefficients   $\alpha$ and $b$.
 We designate it as ${\rm SKLE}_{\alpha, b}$.
 Recall that the functions $\alpha$ and $b$ are homogeneous functions on $\cS $ with degree $0$ and $-1$ respectively, and  satisfy the Lipschitz  condition {\bf (L)} in
 \S 4.
In \S 6.1, we shall give a typical example of such a function $b.$

\medskip

Besides the scaling property of ${\rm SKLE}_{\alpha,b}$ demonstrated in Proposition \ref{P:5.9}, we now present its domain Markov property.
Since ${\rm SKLE}_{\alpha,b}$ depends on the initial value $\w=(\xi,\s)\in \RR\times \cS$, we shall denote it also as ${\rm SKLE}_{\w,\alpha,b}$ or ${\rm SKLE}_{\xi,\s,\alpha,b}$.

Let $\W=(W(t), \P_\w)$ be the diffusion process on $\RR\times\cS$ corresponding to the solution of the SDE
\eqref{e:3.39}-\eqref{e:3.40}.
$\W$ satisfies the Markov property with respect to the augmented filtration $\{\gg_t\}$ of the Brownian motion appearing in the SDE.

Let $g_t(z)$ be the unique solution of the ODE \eqref{e:5.25}.  Define $\wt g_s(\wt z)=g_{t+s}\circ g_t^{-1}(\wt z),\ \wt z\in D_t=D(\s(t)).$  Then $\{\wt g_s(\wt z)\}_{s\ge 0}$ is the solution of the KL-equation
\[\partial_t\wt g_s(\wt z)=-2\pi \Psi_{\s(t+s)}(\wt g_s(\wt z), \xi(t+s)),\quad \wt g_0(\wt z)=\wt z.\]
for the driving process $\{W(t+s)=(\xi(t+s),\s(t+s)): s\ge 0\}$ that is the solution of the SDE
\eqref{e:3.39}-\eqref{e:3.40} with initial value $\W(t).$  Consider the associated growing hulls $\{\wt F_s\}$ in $D_t$ for $\wt g_s$ according to \eqref{e:5.27}. Thus
$\dis \{\wt F_s\}_{s\ge 0}\ \text{is the}\ {\rm SKLE}_{W(t),\alpha,b}.$

Take an arbitrary $\wt z\in D_t$ and set $z=g_t^{-1}(\wt z)\in D\setminus F_t.$  Using the Markov property of $\W$, we have for $s\ge 0$
\begin{eqnarray*}
&&\P_{\W(t)}(\wt z\in \wt F_s)=
\P_{\W(t)}(\text{life time of}\ \wt g_\cdot(\wt z)\le s)\\
&=&\P_{\w}(\text{life time of}\ g_{t+\cdot}(z)\le s\;\big|\; \gg_t)
=\P_{\w}(z\in F_{t+s}\setminus F_t\;\big|\; \gg_t)\\
&=&\P_{\w}(\wt z\in g_t(F_{t+s}\setminus F_t)\;\big|\; \gg_t),\quad \w \in \RR\times \cS.
\end{eqnarray*}
By Theorem 5.5, the set-valued random variable $F_t$ is $\gg_t$-adapted.  Denote by $\gg_t^0$ the sub-$\sigma$-field of $\gg_t$ generated by $\{F_u; u\le t\}.$  In view of Theorem 5.5 and Theorem 5.8, $\W(t)=(\xi(t), \s(t))$ is $\gg_t^0$-adapted so that
\begin{equation}\label{5.25}
\P_{\w}(\wt z\in g_t(F_{t+s}\setminus F_t)\;\big|\; \gg_t^0)=\P_{\W(t)}(\wt z\in \wt F_s),\quad \w\in \RR\times \cS.
\end{equation}
This can be rephrased as follows:

\begin{prop}\label{P:DMP}\
For every $\w\in \RR\times \cS$,
$\P_\w$-a.s. the conditional law of $\{g_t(F_{t+s}\setminus F_t)\}_{s\ge 0}$
given $\gg_t^0 $ has the same distribution as that of $ {\rm SKLE}_{W(t),\alpha,b}$.
\end{prop}

It will be shown in Theorem \ref{T:5.10} below that the half-plane capacity
of ${\rm SKLE}_{{\bf w},\alpha,b}$ is $2t$.

For $\wh D =D\setminus F\in \wh \dd$, where $D\in \dd$ and $F\subset D$ is an $\HH$-hull,
let $\Omega (\wh D)$ denote the collection
of families of increasing
bounded closed subsets $ {\bf F}=\{  {\bf F} (t); t\geq 0\}$ of $\wh D$
such that each $F\cup {\bf F}(t)$ is an $\HH$-hull.
For $D\in \dd$, we
introduce a filtration $\{\gg_t( D); t\geq 0\}$
on $\O( D)$ by
\[
\gg_t(D):= \sigma\{  {\bf F} (s):\;0\le s\le t\} , \qquad
\gg(D):=\sigma\{ {\bf F} (s): s\ge 0\}.\]
For $\wh D \in \wh \dd,$ we then
introduce a $\sigma$-field
$ \gg (\wh D)$
on $\O(\wh D)$ by
$\dis \gg(\wh D)= \Phi^{-1}\gg\left(\Phi(\wh D)\right),$
using the canonical map $\Phi$ from $\wh D$ to $\Phi(\wh D)\in \dd.$
For $D\in \dd$ and $t\geq 0$, define the shift operator
$\theta_t:  \O( D) \mapsto \O( D \setminus {\bf F}(t) )$ by
 \begin{equation}\label{e:5.26b}
(\theta_t  {\bf F})(s) = {\bf F} (t+s)\setminus {\bf F}(t) \qquad \hbox{for }
 s \geq 0 .
\end{equation}

For $D\in \dd$ and $z\in \partial  \HH$, we  use $\P_{ D, z}$
to denote the induced probability measure on $\Omega ( D)$ by
$\P_{\bf w}$, where ${\bf w}=(z, \s (D))$.
Observe that
by Theorem \ref{T:5.8}, $\{g_t (z); t\geq 0\}$ driven by the solution $\W_t=(\xi(t),\s(t))$
 of the SDE \eqref{e:3.39}-\eqref{e:3.40} with initial condition $\W_0 ={\bf w}$
 is the unique conformal map from $D\setminus F_t$
to a standard slit domain for each fixed $t\geq 0$ satisfying the hydrodynamic normalization at infinity, where $\{F_t; t\geq 0\}$ are the associated
${\rm SKLE}_{{\bf w}, a, b}$-hulls.
Thus the probability measures $\P_{\bf w}$ and $\P_{D, z}$ are in one-to-one correspondence.

\begin{thm}\label{T:DMP} The probability measures $\{\P_{ D, z}; D\in \dd, z\in \partial \HH\}$ enjoy the following properties.

\begin{description}

\item{\rm (i)} For any $D\in \dd$ and $z\in\partial \HH$,
$$
\P_{D, z} \left(  \cap_{t>0} {\bf F} (t)=\{z\}  \hbox{ and the half-plane capacity of }
{\bf F} (t) \hbox{ is }  2t  \hbox{ for every } t \geq 0 \right) =1.
$$
Let $\wh g_t(z)$ be the canonical map on $  D\setminus {\bf F}(t)$
and $\wt \s (t):=\s (D_t)$, where $D_t:=\wh g_t (D\setminus {\bf F}(t))\in
\dd$. Then
$$
\P_{ D, z}\left( \bigcap_{\delta >0} \overline{\wh  g_t ({\bf F}(t+\delta)
 \setminus {\bf F}(t))} = \{\wt \xi (t)\} \subset \partial \HH
 \hbox{ for every } t\geq 0\right) =1.
$$
Moreover,  $(\wt \xi (t), \wt \s (t))$ has the same distribution
  as the unique solution $(\xi (t), \s (t))$ of \eqref{e:3.39}-\eqref{e:3.40} with initial condition
 $(  \xi (0),  \s (0)) =(z,  \s (\D))$.

\item{\rm (ii)}  \rm{(}Domain Markov property\rm{)}: For each
 $t\ge 0$,
\begin{equation}\label{e:5.27a}
\P_{D,z}\left( \theta_t^{-1}\Lambda\big| \gg_t(D)\right)=\P_{ D_t,  \wt \xi (t)}( \wh g_t (\Lambda )) \quad
\hbox{for every } \Lambda\in \gg( D \setminus {\bf F}(t)).
\end{equation}

 \item{\rm (iii)} {\rm(}Invariance under linear conformal map{\rm)}: for any $D\in  \dd$ and any linear conformal map $f$ from $D$ onto $f(D)\in \dd,$
\begin{equation}\label{e:5.28a}
\P_{f(D),f(z)}= \P_{D,z}  \circ f^{-1}
\qquad  \hbox{for every } z\in \partial \HH.
\end{equation}
\end{description}
\end{thm}

\pf  (i)  follows immediately from Theorem \ref{T:5.8}.

(ii) Consider a generic event
$\wt \Lambda=\{\wt\F\in \Omega(D_t):\wt z\in \wt\F(s)\}\in \gg(D_t)$
for $\wt z\in D_t,\ s\ge 0$.
 Such sets generate the $\sigma$-field $\gg (D_t)$. Define
$\Lambda=\wh g_t^{-1}(\wt \Lambda)$.
 Clearly,
$\Lambda\in \gg(D\setminus \F(t))$ and $\wt \Lambda=\wh g_t(\Lambda)$.
Observe that
\begin{eqnarray*}
\theta_t^{-1}\Lambda
&=& \{ \F \in \Omega(D): \{ \F(u+t)\setminus  \F(t)\}_{u\ge 0}\in \Lambda\}\\
&=& \big \{ \F \in \Omega(D): \{\wh g_t( \F(u+t)\setminus  \F(t))\}_{u\ge 0}\in \wh g_t (\Lambda)=\wt\Lambda \big\} \\
&=& \{ \F\in \Omega(D): \wt z\in \wh g_t( \F(s+t)\setminus \F(t))\}.
\end{eqnarray*}
Now \eqref{e:5.27a} follows from Proposition \eqref{5.25}
and thus (ii) is established as such $\Lambda = \wh g_t^{-1} (\wt \Lambda)$ generates $\gg (D\setminus
{\bf F}(t))$.

(iii) Let $f(z)= c_1 z + c_2$, $c_1>0$, $c_2\in \RR$,
 be a linear conformal map from $D$ to $f(D)\in \dd$.
Clearly, $f^{-1}(z)= (z-c_2)/c_1$.
It follows from Proposition \ref{P:5.9} that for $\xi \in \partial \HH$ and
$\s \in \cS$,
$\{ c_1^{-1} ( g_{c_1^2t}(c_1z-c_2)+c_2); t\geq 0\}$ under $\P_{(f(\xi), f(\s))}$ has the same distribution as
$\{g_t (z); t\geq 0\}$ under $\P_{(\xi, \s)}$. Consequently,
$\{ F_t, t\geq 0\}$ under $\P_{f(D), f(\xi)}$  has the same distribution as
$\{c_1 F_{c_1^{-2}t}-c_2; t\geq 0\}$ under $\P_{D, \xi}$.
That is, $\P_{f(D),f(z)}= \P_{D,z}  \circ f^{-1}$,
under the $2t$-half-plane capacity parametrization.
\qed

We remark that
 the shift operator $\theta_t$ in \eqref{e:5.26b} is a natural extension of
$\dot \theta_t$ in \eqref{e:3.4a}, and
 the identity \eqref{e:5.27a} is analogous to \eqref{e:3.9} in Section 
\ref{S:3.2}.

\subsection{Half-plane capacity for SKLE}\label{S:5.2}

We return to the general setting made in the beginning of \S \ref{S:5.1},
 and consider the conformal maps $\{g_t(z)\}$ and
$\HH$-hulls
$\{F_t\}$ in Theorem \ref{T:5.5}.  Let  $a_t$ be the half-plane capacity of $F_t$; that is,  $a_t:=\lim_{z\to\infty}z(g_t(z)-z)$.

\begin{thm}\label{T:5.10}\ It holds that
$a_t=2t$ for every $t\ge 0$.
\end{thm}

This theorem follows immediately from the following proposition,
which compared with the equation \eqref{e:5.25}  implies that $a_t$ is differentiable and $\frac{da_t}{dt}=2$.

\begin{prop}\label{P:5.11}\ $a_0=0,$ $a_t$ is strictly increasing and right continuous.  $g_t(z)$ is right differentiable in $a_t$ and
\begin{equation}\label{e:5.33}
\frac{\partial^+g_t(z)}{da_t} = -\pi
\Psi_{\s (t)}
(g_t(z),\xi(t)),\quad g_0(z)=z\in D,\quad t\in [0,t_z).
\end{equation}
Here $\frac{\partial^+g_t(z)}{da_t}$ is the right
derivative of $g_t (z)$ with respect to $a_t$.
\end{prop}

To prove this, we make arguments parallel to \cite[\S 6.2, \S 6.3, \S 8]{CFR}.   Note however that, while $F_t$ is a portion of a given Jordan arc in \cite{CFR}, $F_t$ is now defined by \eqref{e:5.27} for the solution $g_t(z)$ of the equation \eqref{e:5.25}
for a given continuous function $(\xi(t), \s(t))$ satisfying the property {\bf(II)}.

Fix $t_0>0$ and, for $0\le s<t\le t_0$, set
$g_{t,s}=g_s\circ g_t^{-1}$, which
is a  conformal map from $D_t$ onto
$D_s\setminus g_s(F_t\setminus F_s)$.
 Its inverse map $g_{t,s}^{-1}$ is a conformal map from
$D_s\setminus g_s(F_t\setminus F_s)$ onto the standard slit domain $D_t$ and satisfying a hydrodynamic normalization.  Therefore, in view of the proof of \cite[Theorem 7.2]{CFR}, we can draw the following conclusion: let $\ell_{t,s}$
 be the set of all limiting points of $g_{t,s}^{-1}\circ g_s(z)=g_t(z)$ as $z$ approaches to $F_t\setminus F_s,$ then $\ell_{t,s}$ is a compact subset of $\partial\HH$ and $g_{t,s}^{-1}$ sends $\partial\HH\setminus\overline{g_s(F_t\setminus F_s)}$ into $\partial\HH$ homeomorphically.

 Let  $\Lambda=\{x+iy:a<x<b,\ 0<y<c\}$ be a finite rectangle so that $\ell_{t,s}\subset \{x+i0: a<x<b\}$ and $\Lambda \subset \cap_{0\le t\le t_0} D_t$.
 Then $\Im g_{t,s}(z)$ is uniformly bounded in $z\in \Lambda$ and, by the Fatou theorem (cf. \cite{GM}), it admits finite limit
\begin{equation}\label{e:5.36}
\Im g_{t,s}(x+i0+)=\lim_{y\downarrow 0}\Im g_{t,s}(x+iy)\quad\hbox{\rm for a.e.}\quad x\in (a,b).
\end{equation}
In exactly the same way as the proof of \cite[Lemma 6.2, Theorem 6.4]{CFR},
we get the following.

\begin{lem}\label{L:5.13}\ For $0\leq s <t\leq t_0$,
$\displaystyle a_t-a_s= \pi^{-1}
\int_{\ell_{t,s}} \Im g_{t,s}(x+i0+)dx
$ and
$$
g_s(z)-g_t(z)=
\int_{\ell_{t,s}} \Psi_{\s(t)}(g_t(z),x)\Im g_{t,s}(x+i0+)dx,\quad z\in D\setminus F_t.
$$
\end{lem}

By the Schwarz reflection, we can extend $g_{t,s}^{-1}$ to a conformal map on
\[ D_s\cup\Pi D_s\cup \partial\HH \setminus (\overline{g_s(F_t\setminus F_s)}
\cup \Pi g_s(F_t\setminus F_s)).\]

\begin{lem}\label{L:5.14}\begin{description}
 \item{\rm(i)}\ For any compact subset $V$ of $D_s\cup
 \partial\HH\setminus \{\xi(s)\}$,
$\lim_{t\downarrow s} g_{t,s}^{-1}(z)=z$ uniformly in $z\in V\cup \Pi V.$

\item{\rm(ii)}\ $a_t$ is right continuous in $t$.

\item{\rm (iii)}\ $a_t$ is non-negative and strictly increasing in $t$.
\end{description}
\end{lem}

\pf\ (i) \ Without loss of generality, we may assume   $s=0$ and so $g_{t,s}^{-1}=g_t.$  Let $V$ be any relatively compact open subset of $D\cup(\partial\HH\setminus \{\xi(0)\}).$
In Theorem \ref{T:5.5}, we considered the family of solution curves
 $\{(g_t(z),\; 0\le t<t_z): z\in D\}$ of \eqref{e:5.25} parametrized by the initial position $z=g_0(z)\in D.$  We add to this family the solution curve $(g_t(z),\; 0\le t<t_z)$ of \eqref{e:5.25} with initial position $z=g_0(z)\in \partial\HH\setminus \{\xi(0)\}$ satisfying $g_t(z)\in \partial\HH,\; 0\le t<t_z,$ where
\[ t_z=\sup\{t\in [0,\zeta): \inf_{s\in [0,t]}|g_s(z)-\xi(s)|>0\}.\]
By Proposition \ref{P:5.1} (vi) and Lemma \ref{L:5.2} (ii), such a solution exists uniquely and takes values in $\partial \HH$.
  Define $F_t(\partial\HH)=\{z\in \partial\HH\setminus\{\xi(0)\}: t_z\le t\},\ t>0.$  By a general theorem on ODE cited in the proof of Theorem \ref{T:5.5} already, $g_t(z)$ is jointly continuous on ${\cal G}=\{(t,z): z\in D\cup(\partial\HH\setminus\{\xi(0)\}),\ t\in [0,t_z)\}.$

For the set $V$ as above, Theorem \ref{T:5.8} (iii) implies that
there exists $\delta>0$ such that
 $\overline{F_\delta} \cup F_\delta(\partial\HH)$ is disjoint from $\overline V.$
So $[0,\delta]\times \overline V$ is a compact subset of ${\cal G}.$
Hence
$\sup_{t\in [0,\delta], z\in \overline V\cup\Pi \overline V}|g_t(z)|=\sup_{t\in [0,\delta], z\in \overline V}|g_t(z)|$ is finite by the continuity of $g_t(z)$ mentioned above, and accordingly $\{g_t(z): 0\le t\le \delta\}$ is a normal family of analytic functions on $V\cup \Pi V$. This implies that
$\lim_{t\downarrow 0} g_t(z)=z$ uniformly in $z\in V\cup \Pi V$.

 \medskip
(ii)\ This follows from (i) as in the proof of \cite[Theorem 8.4]{CFR}.

\medskip
(iii)\ Choose $R>0$ so large that $F_t\cup K\subset \{|z|<R\}.$ By \eqref{e:7.19} below, we then have
\[a_t=\frac{2R}{\pi}\int_0^\pi h_t(Re^{i\theta})\sin \theta d\theta\quad{\rm for}\quad h_t(z)=\E_z^*\left[\Im Z_{\sigma_{F_t}}^*; \sigma_{F_t}<\infty\right].\]
Here $Z^*=(Z^*_t, \zeta^*, \P_z^*)$ is BMD on $D\cup\{c_1^*,\cdots, c_N^*\}$.
Since, by Theorem \ref{T:5.8},  $F_t$ is
strictly increasing in $t$ and $\HH\setminus F_t$ is simply connected,
$F_t$ is non-polar for the planar Brownian motion and consequently for the absorbing Brownian motion on $D.$  Hence the above expression implies that $a_t>0$ for $t>0.$   As $\{F_t\}$ is strictly increasing, so is $\{a_t\}$ by its additivity under the composite map.
\qed

\noindent
{\bf Proof of Proposition \ref{P:5.11}. }
We now know from
Lemma \ref{L:5.14} that $a_t$ is strictly increasing and right continuous.
For any $\eps_0>0$ with $B (\xi(s), \eps_0)\cap \HH\subset D_s$, there exists $\delta>0$ so that
\[g_s(F_t\setminus F_s)\cup \Pi g_s(F_t\setminus F_s)\subset B(\xi (s), \eps_0)
\quad \hbox{for any } t\in (s,s+\delta)
\]
 by virtue of Theorem \ref{T:5.8}(iii).
 In particular,
 $\ell_{t,s}$ is in the interior of the region bounded by the Jordan curve
 $g_{t,s}^{-1} (\partial B(\xi (s), \eps_0))$.
By Lemma \ref{L:5.14}, we have for sufficiently small $\delta>0$,
\[
|g_{t,s}^{-1}(z)-z|<\eps_0,\ \hbox {\rm for any}\ z\in \partial B(\xi (s), \eps_0) \hbox{ and for any }  t\in (s,s+\delta).
\]
In particular,
the diameter of $g_{t,s}^{-1} (\partial B(\xi (s), \eps_0))$ is less than
$3 \eps_0$.
  Therefore, we get for any $x\in \ell_{t,s}$
\begin{equation}\label{e:5.39}
|\xi(s)-x|\le |\xi(s)-z|+|z-g_{t,s}^{-1}(z)|+|g_{t,s}^{-1}(z)-x|<5\eps_0,
\end{equation}
by taking any $z\in g_{t,s}^{-1} (\partial B(\xi (s), \eps_0))$.
On the other hand, from the Lipschitz continuity of $\Psi$ and the continuity of $\s(t)$,
 we can conclude that $\Psi_{\s(t)}(z,x)$ is jointly continuous in $(t,z,x)$
 as in the proof of \cite[Theorem 9.8]{CFR}.  Fix $z\in D.$  Since $g_t(z)$ is continuous in $t$, $\Psi_{\s(t)}(g_t(z), x)$ is continuous in $t>0$ and $x \in \HH.$  Therefore, for any $\eps>0,$ there exist $\delta>0$ and $\eps_0>0$ such that
\begin{equation}\label{e:5.40}
|\Psi_{s(t)}(g_t(z),x)-\Psi_{\s(s)}(g_s(z), \xi(s))|<\eps
\end{equation}
for any $t\in (s,s+\delta)$ and for any $x\in \partial\HH$ with $|x-\xi(s)|<5\eps_0.$
It now follows from Lemma \ref{L:5.13}, \eqref{e:5.39} and \eqref{e:5.40} that, there exists $\delta>0$ such that, for any $t\in (s,s+\delta),$
\[\left|\frac{g_t(z)-g_s(z)}{a_t-a_s}+\pi\Psi_{\s(s)}(g_s(z),\xi(s))\right|<\eps .
\]
This proves the Proposition. \qed

\section{Locality of SKLE}\label{S:6}

\subsection{BMD domain constant $b_{\rm BMD}$}\label{S:6.1}

For each standard slit domain $D\in \dd$,
 let $\Psi(z,\xi)=\Psi_D(z,\xi)$, $z\in D$, $\xi\in \partial\HH$,
 be the BMD-complex Poisson kernel of $D$,  and define
\begin{equation}\label{e:6.1}
b_{\rm BMD}(\xi;D)=2\pi \lim_{z\to\xi} \left(\Psi_D(z,\xi)+\frac{1}{\pi}\frac{1}{z-\xi}\right),\quad \xi\in \RR.
\end{equation}
Since
$\Psi^\HH(z,\xi)=-\frac{1}{\pi}\frac{1}{z-\xi}$
is the complex Poisson kernel for the ABM on $\HH,$ $b_{\rm BMD}(\xi;D)$ indicates a discrepancy of the slit domain $D$ from $\HH$ relative to BMD.
It follows from Lemma \ref{L:5.6}(ii) that $b_{\rm BMD}(\xi;D)$ is well-defined
by \eqref{e:6.1} as a finite real number.
Sometimes we also write $b_{\rm BMD}(\xi;D)$ as $b_{\rm BMD}(\xi,\s)$ in terms of the slits $\s=\s(D) $ of $D$.
We set $b_{\rm BMD}(\s)=b_{\rm BMD}(0,\s)$ and call it the $BMD$ {\it domain constant}
of $D=D (\s)$.

\begin{lem}\label{L:6.1}\begin{description}
 \item{\rm(i)}\ $b_{\rm BMD}(\s)$, $\s\in \cS$,  is a homogeneous function of degree $-1$ on $\cS$.

\item {\rm(ii)}\ $\dis b_{\rm BMD}(\xi,\s)=b_{\rm BMD}(\s-\wh \xi)$ for
$\s\in \cS$ and $\xi\in \RR$.

\item {\rm(iii)}\ $b_{\rm BMD}(\s)$ satisfies the Lipschitz continuity condition {\bf(L)} {\rm(}see {\rm \eqref{e:4.1})}.
\end{description}
\end{lem}

\pf\ (i)\
By \eqref{e:3.42} in Remark \ref{R:3.12},  for any $\s \in \cS$ and $c>0$,
$$
b_{\rm BMD}(c\s)
= 2\pi\lim_{z\to {\bf 0}} \left( \Psi_{c\s}(cz,{\bf 0})+ (c \pi z)^{-1} \right)
=  \frac{2\pi}{c} \lim_{z\to {\bf 0}} \left( \Psi_{\s}(z,{\bf 0})+ (\pi z )^{-1} \right) =c^{-1} b_{\rm BMD}(\s).
$$
\noindent
(ii)\ By \eqref{e:3.43}, we have for any $\eta\in \RR$
\[2\pi \left( \Psi_\s(z,\xi)+\frac{1}{\pi}\frac{1}{z-\xi} \right)
=2\pi \left( \Psi_{\s+\wh\eta}(z+\eta,,\xi+\eta)+\frac{1}{\pi}\frac{1}{(z+\eta)-(\xi+\eta)}\right).
\]
Taking $z\to {\bf 0}$ yields  $b_{\rm BMD}(\xi,\s)=b_{\rm BMD}(\xi+\eta,\s+\wh\eta)$.

\noindent
(iii) For $\s_1, \s_2 \in \cS$, $ b_{\rm BMD}(\s_1) - b_{\rm BMD}(\s_2)
= 2\pi \lim_{z\to {\bf 0}}\left( \Psi_{D(\s_1)}(z, {\bf 0})- \Psi_{D(\s_2)}(z, {\bf 0} )\right)$.
The Lipschitz continuity of $b_{\rm BMD}(\s)$ in $\s \in \cS$ follows from the Lipschitz continuity of $\Psi_{D}$ in $D\in \dd$ established in \cite[Theorem 9.1]{CFR}.
\qed

 \subsection{Generalized Komatu-Loewner equation for image hulls}\label{S:6.2}

In the rest of this paper, we   make a special choice of the driving process $(\xi(t),\s(t))$ as in the last part of \S 5.1:
let $\W_t=(\xi(t),\s(t))$ be the solution
of the SDE \eqref{e:3.39}-\eqref{e:3.40} in Theorem \ref{T:4.2}
 for a given non-negative homogeneous function $\alpha(\s)$ of $\s\in \cS $ with degree $0$ and a given homogeneous function $b(\s)$ of $\s\in \cS $ with degree $-1$,  both satisfying the condition {\bf (L)}.

We shall use the term ``canonical map" introduced in the second paragraph of
\S \ref{S:3.1}.  Let $\{g_t(z)\}$ and $\{F_t\}$ be the family of the random conformal maps and the random growing hulls in
Theorem \ref{T:5.5}.
Recall that $\{F_t\}$ is called the SKLE driven by the solution of the SDE
\eqref{e:3.39}-\eqref{e:3.40} with coefficients determined by $\alpha$ and $b$, and is designated as ${\rm SKLE}_{\alpha,b}$.  For each $t>0$,   $g_t$ is the canonical map from
$D\setminus F_t$ onto $D_t=D(\s(t))$ where $D$ denotes $D(\s(0)).$

To formulate a locality property of SKLE, take any $\HH$-hull $A\subset D$ and
define
$$ \tau_A=\inf\{t>0: \overline F_t\cap \overline A\neq \emptyset\}.
$$
In what follows, we only consider those parameter $t$ with
$\tau_A$.

\begin{figure}[h]
	\begin{center}
    \vspace{-1em}	\includegraphics[scale=.45]{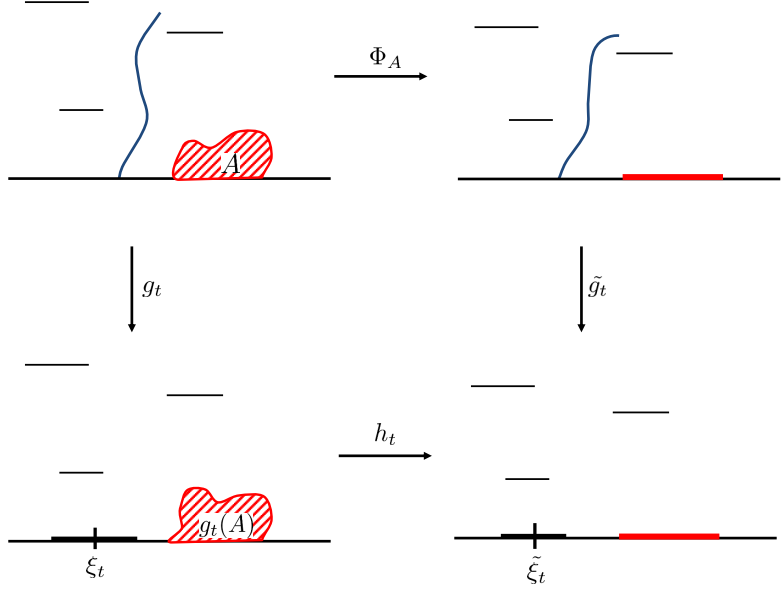}
    \vspace{-2.0em} %adjust the distance between the caption and the graph
    \caption{Conformal mappings $\Phi_A$ and $h_t$}\label{fig:5}
    \vspace{-1em}
	\end{center}
	\end{figure}

  Let $\Phi_A$ be a canonical map from $D\setminus A$ onto $\wt D\in \dd$ and define $\wt F_t=\Phi_A(F_t) $.   Let $\wt g_t $ be the canonical map from $\wt D\setminus \wt F_t$ onto $\wt D_t\in \dd$ and
  $\wt a_t$  the half-plane capacity of $\wt g_t$, that is
  $\wt a_t=\lim_{ z\to\infty}  z(\wt g_t(  z)-  z).$
$\wt a_t$ will be also denoted by $\wt a(t).$
Along with the canonical maps $g_t,\ \Phi_A$ and $\wt g_t,$
we consider
the canonical map $h_t$ from $D_t\setminus g_t(A)$.
Then
\begin{equation}\label{e:6.2}
\wt g_t\circ \Phi_A= h_t\circ g_t
\end{equation}
because both of them are canonical maps from $D\setminus (F_t\cup A)$.
See Figure \ref{fig:5}.
The union of the slits in domains $\wt D$ and $\wt D_s$ are denoted by
$\wt K=\bigcup_{j=1}^N \wt C_j$ and $\wt K(s)=\bigcup_{j=1}^N \wt C_j(s)$, respectively.  Denote by $\wt A$ the set of all limiting points of $\Phi_A(z)$ as $z$ approaches to $A.$

Define
\begin{equation}\label{e:6.3}
\wt \xi(t)=h_t(\xi(t)).
\end{equation}
We further denote by $\wt \Psi_t(z,x),\;z\in \wt D_t,\; x\in \partial\HH,$ the BMD-complex Poisson kernel of $\wt D_t.$

In this subsection, we aim at proving Proposition \ref{P:6.2} for $\{\wt F_t\}$ stated below that is analogous to Proposition \ref{P:5.11} formulated for $\{F_t\}.$
To this end, we prepare three lemmas and a proposition.

\begin{lem}\label{L:6.3}\ $\{\wt F_t\}$ is strictly increasing in $t.$  It is right continuous at $t$ with limit $\wt\xi(t)$
in the following sense:
\begin{equation}\label{e:6.5}
\bigcap_{\delta>0} \overline{
 \wt g_t (\wt F_{t+\delta} \setminus \wt F_t) }=\{\wt\xi(t)\}.
\end{equation}
\end{lem}

\pf\  The first statement follows from the corresponding statement in Theorem \ref{T:5.8}.  The second one follows from \eqref{e:5.31}, \eqref{e:6.1} and \eqref{e:6.2} as
\[
\bigcap_{t>s}\overline{\wt g_s(\wt F_t\setminus \wt F_s)}=\bigcap_{t>s} \overline{ h_sg_s\Phi_A^{-1}(\wt F_t\setminus \wt F_s)}
= \bigcap_{t>s}\overline{h_sg_s(F_t\setminus F_s)}=h_s(\xi(s))=\wt \xi(s).
\]
\qed

For $0\le s<t<\tau_A$,
set $\wt g_{t,s}=\wt g_s\circ \wt g_t^{-1}$.
Denote by  $\wt \ell_{t,s}$ the set of all limiting points of $\wt g_{t,s}^{-1}\circ\wt g_s(\wt z)=g_t(\wt z)$ as $\wt z$ approaches to $\wt F_t\setminus \wt F_s.$

\begin{lem}\label{L:6.4}\begin{description} \item {\rm(i)}\ $\wt \ell_{t,s}$ is a compact subset of $\partial\HH$ and
\begin{equation}\label{e:6.7}
\wt a_t-\wt a_s=\frac{1}{\pi}
\int_{\wt\ell_{t,s}} \Im \wt g_{t,s}(x+i0+)dx,
\end{equation}
\begin{equation}\label{e:6.8}
\wt g_s( z)-\wt g_t ( z)=
\int_{\wt \ell_{t,s}} \wt \Psi_t(\wt g_t(  z),x)\Im \wt g_{t,s}(x+i0+)dx,\quad
 z\in \wt D\cup \partial_p\wt K\setminus \wt F_t,
\end{equation}
where $\Im \wt g_{t,s}(x+i0)$ is the Fatou boundary limit existing {\rm a.e.} on $\wt\ell_{t,s}.$

\item{\rm(ii)}\ $\wt a_t>0$ and $\wt a_t$ is strictly increasing.

\item{\rm(iii)}\ For each $t>0$ and
$ z\in \wt D \setminus \wt F_t,$\quad
$\sup_{0\le s\le t} |\wt g_s(  z)|<\infty.$
\end{description}
\end{lem}

\pf\ (i) \ This can be shown in the same way as that for Lemma \ref{L:5.13}.  The identity \eqref{e:6.8} can be obtained first for
$z\in \wt D\setminus \wt F_t$ and then extended to
$z\in \wt D\cup\partial_p\wt K\setminus \wt F_t.$

(ii)\ This can be proved exactly in the same way as that for Lemma \ref{L:5.14} (iii) by the probabilistic expression \eqref{e:7.19} for $\wt a_t.$

(iii)\  For $0\leq s\leq t$,   \eqref{e:6.7} and \eqref{e:6.8} imply that
$\dis  |\wt g_s( z)| \le |\wt g_t(  z)|+\pi \sup_{x\in \wt \ell_{t,0}}
|\wt \Psi_t(\wt g_t(z),x)| \, \wt a_t.$   \qed

We next present a probabilistic representation of $\Im\wt g_t(z)$ which enables
us to derive the joint continuity of $\Im\wt g_t(z)$
with a uniform bound from those of  $\Im g_t(z)$.

For $D=\HH \setminus \bigcup_{j=1}^N C_j$,
we consider Jordan curves $\eta_j$ surrounding $C_j$ that are mutually disjoint and disjoint from $F_t\cup A\cup \partial \HH.$   Denote by $Z^{D,*}=(Z_\cdot^{D,*},\P_z^{D,*})$ the BMD on $D\cup\{c_1^*,\cdots, c_N^*\}$ obtained from the absorbing Brownian motion (ABM)
$Z^\HH=(Z_\cdot^\HH, \P_z^\HH)$ on $\HH$ by rendering each slit $C_j$ into a single point $c_j^*$, and set
$K=\bigcup_{j=1}^N C_j$.
Notice that $Z^{D,*}$ was denoted as $Z^{\HH,*}$ in \cite{CFR}.
The notation $Z^{D,*}$ is more convenient for later discussions.  Define a measure $\nu_j$ on $\eta_j$ by
\begin{equation}\label{e:6.9}
\nu_j(B)=\P_{c_j^*}^{D,*}(Z_{\sigma_{\eta_j}}^{D,*}\in \Gamma),\quad \Gamma\in {\cal B}(\eta_j),\quad 1\le j\le N.
\end{equation}

\begin{prop}\label{P:6.5} \begin{description}
\item{\rm(i)}\ Define for $z\in D\setminus ( F_t\cup A)$,
\begin{equation}\label{e:6.10}
q_t(z) =  \Im g_t(z)-\sum_{j=1}^N\kappa_j(t)\;\P_z^\HH\left(Z_{\sigma_K}^\HH\in C_j;\; \sigma_K<\sigma_A\right)
-  \E_z^\HH\left[\Im g_t(Z_{\sigma_A}^\HH);\; \sigma_A<\sigma_K\right] ,
\end{equation}
where $\kappa_j(t)$ is the $y$-coordinate of the $j$th slit of $D_t.$
It  holds for
$  z\in \wt D\setminus \wt F_t$ that
\begin{equation}\label{e:6.11}
 \Im \wt g_t( z)= q_t(\Phi_A^{-1}  z)
+
\sum_{i=1}^N\P_{\Phi_A^{-1}  z}^\HH\left( Z_{\sigma_K}^\HH\in C_i;\; \sigma_K<\sigma_A\right)\;\sum_{j=1}^N \gamma_{ij}\int_{\eta_j}\; q_t(z)\nu_j(dz),
\end{equation}
for some positive constants $\gamma_{ij}$, $1\le i, j\le N$, independent of $t$.

\item{\rm(ii)}
\ For each $T\in (0,\tau_A),$ the function $\Im\wt g_t(z)$ is extended to be jointly continuous in $(t,z)\in [0,T]\times \overline \HH\setminus \overline{\wt F_T}\setminus\overline{\wt A}$ and has a bound, for some constant $\gamma>0,$
\begin{equation}\label{e:6.12}
0 \le  \Im\wt g_t(z)\le \Im \Phi_A^{-1} z  + \gamma\quad{\rm for}\quad t\in [0,T], \ z\in \overline\HH\setminus\overline{\wt F_T}\setminus \overline{\wt A}.
\end{equation}
\end{description}
\end{prop}

\pf (i)\  For $D_t=\HH\setminus K(t),\ K(t)=\bigcup_{j=1}^N C_j(t),$
$g_t(\eta_j)$ are Jordan curves surrounding $C_j(t)$ that are mutually disjoint and disjoint from $ g_t(A)\cup \partial \HH.$   Let $Z^{D_t,*}=(Z_\cdot^{D_t,*},\P_z^{D_t,*})$ be the BMD on $D_t\cup\{c_1^*(t),\cdots, c_N^*(t)\}$ obtained from the ABM $Z^\HH=(Z_\cdot^\HH, \P_z^\HH)$ on $\HH$ by rendering each slit $C_j(t)$ into a single point $c_j^*(t)$.
Analogously to \eqref{e:6.9},
define a measure $\nu_j^t$ on $g_t(\eta_j)$ by
\begin{equation}\label{e:6.13}
\nu_j^t(\Gamma)=\P_{c_j^*(t)}^{D_t,*} \left(Z_{\sigma_{g_t(\eta_j)}}^{D_t,*}\in \Gamma \right),\quad \Gamma\in {\cal B}(g_t(\eta_j)),\quad 1\le j\le N.
\end{equation}
Owing to the conformal invariance of the BMD
(see \cite[Remark 7.8.2]{CF1}),
we have
\begin{equation}\label{e:6.14}
\nu_j^t(g_t(\Gamma))=\nu_j(\Gamma)
\quad \hbox{for any } \Gamma\in {\cal B}(\eta_j).
\end{equation}
Applying
\cite[Theorem 7.1]{CFR} to the canonical map $h_t$ from $D_t\setminus g_t(A)$ with $g_t(\eta_j)$ in place of $\eta_j,\ 1\le j\le N$,
we get
\begin{equation}\label{e:6.15}
\begin{cases}
\Im h_t(z)=v_t(z)
+ \sum_{j=1}^N f_j(t,z) v_t^*(c_j^*(t))\\
v_t(z)=\Im z-\E_z^\HH\left[\Im Z_{\sigma_{K_t\cup g_t(A)}}^\HH;
\sigma_{K_t\cup g_t(A)}<\infty\right]\\
v_t^*(c_j^*(t))=\sum_{k=1}^N \frac{M_{jk}(t)}{1-R_k^*(t)}\int_{g_t(\eta_k)}v_t(z)\nu_k^t(dz).
\end{cases}
\end{equation}
Here
\begin{equation}\label{e:6.16}
f_j(t,z)=\P_z^\HH \left( Z_{\sigma_{K_t}}^\HH\in C_j(t); \sigma_{K_t}<\sigma_{g_t(A)} \right)
,\quad R_k^*(t)= \int_{g_t (\eta_k) } f_k (t, z) \nu^t_j (dz) ,
\end{equation}
and $M_{jk}(t)$ is the entry of $\sum_{n=0}^\infty Q^*(t)^n$ for the matrix $Q^*(t)$ with zero diagonal entry and off-diagonal entry given by
\begin{equation}\label{e:6.17}
q_{ij}^*(t)= \P_{c_i^*(t)}^{D_t,*}(Z_{\sigma_{K^*(t)}}^{D_t,*}=c_j^*(t),\; \sigma_{K^*(t)}<\sigma_{g_t(A)})/(1-R_i^*(t)),\ i\neq j.
\end{equation}
By the conformal invariance
of the ABM $Z^\HH$ under the map $g_t$ (see \cite[Theorem 5.3.1]{CF1}),
\begin{equation}\label{e:6.18}
f_j(t,g_t(z))=\P_z^\HH \left( Z_{\sigma_K}^\HH\in C_j,\; \sigma_K<\sigma_A \right)
\end{equation}
and
\begin{equation}\label{e:6.19}
v_t(g_t(z))=\Im g_t(z)-\E_z^\HH \left[\Im g_t(Z_{\sigma_{K\cup A}}^\HH);\; \sigma_{K\cup A}<\infty \right]=:q_t(z).
\end{equation}
Thus by \eqref{e:6.14},
\begin{equation}\label{e:6.21}
\begin{cases}
R_k^*(t)=\int_{\eta_k}\P_z^\HH(Z_{\sigma_K}^\HH\in C_k; \sigma_K<\sigma_A)\nu_k(dz)=:R_k^*,\\
\int_{g_t(\eta_k)}v_t(z)\nu_k^t(dz)
=\int_{\eta_k} q_t(z) \nu_k(dz).
\end{cases}
\end{equation}
Finally we use again the conformal invariance of the BMD $Z^{D,*}$ under the map $g_t$ to get from \eqref{e:6.17}
\begin{equation}\label{e:6.22}
q_{ij}^*(t)= \P_{c_i^*}^{D,*}(Z_{\sigma_{K^*}}^{D,*}=c_j^*,\; \sigma_{K^*}<\sigma_A)/(1-R_i^*)=: q_{ij}^*,\ i\neq j.
\end{equation}
Denote by $M_{ij}$ is the entry of $\sum_{n=0}^\infty (Q^*)^n$ for the matrix $Q^*$ with zero diagonal entry and off-diagonal entry $q_{ij}^*.$
It follows
from \eqref{e:6.15} and \eqref{e:6.18}-\eqref{e:6.22} that
\begin{eqnarray}
 \Im h_t\circ g_t(z)= q_t(z)
+
\sum_{i=1}^N\P_z^\HH\left( Z_{\sigma_K}^\HH\in C_i;\; \sigma_K<\sigma_A\right)\;\sum_{j=1}^N \gamma_{ij}\int_{\eta_j}\; q_t(z)\nu_j(dz), \label{e:6.23}
\end{eqnarray}
for $z\in D\setminus F_t\setminus A,$ where $\gamma_{ij}=\frac{M_{ij}}{1-R_j^*}.$
This together with  \eqref{e:6.2} establishes \eqref{e:6.11}.

(ii)\ As $g_t(z)$ is a solution of the K-L equation \eqref{e:5.25},
$\Im g_t(z)$
  is jointly continuous and satisfies
$0\le \Im g_t(z)\le \Im z.$
The functions $\{\kappa_j (t);
1\leq j \leq N\}$ are continuous due to the continuity of $\s(t).$  Therefore by \eqref{e:6.10},
$q_t(z)$  is jointly continuous.
Since the function $u(z)=\Im z$ is excessive with respect to $Z^\HH$,
 $v_t(z)$ defined by \eqref{e:6.15} is
 non-negative.  Hence $0\le q_t(z)=v_t (g_t (z)) \le \Im g_t(z)\le \Im z$.
It follows from \eqref{e:6.11}
that  $\Im \wt g_t(z)$ is jointly continuous in $(t,z)\in [0,T]\times (\wt D\setminus \wt F_T)$ and
 $q_t(\Phi_A^{-1}z)\le \Im \Phi_A^{-1}z$.  Thus we readily obtain the stated joint continuity with a bound \eqref{e:6.12}.
\qed

\begin{lem}\label{L:6.6}\begin{description}
\item{\rm (i)}
\ For each $T\in (0,\tau_A)$, the functions $\{\wt g_t(z), t\in [0,T]\}$ are extended to be locally equi-continuous and locally uniformly bounded in $z\in (\wt D\cup\partial_p\wt K\cup\partial \HH)\setminus \overline{\wt F_T}\setminus \overline{\wt A}.$

\item{\rm(ii)}\ For $s\ge 0,$
 $\lim_{t\downarrow s} \wt g_{t,s}^{-1}(z)= z$ uniformly on each compact subset of
$\wt D_s\cup \partial_p \wt K(s)\cup (\partial\HH\setminus \{\wt\xi(s)\}\setminus \wt g_s(\wt A)).$

\item{\rm(iii)}
\ For $T\in [0,\tau_A),$
$\wt g_t(\wt z)$ is jointly continuous in $(t,\wt z)\in [0,T]\times [(\wt D\cup \partial_p\wt K\cup\partial\HH)\setminus \overline{\wt F_T}\setminus \wt A].$

\item{\rm(iv)}\ $\wt a_t$ is  right continuous in $t$ and
 $\wt D_t$ is  continuous in $t$.

\item{\rm(v)}
\ $\wt \Psi_t(z, x)$ is jointly continuous in
$(t, z, x) \in \bigcup_{t\in [0,\tau_A)}\{t\}
\times [\wt D_t\cup \partial_p \wt K(t)\cup
(\partial\HH\setminus \{x\})]\times \partial\HH$.

\end{description}
\end{lem}
\pf\ (i) follows from Proposition \ref{P:6.5} (ii) together with Lemma \ref{L:6.4} (iii) exactly in the same way as the proof of \cite[Theorem 7.4]{CFR}.
(ii) follows from (i) and Proposition \ref{P:6.5} (ii) as the proof of \cite[Theorem 8.2]{CFR}.
(iii) can be shown in a quite similar way to (ii).
The right continuity of $\wt a_t$ follows from (ii) as the proof of \cite[Theorem 8.4]{CFR}.  The continuity of $\wt D_t$ is a consequence of (iii).
(v) follows from the continuity of $\wt D_t$ in (iv)
and the Lipschitz continuity of $\wt \Psi$ as the proof of \cite[Theorem 9.8]{CFR}.
\qed

\begin{prop}\label{P:6.2}\ $\wt a_0=0,$ $\wt a_t$ is strictly increasing and
 right continuous.
For each $T\in (0,\tau_A)$ and $z\in \wt D\cup \partial_p\wt K\setminus \wt F_T,$
$\wt g_t(z)$ is right differentiable in $\wt a_t$ and
\begin{equation}\label{e:6.4}
\frac{\partial^+\wt g_t(z)}{\partial\wt a_t}
 = -\pi \wt\Psi_t(\wt g_t(  z),\wt\xi(t)),\quad
\wt g_0(  z)=  z,\quad {\rm for}\quad
t\in [0,T].
\end{equation}
Here the left hand side indicates the right derivative.
\end{prop}
\pf \ This follows from Lemma \ref{L:6.3}, Lemma \ref{L:6.4} and
Lemma \ref{L:6.6}
just as in the proof of Proposition \ref{P:5.11}.
\qed

Note that equation \eqref{e:6.4} does not characterize 
the conformal map $\wt g_t$ since its left hand side involves
only the right derivative. To characterize $\wt g_t$ uniquely,
we need to show that $\wt g_t$ is differentiable in $t$;
see \cite[Remark  2.7]{CFS}. 
The first assertion of the next proposition is crucial not only for this purpose but also in legitimating the stochastic calculus in the next subsection.

The conformal map $h_t(z)$ (resp. $\Phi_A(z)$) from $D_t\setminus g_t(A)$ (resp. $D\setminus A$) onto $\wt D_t$ (resp. $\wt D$) is extended to a conformal map
on
\begin{equation}\label{eq:domain}
\left( D_t\cup \Pi D_t\cup \partial \HH \right) \setminus
\left( g_t(\overline A)\cup \Pi g_t(A) \right)
\quad{\rm(resp.}\left( D\cup \Pi D\cup \partial \HH \right) \setminus
\left( \overline A\cup \Pi A \right){\rm)}
\end{equation}
by the Schwarz reflection.
Note that $h_0(z)=\Phi_A(z)$.
 
\begin{prop}\label{P:6.7}
\ {\rm (i)}\ For any $t\in (0,\tau_A)$ and $z\in D_t\cup\partial\HH\setminus \overline{g_t(A)},$ $h_t(z),\;h_t'(z),\; h_t''(z)$ are jointly continuous
in $(t,z)$.

\noindent
{\rm (ii)}\ Locally uniformly in $z\in (D\cup\partial\HH)\setminus \overline A,$\begin{equation}\label{e:6.23a}
\lim_{t\downarrow 0} h_t(z)=\Phi_A(z),\quad \lim_{t\downarrow 0} h_t'(z)=\Phi_A'(z),\quad\lim_{t\downarrow 0} h_t''(z)=\Phi_A''(z).
\end{equation}
\end{prop}

\pf\ (i)\ It follows from \eqref{e:6.2} that for
$t\in [0,\tau_A),\ 0\le s<t$ and $z\in D_t\setminus g_t(A),$
\begin{equation}\label{e:6.24}
h_t(z)=\wt g_{t,s}^{-1}\circ h_s\circ g_{t,s}(z)
\quad \hbox{ where } \ g_{t,s}=g_s\circ g_t^{-1} \hbox{ and } \  \wt g_{t,s}=\wt g_s \circ \wt g_t^{-1}.
\end{equation}
For $t>0$ and $z\in D_t\cup (\partial\HH\setminus \{\xi(t)\}),$
let $\varphi(u;t,z),\ u\in I_{t,z},$ be the unique solution
of the ODE \eqref{e:5.1} in variable $u$ with initial condition $\varphi (t;t,z)=z$
and with
the maximal time interval $I_{t,z}$ of existence.
If $z\in D_t,$ then $I_{t,z}=[0,t_z)$ by 
Proposition \ref{P:5.4}
and it holds that $\varphi(u;t,z)=g_{t,u}(z)$ for $u\in [0,t].$   But, if $z\in \partial\HH\setminus \{\xi(t)\},$ then
$\varphi(\cdot; t,z)$ is a continuous motion on $\partial\HH$ and
it could be that $I_{t,z}=(\alpha_{t,z},\beta_{t,z})$ with
$0\le \alpha_{t,z}<t<\beta_{t,z}$. 
 Our strategy for the proof of (i) is to use the identity \eqref{e:6.24} for some fixed $s\in (0,t)$
 along with the joint continuity of $\varphi(s;t,z)$
 and that of $\wt g_{t,s}^{-1}(\wt z)=\wt g_t\circ \wt g_s^{-1}(\wt z)$ basically shown in Lemma 6.5.
 
Recall that  $\tau_A =\inf\{u>0:\overline F_u\cap \overline A\neq \emptyset\})$ and
define $\Pi z=\bar z,\; z\in \HH.$. Fix $T\in (0, \tau_A)$.
 Take any smooth Jordan curve $\Gamma\subset \CC$ with $\Pi\Gamma=\Gamma$ such that $\Gamma$ surrounds $\overline F_T,$
the sets $\overline A$ and $K=\cup_{j=1}^N C_j$ are located outside $\Gamma,$ and $\Gamma$ intersects $\partial\HH$ at only two points.  For $t\in (0,T),$ we extend $g_t$ by the Schwarz reflection and let $\Gamma_t=g_t(\Gamma).$   Then $\Gamma_t$ surrounds $\overline{g_t(F_T\setminus F_t)}$ and the sets $\overline{g_t(A)}$ and $K(t)$ are located outside $\Gamma_t.$  In particular, $\xi(t)\notin \Gamma_t$ in view of \eqref{e:5.31}.

From now we fix an arbitrary $t\in (0,T)$ and let $\Gamma_t\cap \partial\HH=\{z_1,z_2\}.$ 
The ODE \eqref{e:5.1} and its solution $\varphi(u;t,z), u\in I_{t,z},$ are extended to $\Pi D_t$ by mirror reflection.
We then choose any $s\in (\alpha_{t,z_1}\vee \alpha_{t,z_2}, t)$ so that $\varphi(s;t,z)$ is well defined for all $z\in \Gamma_t.$
According to a general theorem \cite[Theorem V.2.1]{H} on ODE,
$(t, z)\mapsto \varphi(s;t,z)$ is joint continuous
in the following sense: \ for any
$\eps>0,$ there exists $\delta=\delta(\eps,t,z)>0$ with $s<t-\delta<t+\delta<T$ such that, for any $u>0,\ w\in \CC$ with $|u-t|<\delta,\ |w-z|<\delta,$ we have $\alpha_{u,w}<s$ and
$|\varphi(s;u,w)-\varphi(s;t,z)|<\eps$ for any $z\in \Gamma_t.$
A covering argument then yields the existence of $\delta\in (0,t-s)$ such that
\begin{equation}\label{e:jc1}
\alpha_{u,z}<s\ {\rm and}
\ |\varphi(s;u,z)-\varphi(s;t,z)|<\eps\ \text{\rm for any}\ u\in [t-\delta,t+\delta]\ \text{\rm and for any}\ z\in \Gamma_t.
\end{equation}

Observe that $\varphi(s;t,z)=g_{t,s}(z)$ for $z\in D_t.$  Hence, by the continuity of $\varphi(s;t,z)$ in $z$, we get the identity $\{w=\varphi(s;t,z):z\in \Gamma_t\}=\Gamma_s.$   We can choose $\eps>0$ so that the $\eps-$neighborhood $\Gamma_{s,\eps}$ of $\Gamma_s$ is disjoint from $\overline{g_s(F_T\setminus F_s)}\cup \overline{g_s(A)}.$   On account of the relation $h_s(w)=\wt g_s\circ \Phi_A\circ g_s^{-1}(w),\ w\in D_s,$
we have
\begin{equation}\label{e:jc2}
h_s(\Gamma_{s,\eps})\subset \wt D_s\cup\Pi\wt D_s\cup \partial\HH\setminus
\left(\overline{\wt g_s(\wt F_T\setminus \wt F_s)}\cup \wt \Pi \wt g_s(\wt F_T\setminus \wt F_s))\cup \wt g_s(\wt A)\cup \Pi\wt g_s(\wt A) \right).
\end{equation}

On the other hand, we have the
following
variant of Lemma \ref{L:6.6} (iii):\
\begin{eqnarray}\label{e:jc3}
&&{\rm For}\ T\in (0,\tau_A)\ {\rm and}\ s\in [0, T),
\ \wt g_{u,s}^{-1}(\wt z)=\wt g_u(\wt g_s^{-1}(\wt z)) \ \text{\rm is jointly continuous in}\nonumber \\
&& (u,\wt z)\in [s,T]\times [(\wt D_s\cup \partial\HH)\setminus \overline{\wt g_s(\wt F_T\setminus \wt F_s)}\setminus \wt g_s(\wt A)].
\end{eqnarray}
This can be proved as follows.
 By using the relation \eqref{e:6.24}, we first express $\Im \wt g_{u,s}^{-1}$, $u\ge s$, 
   in terms of the BMD on $\wt D_s$ and the ABM on $\HH$
 in analogy to \eqref{e:6.11}, which yields the joint continuity of 
$\Im \wt g_{u,s}^{-1}(\wt z)$ in
$(u,\wt z)\in [s,T]\times [\overline \HH\setminus \overline{\wt g_s(\wt F_T\setminus \wt F_s)}\setminus \wt g_s(\wt A)].$
This combined with Lemma \ref{L:6.4} (iii) (replacing $(s,t)$ by $(u,T)$) implies, in the same way as the proof of \cite[Theorem 7.4]{CFR}, the local uniform boundedness of the family $\{\wt g_{u,s}^{-1}(\wt z); u\in [s,T]\}$ in
$\wt z\in \wt D_s\cup \partial \HH\setminus \overline{\wt g_s(\wt F_T\setminus \wt F_s)}\setminus \wt g_s(\wt A)].$  Thus we can get \eqref{e:jc3} as the proof of \cite[Theorem 8.2]{CFR}.

By \cite[Theorem V.2.1]{H} again,
$\varphi(s;u,z)$ is jointly continuous in
$(u,z)\in [t-\delta,t+\delta]\times \Gamma_t.$
Since $h_s(\varphi(s;u,z))\in h_s(\Gamma_{s,\eps})$ by \eqref{e:jc1}, we conclude from \eqref{e:jc2} and \eqref{e:jc3} that the relation \eqref{e:6.24} extends to $h_u(z)=\wt g_{u,s}^{-1}(h_s(\varphi(s;u,z))$ to be jointly continuous at each
$(u,z)\in [t-\delta,t+\delta]\times \Gamma_t.$
In particular $\sup_{u\in [t-\delta,t+\delta],\; z\in \Gamma_t} |h_u(z)|$
is finite. Moreover, by the joint continuity of the solution of \eqref{e:5.1},
we may assume that $\Gamma_t \subset \bigcap_{u\in [t-\delta,t+\delta]}(D_u \setminus \overline{g_u(A)})$. 

As $h_u$ is analytic, the Cauchy integral formula yields that $h_u(z),\;h_u'(z),\;h_u''(z)$ are jointly continuous in 
$(u,z)\in [t-\delta,t+\delta] \times U(t)$, where $U(t)$ is an open set  enclosed by $\Gamma_t.$

\smallskip
\noindent
(ii)\ We continue to work with the function $\varphi(u; t,z),\ u\in I_{t,z},$ as  above
and claim the following: for any $\eps>0,$ there exists $\delta>0$ such that, for any $t\in (0,\delta]$ and any $z\in \partial\HH\setminus [\xi(0)-\eps, \xi(0)+\eps],$ $I_{t,z}=[0,\beta)$ for some $\beta>t.$

To see this, we fix $\eps_1\in (0,\eps)$ and take $t_0>0$ with $\{\xi(u):u\in [0,t_0]\}\subset (\xi(0)-\eps_1,\xi(0)+\eps_1).$   Since the solution $\varphi(t,u,\xi(0)\pm \eps_1)$ of \eqref{e:5.1} with $\varphi(u,u,\xi(0)\pm \eps_1)=\xi(0)\pm \eps_1$ is jointly continuous in $(t,u),$ there is $\delta\in (0,t_0]$ with
$
\xi(0)-\eps<\inf_{0\le u\le t\le \delta}\varphi(t,u,\xi(0)-\eps_1),
\ \sup_{0\le u\le t\le \delta}\varphi(t,u,\xi(0)+\eps_1)<\xi(0)+\eps.
$
Take any $t\in (0,\delta]$ and any $z\in \partial\HH$ with $z>\xi(0)+\eps.$  Suppose
$I_{t,z}=(\alpha,\beta)$ for some $\alpha\in (0,t).$
As $\liminf_{u\downarrow \alpha} |\varphi(u,t,z)-\xi(u)|=0,$ we find $u_1\in (\alpha,t)$ with $\varphi(u_1,t,z)=\xi(0)+\eps_1,$ arriving at a contradiction
$z=\varphi(t,u_1, \xi(0)+\eps_1)<
 \xi(0)+\eps$.
 Hence $I_{t,z}=[0,\beta).$  The same is true for $z\in \partial\HH$ with $z<\xi(0)-\eps.$

Observe that $g_t^{-1}(z)=\varphi(0,t,z)$ for
$(t,z)\in [0,\delta]\times (D \cup [(\partial\HH\setminus[\xi(0)-\eps,\xi(0)+\eps])],$
and it is jointly continuous in $(t,z)$ there by the theorem cited above.
 Let $V$ be a compact subset of $D\cup (\partial\HH\setminus [\xi(0)-\eps,\xi(0)+\eps])\setminus \overline A.$  We may assume that $\delta<\tau_A$ and $V$ is disjoint from $\bigcup_{s\in [0,\delta]}\overline{g_s(A)}.$

Combining this
with the identity $h_t=\wt g_t\circ\Phi_A\circ g_t^{-1}$ from \eqref{e:6.2} and  with Lemma \ref{L:6.6} (ii), (iii), we see that $h_t(z)$ is jointly continuous at each $(t,z)\in [0,\delta]\times V$ and consequently $\sup_{t\in [0,\delta]\;z\in V}|h_t(z)|$ is finite and $\lim_{t\downarrow 0}h_t(z)=\Phi_A(z)$ for each $z\in V.$  By taking appropriate circles as $V$,
we get the local uniform convergence \eqref{e:6.23a} in a similar
way as in the proof of (i).
\qed

In the remaining part of this paper, the derivative of a function $f$ in the time parameter will be designated by $\dot f $.

\begin{thm}\label{T:6.9}
\ For $s\in (0,\tau_A)$ and $z\in \wt D\cup \partial_p\wt K\cup\wt F_t$,  $\wt g_s(z)$ is continuously
differentiable in $s\in [0.t]$ and
\begin{equation}\label{e:6.25}
\frac{d\wt g_s(z)}{ds}
= -2\pi |h_s'(\xi(s))|^2\;\wt\Psi_s(\wt g_s(  z),\wt\xi(s)),
\quad g_0( z)= z.
\end{equation}
\end{thm}

\pf\ It suffices to prove
\begin{equation}\label{e:6.26}
\dot{\wt a_s}\;=\; 2|h_s'(\xi(s))|^2.
\end{equation}
This is because \eqref{e:6.26} together with \eqref{e:6.4} implies that
\eqref{e:6.25} holds with the right derivative $\frac{\partial^+\wt g_s(z)}{ds}$
in place of $\frac{\partial\wt g_s(z)}{ds}$. But since the right hand side of \eqref{e:6.25} is continuous in $s$
in view of Lemma \ref{L:6.6} and
Proposition \ref{P:6.7},
$\wt g_s(\wt z)$ is actually continuously differentiable in $s$.

For $D\in \dd$ and an $\HH$-hull $K\subset D,$ we denote by ${\rm Cap}^D(K)$ (resp. ${\rm Cap}^\HH(K)$) the half-plane capacity of $K$ for the canonical map from $D\setminus K$ (resp. the Riemann map from $\HH \setminus K$ onto $\HH$).  For a set $A\subset \HH$, we put ${\rm rad}(A)=\sup_{z\in A}|z|.$
Fix $s\ge 0$ and let
$K_\epsilon=g_s(F_{s+\epsilon}\setminus F_s)$ and $\wt K_\epsilon=\wt g_s(\wt F_{s+\epsilon}\setminus \wt F_s)$.  Then
${\rm rad}(K_\epsilon \setminus \{\xi(s)\})=o(\epsilon)$ and ${\rm rad}(\wt K_\epsilon \setminus \{\wt\xi(s)\})=o(\epsilon)$ by \eqref{e:5.31} and \eqref{e:6.5}, respectively.  Consequently we have by Theorem \ref{T:7.1} of Appendix
\begin{equation}\label{e:cap}
{\rm Cap}^{D_s}(K_\epsilon)-{\rm Cap}^\HH(K_\epsilon)= o(\epsilon),\quad
{\rm Cap}^{\wt D_s}(\wt K_\epsilon)-{\rm Cap}^\HH(\wt K_\epsilon)= o(\epsilon).
\end{equation}

Since $\wt K_\epsilon=h_s(K_\epsilon),$ we get from \eqref{e:cap} and \cite[(4.15)]{L1} that
\begin{eqnarray*}
&& \wt a_{s+\epsilon}-\wt a_s
= {\rm Cap}^{\wt D_s}(\wt K_\epsilon)
= {\rm Cap}^\HH(h_s(K_\epsilon)) + o(\epsilon)\\
&&= \Phi'_s(\xi(s))^2 {\rm Cap}^\HH(K_\epsilon) + o(\epsilon)
= \Phi'_s(\xi(s))^2 {\rm Cap}^{D_s}(K_\epsilon) + o(\epsilon)
= \Phi'_s(\xi(s))^2 (a_{s+\epsilon} - a_s) + o(\epsilon),
\end{eqnarray*}
which yields \eqref{e:6.26} as $a_{s+\epsilon}-a_s=2\epsilon$ by Theorem \ref{T:5.10}. \qed

\subsection{Characterization of locality of ${\rm SKLE}_{\alpha, -b_{\rm BMD}}$}\label{S:6.3}

We continue to operate under the setting in the preceding subsection.
To investigate the locality, we need to compute the driving processes
for $\{\wt F_t; t < \tau_A\}$.
It follows from \eqref{e:5.25} that the inverse map $g_t^{-1}$ of $g_t$ satisfies
 \begin{equation}\label{e:6.28}
\dot g_t^{-1}(z)=2\pi (g_t^{-1})'(z)\Psi_{\s (t)}(z, \xi(t)),\quad g_0^{-1}(z)=z .
\end{equation}
From \eqref{e:6.2}, we have
$$
\dot h_t(z)=\dot{\wt g}_t(\Phi_A\circ g_t^{-1}(z))+
(\wt g_t\circ\Phi_A)'(g_t^{-1}(z))
\dot g_t^{-1}(z),   \quad z\in D_t\setminus g_t(A).
$$
This together with \eqref{e:6.28}, Theorem \ref{T:6.9}, and then by  \eqref{e:6.2} again yields
that for $ z\in D_t\setminus g_t(A)$,
\begin{eqnarray}
\dot h_t(z) &=& -2\pi |h_t'(\xi(t))|^2\wt \Psi_t(\wt g_t\circ\Phi_A\circ g_t^{-1}(z),\wt \xi(t))
+ (\wt g_t\circ\Phi_A)'(g_t^{-1}(z)) 2\pi (g_t^{-1})'(z) \Psi_{\s(t)}(z,\xi(t))
\nonumber  \\ \label{e:6.29}
&= & -2\pi |h_t'(\xi(t))|^2\wt \Psi_t(h_t(z),h_t(\xi(t)))
 + 2\pi h_t'(z) \Psi_{\s (t)} (z,\xi(t)).
\end{eqnarray}

Functions $h_t(z)$ and $h_t'(z)$ are
extended  to the region \eqref{eq:domain}, call it $G_t$,
 by the Schwarz reflection.   
 Fix $t_0>0$ and take a disk $B$ centered at $\xi(t_0)$ with $\overline B\subset \cap_{|t-t_0|\le\delta}G_t$ and $\{\xi(t): |t-t_0|\le \delta\} \subset B$ for some $\delta>0.$
 Denote the right hand side of \eqref{e:6.29} by $f(t,z).$
By virtue of Proposition \ref{P:6.7}(i), Lemma \ref{L:6.6}(v) and Proposition \ref{P:5.1}(i), $f(t,z)$ is jointly continuous 
and hence uniformly bounded in $(t,z)\in [t_0-\delta, t_0+\delta]
\times (\partial B \cap \overline \HH)$.
By taking Schwarz reflections of $\wt\Psi_t(z,h_t(\xi(t)))$ and $\Psi_{\s(t)}(z,\xi(t))$ in $z$, $f(t,z)$ admits an extension to $[t_0-\delta,t_0+\delta]\times \partial B$ to be jointly continuous and uniformly bounded there, and the identity $\dot h_t(z)=f(t,z)$ extends to $(t,z)\in (t_0-\delta, t_0+\delta)\times (\partial B \setminus \partial \HH).$

Expressing 
$(h_u(z)-h_t(z))/(u-t),\ z\in B,\; t\in (t_0-\delta,t_0+\delta),$ 
by the Cauchy integral formula and letting $u\to t$, we see that $h_t(z)$ is differentiable in $t$ for any $z\in B$
with $\dot h_t(z)$ being analytic in $z\in B$ and jointly continuous in $(t,z)\in (t_0-\delta,t_0+\delta)\times B.$
In particular, 
$\dot h_t(\xi(t))$ can be computed by $\lim_{z\to\xi(t),\;z\in \HH}\dot h_t(z)$ explicitly.
Indeed, by the definition \eqref{e:6.1} of $b_{\rm BMD}(\s,\xi)$, we get from \eqref{e:6.29}
\begin{eqnarray}
\dot h_t(\xi(t)) &=&
h_t'(\xi(t)) \, b_{\rm BMD}(\xi(t),\s(t))
 - |h_t'(\xi(t))|^2 \, b_{\rm BMD}(h_t(\xi(t)), h_t(\s(t)))
 \nonumber \\
&& + \lim_{z\to \xi(t)} \left( \frac{2|h_t'(\xi(t))|^2}{h_t(z)-h_t(\xi(t))}
-\frac{2h_t'(z)}{z-\xi(t)} \right) \nonumber \\
&=& h_t'(\xi(t)) \, b_{\rm BMD}(\xi(t),\s(t))
 - |h_t'(\xi(t))|^2 \, b_{\rm BMD}(h_t(\xi(t)), h_t(\s(t)))
 -3h_t''(\xi(t)) .   \label{e:6.31}
\end{eqnarray}
Thus $h_t(z)$ is differentiable in $t$ for
 each $z\in \partial\HH\cap G_t$ with $\dot h_t(z)$ being jointly continuous in $t>0,\ z\in \partial\HH\cap G_t.$ 
Moreover,  $h_t'(z)$ and $h_t''(z)$ are jointly continuous by Proposition \ref{P:6.7}.
 Since $\wt \xi(t)=h_t(\xi(t))$ and $\xi(t)$ is the solution of the SDE \eqref{e:3.39}, 
 we can readily 
  apply a generalized It\^o formula 
 (see \cite[Remark 2.9]{CFS} and \cite[(IV.3.12)]{RY})  to  get
$$
d\wt \xi(t)=\left( \dot h_t(\xi(t)) +h_t'(\xi(t))b(\s(t)-\wh \xi (t))
 + \frac12 h_t''(\xi(t))\alpha(\s(t)-\wh \xi (t))^2\right) dt
+h_t'(\xi(t))\alpha(\s(t)-\wh \xi (t))dB_t.
$$
This combined with \eqref{e:6.31} gives the following.

\begin{thm}\label{T:6.10}\ It holds that
\begin{eqnarray}
d\wt \xi(t) &=& h_t'(\xi(t))\left( b(\s(t)-\wh \xi (t))+b_{\rm BMD}(\xi(t),\s(t))\right) dt
 +\frac12 h_t''(\xi(t))\left( \alpha(\s(t)-\wh \xi (t))^2 -6\right) dt\nonumber\\
&&-|h_t'(\xi(t))|^2 b_{\rm BMD}(\wt{\xi}(t), h_t(\s(t)))dt
+h_t'(\xi(t))\alpha(\s(t)-\wh \xi (t))dB_t.\label{e:6.33}
\end{eqnarray}
\end{thm}

Let $\{\c F_t\}_{t<\c\tau_A}$ be the half-plane capacity reparametrization of the image hulls $\{\wt F_t\}_{t<\tau_A}$,
namely,
 \begin{equation}\label{e:time}
\c F_t=\wt F_{\wt a^{-1}(2t)},\qquad
\c \tau_A=  \wt a(\tau_A) / 2.
\end{equation}
where $\wt a(t)$ is the half-plane capacity  of $\wt F_t$ and $\wt a^{-1}$ is its inverse function.
Accordingly, the processes $\wt \xi(t)=h_t(\xi(t)) =\wt g_t\circ \Phi_A(\xi) $ and
$\wt \s_j(t)=h_t(\s_j(t)) =\wt g_t\circ \Phi_A(\s_j)$ are time-changed into
\begin{equation}\label{e:6.35}
\c\xi(t )=\wt \xi(\wt a^{-1}(2t)) \quad{\rm and}\quad
\c\s_j( t)=\wt \s_j(\wt a^{-1}(2t)),\quad 1\le j\le 3N,\quad t<\c\tau_A,
\end{equation}
respectively.

Set
$\c g_{t}=\wt g_{\wt a^{-1}(2t)}$ and $\c \Psi_{t}=\wt \Psi_{\wt a^{-1}(2t)}$.
It follows from \eqref{e:6.25}, \eqref{e:6.26}, Lemma \ref{L:6.6} (v) and
Proposition \ref{P:6.7}(i)
that, for $T\in (0,\c\tau_A),$ $\c g_t(z)$ is continuously differentiable in $t\in [0,T]$ and
\begin{equation}\label{e:6.38}
\frac{d\c g_t(z)}{dt} = -2\pi \c\Psi_t(\c g_t(z),\c\xi(t)),\quad
\c g_0(z)= z\in \wt D\cup \partial_p\wt K\setminus \c F_t.
\end{equation}

\begin{lem}\label{L:6.11}\ It holds under $\P_{(\xi,\s)}$ that
\begin{equation}\label{e:6.36}
\c \s_j(t)=\Phi_A(\s_j)+\int_0^t \c b_j(\c\xi(s),\c\s(s))ds,
\quad t\in [0,\c \tau_A),\quad 1\le j\le 3N,
\end{equation}
where $\c b_j(\w)=\c b_j(\xi,\s)$ is defined by {\rm\eqref{e:3.38}} with $\Psi_s$ being replaced by $\c \Psi_s.$
\end{lem}

\pf\ We can get \eqref{e:6.36} from the K-L equation \eqref{e:6.38} exactly in the same way as the proof of Theorem \ref{T:2.3}, if Lemma \ref{L:2.1} for $\c g_t,\; \c\tau_A,\; \c\Psi_t$ in place of $g_t,\; t_\gamma,\; \Psi_t$ is once established.  Let us call Lemma 2.1' such a counterpart of Lemma \ref{L:2.1}.

The first and second assertions of Lemma 2.1' follow from Lemma \ref{L:6.6} (v),
Proposition \ref{P:6.7}(i),
\eqref{e:6.26} and \eqref{e:6.38} as in the proof of those of Lemma \ref{L:2.1}.   The third assertion of Lemma 2.1' can be obtained by proving an analogue to \eqref{e:2.5} using a similar method to the proof of \eqref{e:6.4} combined with \eqref{e:6.26}.  The rest of assertions of Lemma 2.1' can be proved quite similarly.
\qed

Let $M_t=\int_0^t h_s'(\xi(s))dB_s$. Clearly by \eqref{e:6.26},
$\< M\>_t=\int_0^t h_s'(\xi(s))^2ds= \wt a(t)/2$.
Hence
$\check B_t:=M_{\wt a^{-1}(2t)}$ is a Brownian motion.
 The formula \eqref{e:6.33} can be rewritten as
\begin{eqnarray}\label{e:6.36n}
 \check \xi (t)&=&
 \Phi_A(\xi(0)) + \int_0^t\c h_s'(\st\xi(s))^{-1}\left( b(\st\s(s)-\wh{\st\xi}(s))+b_{\rm BMD}(\st\xi(s),\st\s(s))\right) ds\nonumber\\
&&+\frac12 \int_0^t \c h''_s(\st\xi(s))\cdot \c h_s'(\st\xi(s))^{-2}  \left( \alpha(\st\s(s)-\wh{\st\xi}(s))^2-6 \right) ds\nonumber\\
&& -\int_0^t  b_{\rm BMD} (\check \xi (s), \check \s (s)) ds + \int_0^t \alpha(\st \s(s)-\wh{\st\xi}(s)) d\check B_s,\quad
t\in [0,\c\tau_A),
\end{eqnarray}
where
$\c h_s'(z):=h_{\wt a^{-1}(2s)}'(z)$, $\c h_s''(z):=h_{\wt a^{-1}(2s)}''(z)$,
$\st\xi(t ):= \xi(\wt a^{-1}(2t))$ and $\st\s_j( t)= \s_j(\wt a^{-1}(2t))$
for  $1\le j\le 3N$.
Note that since $h_t(z)$ is univalent in $z$ on the region \eqref{eq:domain},  $h_t'(z)$ never vanishes there.

Let $\{F_t\}$ be a
${\rm SKLE}_{\alpha,b}$.
Since $\{F_t\}$ depends also on the initial value $(\xi,\s)$ for  SDE \eqref{e:3.39}-\eqref{e:3.40},
 we shall write ${\rm SKLE}_{\alpha,b}$ occasionally
    as ${\rm SKLE}_{\xi,\s,\alpha,b}$ for emphasis on its dependence on the initial
    position $(\xi, \s)$.
   Recall that, for an $\HH$-hull $A\subset D(\s),$
$\tau_A=\inf\{t>0: \overline F_t\cap \overline A\neq \emptyset\}.$
Let $\{\c F_t\}_{\{t<\c\tau_A\}}$ be the half-plane capacity reparametrization of the image hulls $\{\wt F_t=\Phi_A(F_t)\}_{\{t<\tau_A\}}$
specified by \eqref{e:time}.

${\rm SKLE}_{\alpha,b}$ is said to
have  the {\it locality} property
if, for the ${\rm SKLE}_{\xi,\s,\alpha,b}$ $\{F_t\}$ with an arbitrarily fixed $(\xi,\s)\in \RR\times \cS$ and for any $\HH$-hull $A\subset D(\s),$
$\{\c F_t,\; t<\c\tau_A \}$ has the same distribution as ${\rm SKLE}_{\Phi_A(\xi),\Phi_A(\s),\alpha, b}$ restricted to $\{t< \tau_{\Phi_A (A)}\}$.
Here ${\rm SKLE}_{\alpha,b}$ and  ${\rm SKLE}_{\Phi_A(\xi),\Phi_A(\s),\alpha, b}$  can live on two different probability spaces.

\begin{thm}\label{T:6.12}\
 ${\rm SKLE}_{\alpha, -b_{\rm BMD}}$
for a constant $\alpha > 0$
enjoys the locality if and only if $\alpha=\sqrt{6}.$

\end{thm}

\pf \ {\bf ``If" part.}
Assume that
$\alpha = \sqrt{6}$ and
$b(\xi,\s):=b(\s -\wh  \xi)=-b_{\rm BMD}(\xi,\s).$
Then \eqref{e:6.36n} is reduced to
\begin{equation}\label{e:6.42}
d\c \xi(t)=-b_{\rm BMD}(\c \xi(t), \c\s(t))dt
+ \sqrt{6}\;d\c B_t.
\end{equation}
Thus
$\{\c F_t\}$ is an increasing sequence of $\HH$-hulls associated with the unique solution $\c g_t$ of the Komatu-Loewner equation \eqref{e:6.38},
driven by $(\c \xi(t), \c \s(t))$, which is
the unique solution of \eqref{e:6.42} and \eqref{e:6.36}.  Therefore
$\{\c F_t\}_{\{t<\c\tau_A\}}$ is
${\rm SKLE}_{\Phi_A(\xi),\Phi_A(\s),
\sqrt{6}, -b_{\rm BMD}}$ restricted to
$\{t<\tau_{\Phi_A(A)}\} $,
 yielding the `if' part of the theorem.

 \medskip
\noindent
{\bf ``Only if" part.}
Assume that
$\dis 
\alpha\ \text{\rm is a positive constant} 
\ \hbox{and}\ b(\xi,\s)=-b_{\rm BMD}(\xi,\s).$
Then \eqref{e:6.36n} is reduced to
\begin{equation}\label{e:6.39}
\c\xi(t)=\Phi_A(\xi)+
\frac{\alpha^2-6}{2} \int_0^t \c h''_s(\st\xi(s))\cdot \c h_s'(\st\xi(s))^{-2} ds
 -\int_0^t  b_{\rm BMD} (\check \xi (s), \check \s (s)) ds + \alpha \check B_t,\quad t>0,
\end{equation}
Let $\{F_t\}$ be a ${\rm SKLE}_{\xi, \s,\alpha, -b_{\rm BMD}}$,  $A\subset D(\s)$
an $\HH$-hull, and $\{\c F_t\}$
be defined by \eqref{e:time}.
The equations \eqref{e:6.39} for $\check \xi$ and \eqref{e:6.36} for $\check\s$ describe the evolution of $\{\c F_t\}$ through \eqref{e:6.38}.

Assume now the locality of  ${\rm SKLE}_{\alpha, -b_{\rm BMD}}.$  Then
$\{(\c \xi (t)), \c \s (t); t\in (0, \c \tau_A )\}$
has the same distribution as
the solution $\{(\bar \xi (t), \bar \s (t));
t\in [0,   \bar \tau_{\Phi_A(A)} )\}$
of the equation
\begin{equation}\label{e:6.38n}
  \bar \xi(t)= \Phi_A(\xi)
  - \int_0^t b_{\rm BMD} (\bar\s(s)-\wh{\bar \xi} (s)) ds
  + \alpha  \bar B_t
\end{equation}
for some Brownian motion $\bar B_t$ coupled with the equation \eqref{e:6.36}  with $(\bar\xi(t),\bar\s(t))$ in place of $(\c\xi(t),\c s(t))$.

On the other hand, if we let
$$
\eta(t):= \frac{\alpha^2-6}{2} \int_0^t \c h''_s(\st\xi(s))\cdot \c h_s'(\st\xi(s))^{-2} ds,
$$
then we see from \eqref{e:6.39} that $\c\xi(t)$ is, under the Girsanov transform generated by the local martingale $-\alpha^{-1} \eta(t)\; d\c B_t$, locally equivalent in law to $\bar \xi(t).$  It follows that $\eta(t)=0,\ t<\c \tau_A,$ almost surely, and accordingly
\begin{equation}\label{e:6.40}
(\alpha^2-6) \int_0^{\wt a^{-1}(2t)\wedge\tau_A} h_s''(\xi(s))ds=0,\quad t>0.
\end{equation}

Dividing \eqref{e:6.40} by $\wt a^{-1}(2t)$ and then letting $t\downarrow 0$, we get $\dis (\alpha^2-6)\Phi_A''(\xi)=0$ for every
$ \xi\in \partial\HH\setminus \overline A$ by virtue of Proposition \ref{P:6.7}.  If $\alpha^2 \neq 6$, then $\Phi_A''(\xi) = 0$ for every $\xi\in \partial \HH\setminus A,$
This would imply that $\Phi_A$ is an identity map,
which is impossible unless $A=\emptyset$.
\qed

\begin{remark}\rm\  {\bf (An effect of the second order BMD domain constant)}\\
Along with the BMD domain constant $b_{\rm BMD}$ introduced in Section
\ref{S:6.1}, we 
define for $D\in \dd$
\begin{equation}
c_{\rm BMD}(\xi;D)=2\pi \lim_{z\to\xi}\left(\Psi_D'(z,\xi)-\frac{1}{\pi} \frac{1}{(z-\xi)^2}\right),\quad \xi\in \partial \HH,
\end{equation}  
which is a well defined real number by Lemma 5.6 (ii).  
 We also denote it by 
$c_{\rm BMD}(\xi;\s)$ for $\s=\s(D).$  We set $c_{\rm BMD}(\s)=c_{\rm BMD}(0,\s),\ \s\in {\cal S},$ and call it the {\it second order BMD domain constant}.  On account of (3.31), we then have $c_{\rm BMD}(\xi;\s)=c_{\rm BMD}(\s-\wh \xi)
$ for $s\in {\cal S}$ and $\xi\in \RR$.  

For a constant $\alpha\in (0,2),$ ${\rm SKLE}_{\alpha,b}$ is generated by a simple curve just as ${\rm SLE}_{\alpha^2}$ (\cite{CFS}). 
 As is well known, ${\rm SLE}_{8/3}$ enjoys the so called {\it restriction property} 
that was established by showing that
  $h'_t(\xi_t)^{5/8}$ is a local martingale in \cite{LSW2}. 
 Here $h$ and $\xi$ were defined for the SLE in exactly the same manner as above
for the SKLE.  But we can hardly expect a straightforward generalization of this martingale property to ${\rm SKLE}_{\sqrt{8/3}, -b_{\rm BMD}}$ due to the effect of the second order BMD domain constant $c_{\rm BMD}$ as will be seen below.  See \cite[\S 6]{CFS} for some related literatures.

 It follows from the identity \eqref{e:6.29} that
\begin{eqnarray*}
\dot h_t'(z)= -2\pi|h_t'(\xi(t))|^2 \wt \Psi_t'(h_t(z), h_t(\xi(t)))h_t'(z)
+ 2\pi h_t'(z)\Psi_{\s(t)}'(z,\xi(t))+2\pi h_t''(z)\Psi_{\s(t)}(z, \xi(t)).
\end{eqnarray*}
We then have analogously to \eqref{e:6.31}
\begin{eqnarray*}
\dot h_t'(\xi(t)) &=&
-|h_t'(\xi(t))|^2 c_{\rm BMD}(h_t(\xi(t)),h_t(\s(t)))h_t'(\xi(t))\\
&&+h_t'(\xi(t)) c_{\rm BMD}(\xi(t),\s(t))+h_t''(\xi(t))b_{\rm BMD}(\xi(t),\s(t))+\lim_{z\to\xi(t)} II(z,t),
\end{eqnarray*}
where
\[
II(z,t)= -2|h_t'(\xi(t))|^2\frac{1}{(h_t(z)-h_t(\xi(t))^2}h_t(z)\\
+ h_t'(z)\frac{2}{(z-\xi(t))^2}-h_t''(z)\frac{2}{z-\xi(t)}.
\]
It holds as in \cite[\S 5]{LSW2} that
$\dis \lim_{z\to\xi(t)}II(z,t)=\frac{h_t''(\xi(t))^2}{2h_t'(\xi(t))}-\frac43 h_t'''(\xi(t)).$

Consider the process $\eta(t)=h_t'(\xi(t))^\delta$ for $\delta>0.$  Using a generalized It\^o formula, 
we have 
\begin{eqnarray*}
\frac{1}{\delta}\;\frac{d\eta(t)}{\eta(t)}
&=& -|h_t'(\xi(t))|^2c_{\rm BMD}(h_t(\xi(t)),h_t(\s(t))dt+c_{\rm BMD}(\xi(t),\s(t))dt\\
&&+\frac{h''_t(\xi(t))}{h_t'(\xi(t))}\left\{ b(\xi(t),\s(t))+b_{\rm BMD}(\xi(t), \s(t))\right\}dt\\
&&+\frac12\{(\delta-1)\alpha(\xi(t),\s(t))^2+1\}\frac{h_t''(\xi(t))^2}{h_t'(\xi(t))^2}dt\\
&&+\left( \frac12\alpha(\xi(t),\s(t))^2-\frac43\right) \frac{h'''_t(\xi(t))}{h_t'(\xi(t))}dt
+\frac{h_t''(\xi(t))}{h_t'(\xi(t))}\alpha(\xi(t),\s(t))dB_t.
\end{eqnarray*}
When
$\dis \alpha=\sqrt{8/3}$, $b=-b_{\rm BMD}$ and $\delta=5/8$, 
we get from the above identity
\begin{eqnarray}
\frac{d\eta(t)}{\eta(t)}&=&\sqrt{\frac83}\frac58 \frac{h_t''(\xi(t))}{h_t'(\xi(t))}dB_t
 +\frac58\left( c_{\rm BMD}(\xi(t),\s(t))-|h_t'(\xi(t))|^2c_{\rm BMD}(h_t(\xi(t)),h_t(\s(t)))\right) dt. \nonumber \\
\end{eqnarray}
The drift term of the right hand side does not vanish  
unless either $h_t$ is the identity map or $c_{\rm BMD}$ is vanishing. 
\end{remark}

\begin{remark} \rm{\bf  ERBM and BMD}
\quad As is explained in Introduction, the derivation of the Komatu-Loewner equation and its fundamental properties in \cite{CFR} is partly based on the probabilistic considerations in terms of the Brownian motion with darning (BMD).  We had constructed and characterized the darning of a general symmetric Markov process (\cite[\S 7.7]{CF1}) when we encountered an article of G. Lawler \cite{L2} where the Komatu-Loewner equation on a standard slit domain
previously obtained analytically
by Bauer-Friedrich \cite[Theorem 3.1]{BF3} was investigated in terms of the excursion reflected Brownian motion (ERBM).  We were strongly motivated by these papers.  In the present paper, the BMD is also used crucially in Section \ref{S:6.2} to derive  
the generalized Komatu-Loewner equation \eqref{e:6.25} 
for the image hulls and in Appendix (\S 7) to extend 
Drenning's result \cite{D} on the comparison of half-plane capacities.

\cite[\S 6]{CF2} gives a detailed proof of the identification of ERBM with BMD (especially in the doubly connected case).  Some comprehensive account on BMD and BMD-harmonic functions can be found in \cite{C, CFR, FK}.

\end{remark}

\section{Appendix: Comparison of half-plane capacities}\label{S:7}

We fix a standard slit domain $D=\HH\setminus K,\ K=\cup_{j=1}^N C_j.$
For $r>0$, define $B_r=\{z\in \CC: |z|<r\}$.  Let $T>0$.
 We consider an increasing  family $\{F_t; t\in (0,T]\}$ of $\HH$-hulls such that
there is an increasing sequence of positive numbers  $r_t$ so that
\begin{equation}\label{e:7.1}
\lim_{t\to 0} r_t =0 \quad \hbox{and} \quad F_t \subset B_{r_t} \hbox{ for }
t\in (0, T].
\end{equation}
 Let $a_t$ be the half-plane capacity of the hull $F_t$.
Let $g_t^0$ be the unique Riemann map from $\HH\setminus F_t$ onto $\HH$ satisfying the hydrodynamic normalization $g_t^0(z)=z+\frac{a_t^0}{z}+o\left( {1}/{|z|}\right)$
at infinity. Clearly,  $a_t^0=\lim_{z\to\infty}z(g_t^0(z)-z)$.

\begin{thm}\label{T:7.1}
$\lim_{t\downarrow 0} {a_t}/{t}$
exists if and only if
$\lim_{t\downarrow 0} {a_t^0}/{t}$ exists.
If both limits exist,
they have the same value.
\end{thm}

\bigskip
When $\{F_t\}$ are Jordan subarcs, such a statement of comparison has appeared in   S. Drenning \cite[Lemma 6.24]{D}.
Its proof uses a probabilistic expression of $a_t$ in terms of the excursion reflected Brownian motion (ERBM) for $D$.  A key step of its proof is \cite[Proposition 4.5]{D}, where an estimate of the ERBM-Poisson kernel under a small perturbation of the standard slit domain $D$ is obtained
using an expression of the ERBM-Poisson kernel
that involves the boundary Poisson kernel, excursion measures and an induced finite Markov chain among the holes.
But BMD counterpart of \cite[Proposition 4.5]{D}
to be
formulated in Proposition \ref{P:7.2} below admits a more straightforward proof due to a simpler expression of the BMD-Poisson kernel in \cite{CFR}.

 Denote by $D^*=D\cup \{c_1^*,\cdots, c_N^*\}$  the space obtained from $\HH$ by rendering each hole $C_i$ into a single point $c_i^*$.
Fix $\eps_0>0$ with $B_{\eps_0}\cap \HH\subset D.$  For any $\eps\in (0,\eps_0),$ we consider perturbed domains
\[D_\eps=D\setminus \overline B_\eps,\qquad D_\eps^*=D^*\setminus \overline
B_\eps=D_\eps \cup \{c_1^*,\cdots, c_N^*\}.
\]
Let $K^*_{D}(z,\zeta),\; z\in D^*,\; \zeta\in \partial\HH,$
(resp. $K^*_{D_\eps}(z,\zeta),\; z\in D^*_\eps,\;
\zeta\in \partial(\HH\setminus B_\eps),$) be the Poisson kernel of BMD on $D^*$ (resp. $D^*_\eps$).

\begin{prop}\label{P:7.2}
It holds that
\begin{equation}\label{e:7.2}
K^*_{D_\eps}(z, \eps e^{i\theta})=  2K^*_D(z,0) \,
 \sin \theta \left(1 + O(\eps)\right),
\end{equation}
where $O(\eps)$ is a function whose absolute value is bounded by $c(z, \theta) \eps$
with $c(z, \theta)$ being uniformly bounded
in $0\le \theta\le \pi$ and $|z|>\eps_0$.
\end{prop}

\pf {\bf (i)}\ Put $\HH_\eps=\HH\setminus \overline B_\eps,\ \eps>0,$ and consider the Poisson kernel
$K_\HH(z,\zeta)=\frac{1}{\pi}\frac{\Im z}{|z-\zeta|^2}$
(resp. $K_{\HH_\eps}(z,\zeta)$)
of $\HH$ (resp. $\HH_\eps$).
Then
\begin{equation}\label{e:7.3}
K_{\HH_\eps}(z,\eps e^{i\theta})=2K_\HH(z,0)\sin \theta (1+ O(\eps)),
\quad \hbox{\rm uniformly in }
0\le \theta\le \pi,\ {\rm and}\ |z|>\eps_0.
\end{equation}
In fact, if we denote by
$Z^\HH=(Z_t,\zeta, \P_z^\HH)$
 (resp. $Z^{\HH_\eps}=(Z_t,\zeta, \P_z^{\HH_\eps})$)
the absorbing Brownian motion (ABM) on $\HH$ (resp. $\HH_\eps$), then according to \cite[p 50]{L1},
\begin{equation}\label{e:7.a}
\P_z^\HH(Z_{\sigma_{\partial B_R\cap\HH}}\in R e^{i\theta}d\theta)\;=\;\frac{2R}{\pi}\frac{\Im z}{|z|^2} \,
\sin \theta \left( 1+O ( R/|z| )\right) d\theta,\quad R>0,
\end{equation}
$O(R/|z|)$ being uniform in $R>0,\; z\in \HH\setminus B_R,$ which yields \eqref{e:7.3}.  We note that,
 for $z\in \HH_\eps,$
\begin{equation}\label{e:7.4}
\E_z^\HH[f(Z_{\sigma_\hc});\;\sigma_\hc<\infty]
=\eps\int_0^{\pi} K_{\HH_\eps}(z,\eps e^{i\theta})f(\eps e^{i\theta})d\theta.
\end{equation}

 Let $G_D(z, z')$
 be the Green function of $D$, namely, $0$-order resolvent density of ABM on $D$ (see \S 4 of [CFR]), and
$K_D(z,\zeta)$, $z\in D$, $\zeta\in \partial \HH$,
be the Poisson kernel of $D$.  The corresponding quantities for $D_\eps$ are designated by
$G_{D_\eps} (z, z')$ and
$K_{D_\eps}(z,\zeta)$, $z\in D$, $\zeta\in \partial (\HH\setminus B_\eps)$.
In view of \cite[\S 4]{CFR}, we have, for the outer normal $\n_\zeta$ at $\zeta.$
\begin{equation}\label{e:7.5}
K_D(z,\zeta)=-\frac12 \frac{\partial}{\partial \n_\zeta} G_D(z,\zeta),\ \
K_{D\eps}(z,\zeta)=-\frac12 \frac{\partial}{\partial \n_\zeta} G_{D_\eps}(z,\zeta).
\end{equation}

{\bf (ii)}\ We next show
\begin{equation}\label{e:7.6}
K_{D_\eps}(z,\eps e^{i\theta})=
2K_D(z,0) \;\sin \theta\; (1+O(\eps)),
\end{equation}
$O(\eps)$  being uniformly in $\theta\in [0,\pi]$ and $|z|>\eps_0.$
 By the strong Markov property of $Z^\HH$
\begin{equation}\label{e:7.7}
K_D(z,0)=K_\HH(z,0)-\E_z^\HH [ K_\HH(Z_{\sigma_K},0);\sigma_K<\infty].
\end{equation}
By the strong Markov property of $Z^{\HH_\eps}$ and then by \eqref{e:7.3} and  \eqref{e:7.7},
\begin{eqnarray*}
 K_{D_\eps}(z,\eps e^{i\theta})
&=& K_{\HH_\eps}(z,\eps e^{i\theta})
-\E_z^{\HH_\eps}\left[ K_{\HH_\eps}(Z_{\sigma_K},\eps e^{i\theta});\sigma_K<\infty\right]\\
&=& \left( 2K_{\HH}(z,0) -2\E_z^{\HH_\eps}[K_\HH(Z_{\sigma_K},0); \sigma_K<\infty]\right)  \sin \theta(1 +O(\eps)) \\
&=& ( 2K_D(z,0)+2A ) \sin \theta (1+O(\eps)),
\end{eqnarray*}
where
\begin{eqnarray*}
A&=&
\E_z^\HH [ K_\HH(Z_{\sigma_K},0),\sigma_K<\infty]
-\E_z^{\HH_\eps} [K_\HH(Z_{\sigma_K},0),\sigma_K<\infty]\\
&=&\E_z^\HH [K_\HH(Z_{\sigma_K},0),\sigma_K<\infty]
-\E_z^{\HH} [K_\HH(Z_{\sigma_K},0),\sigma_K<\sigma_\hc]\\
&=&\E_z^{\HH} [K_\HH(Z_{\sigma_K},0),\sigma_\hc <\sigma_K] \\
&=&\E_z^\HH\left[\E_{Z_{\sigma_\hc}}^{\HH}[K_\HH(Z_{\sigma_K},0),\sigma_K<\infty]; \sigma_\hc<\infty\right].
\end{eqnarray*}
Since $K_\HH(Z_{\sigma_K},0)\le C$ for some constant $C>0,$ $0\le A\leq C \P_z^\HH(\sigma_\hc<\infty)$.
It then follows from \eqref{e:7.3}, \eqref{e:7.4} and \eqref{e:7.7} that
$A \leq  4\eps C K_\HH(z,0)(1+O(\eps)) =O(\eps)K_D(z,0)$ uniformly for
$|z|>\eps_0$, proving (7.7).

\medskip
{\bf (iii)}\ Define
\[
\begin{cases}
\varphi_i(z)=\P_z^\HH\left(\sigma_K<\infty,\; Z_{\sigma_K}\in C_i\right),
\quad z\in D,\ 1\le i\le N , \\
\varphi_i^\eps(z)=\P_z^{\HH_\eps}\left(\sigma_K<\infty,\; Z_{\sigma_K}\in C_i\right),
\quad z\in D_\eps,\ 1\le i\le N.
\end{cases}
\]
By the strong Markov property of $Z^\HH$,
\begin{equation}\label{e:7.8}
\varphi_i^\eps(z)=\varphi_i(z)-
\E_z^\HH\left[ \varphi_i(Z_{\sigma_\hc}); \sigma_\hc<\infty\right] .
\end{equation}
Since $\varphi$ can be extended to be a differentiable function up to $\partial \HH$,
 we get from \eqref{e:7.3},\eqref{e:7.4} and \eqref{e:7.8}
\begin{equation}\label{e:7.9}
\varphi_i^\eps(z)=\varphi_i(z)+O(\eps^2)\quad\hbox{\rm uniformly for}\ |z|>\eps_0.
\end{equation}

{\bf (iv)}\ By virtue of \cite[(5.2)]{CFR}, the BMD-Poisson kernels $K^*_D$ and $K^*_{D_\eps}$ admit the expressions
\begin{equation}\label{e:7.10}
\begin{cases}
K^*_D(z,\zeta)=K_D(z,\zeta)-\sum_{i,j=1}^N b_{ij}\;\varphi_i(z)\;\frac{\partial}{\partial \n_\zeta}\varphi_j(\zeta),\\
K^*_{D_\eps}(z,\zeta)= K_{D_\eps}(z,\zeta)-\sum_{i,j=1}^N b_{ij}^\eps\;\varphi_i^\eps(z)\;\frac{\partial}{\partial \n_\zeta}\varphi_j^\eps(\zeta).
\end{cases}
\end{equation}
Here $(b_{ij})_{1\le i,j\le N}$
(resp. $B=(b_{ij}^\eps)_{1\le i,j\le N}$)
is the inverse matrix of $(a_{ij})_{1\le i,j\le N}$
(resp. $(a_{ij}^\eps)_{1\le i,j\le N}$)
whose entry is the period of $\varphi_i(z)$ (resp. $\varphi_i^\eps(z)$) around $C_j$, namely,
\begin{equation}\label{e:7.11}
a_{ij}=\int_\gamma \frac{\partial\varphi_i(\zeta)}{\partial \n_\zeta} ds(\zeta),
\qquad a_{ij}^\eps=\int_\gamma \frac{\partial\varphi_i^\eps(\zeta)}{\partial \n_\zeta} ds(\zeta),
\end{equation}
for any smooth Jordan curve $\gamma$ surrounding $C_j$ so that
${\rm ins}\; \gamma \supset C_j$ and $\overline{{\rm ins}\; \gamma}\cap C_k=\emptyset$
for $ k\neq j$.

We claim that
\begin{equation}\label{e:7.12}
b_{ij}^\eps=b_{ij}+O(\eps^2),\qquad 1\le i,j \le N.
\end{equation}
It suffices to show
\begin{equation}\label{e:7.13}
a_{ij}^\eps=a_{ij}+O(\eps^2),\qquad 1\le i,j \le N.
\end{equation}
By \eqref{e:7.8}, $\varphi_i^\eps(z)=\varphi_i(z)-h(z)$ for
$h(z)=\E_z^\HH\left[ \varphi_i(Z_{\sigma_\hc}); \sigma_\hc<\infty\right],$
For $\gamma$ in \eqref{e:7.11}, take $\wt\gamma$ surrounding $\gamma$
with $\wt\gamma\cap B_{\eps_0}=\emptyset$
and let $G={\rm ins} \wt\gamma.$\quad
Since $h$ is harmonic on $\overline G,$
$h(z)=\int_{\wt\gamma} p_G(z,\xi)h(\xi)s(d\xi),\; z\in G,$
for the Poisson kernel $p_G$ of $G.$
 Then
$\frac{\partial h(\zeta)}{\partial \n_\zeta}=\int_{\wt\gamma} \frac{\partial p_G(\zeta,\xi)}{\partial \zeta}h(\xi)s(d\xi),\; \zeta \in \gamma.$
As
$\sup_{\zeta\in \gamma,\;\xi\in \wt \gamma}\left|\frac{\partial p_G(\zeta,\xi)}{\partial \zeta}\right|$ is finite and $h(\xi)=O(\eps^2),\; \xi\in\wt \gamma,$ we have
$ \int_\gamma \frac{\partial h(\zeta)}{\partial \n_\zeta}ds(\zeta)=O(\eps^2),$
and hence \eqref{e:7.13} follows from \eqref{e:7.11}.

{\bf (v)}\ We finally show that
\begin{equation}\label{e:7.14}
-\frac12 \frac{\partial}{\partial \n_\zeta} \varphi_j^\eps(\zeta)\big|_{\zeta=\eps e^{i\theta}}
= -\frac{\partial}{\partial \n_\zeta} \varphi_j(\zeta)\big|_{\zeta=0} \; \sin \theta\; (1 +O(\eps)).
\end{equation}
We put $D^j=D\cup C_j$ and let $Z^{D^j}=\{Z_t, \P_z^{D^j},\ z\in D^j\}$  be the ABM on $D^j.$
$Z^{D^j}$ is obtained from $Z^\HH$ by killing upon hitting $\bigcup_{k\neq j}C_k.$  Then
$\varphi_j(z)=\P_z^{D^j}(\sigma_{C_j}<\infty)$ for $ z\in D^j$.
Let $G_{D^j}(z,z')$ be the Green function ($0$-order resolvent density) of $Z^{D^j}.$  By Corollary 3.4.3 and the $0$-order version of Lemma 2.3.10 of [FOT], there exists a finite measure $\nu$ concentrated on $C_j$ such that
\begin{equation}\label{e:7.15}
\varphi_j(z)=\int_{C_j} G_{D^j}(z,z')\nu(dz'),\quad z\in D^j.
\end{equation}

Analogously we put $D^j_\eps=D_\eps\cup C_j$ and let $Z^{D^j_\eps}=\{Z_t, \P_z^{D^j_\eps},\; z\in D^j\}$  be the ABM on $D^j_\eps.$  Then
$\varphi_j^\eps(z)=\P_z^{D^j_\eps}(\sigma_{C_j}<\infty),\; z\in D^j_\eps.$
By the strong Markov property of $Z^{D^j}$, we have
$$\varphi_j^\eps(z)=\varphi_j(z)-
\E_z^{D^j}\left[ \varphi_i(Z_{\sigma_\hc}); \sigma_\hc<\infty\right]
\quad \hbox{for } z\in D_\eps^j,
$$
and also, for the Green function $G_{D_\eps^j}(z,z')$ of $Z^{D_\eps^j},$
$$
G_{D_\eps^j}(z,z')= G_{D^j}(z,z')-
\E_z^{D^j}\left[ G_{D^j}(Z_{\sigma_\hc},z'); \sigma_\hc<\infty\right].
$$
Therefore we can deduce from \eqref{e:7.15} that
\begin{equation}\label{e:7.16}
\varphi_j^\eps(z)=\int_{C_j} G_{D^j_\eps}(z,z')\nu(dz'),\quad z\in D^j_\eps.
\end{equation}
Thus we have by \eqref{e:7.5}, \eqref{e:7.15} and \eqref{e:7.16} with $D^j$ and
$D^j_\eps$ in place of $D$, respectively, that
\begin{eqnarray}\label{e:7.17}
-\frac12 \frac{\partial}{\partial\n_\zeta}\varphi_j(\zeta)
&=& \int_{C_j} K_{D^j}(z',\zeta)\nu(dz'),\quad \zeta\in \partial \HH, \\
\label{e:7.18}
-\frac12 \frac{\partial}{\partial\n_\zeta}\varphi_j^\eps(\zeta)
&=& \int_{C_j} K_{D^j_\eps}(z',\zeta)\nu(dz'),\quad \zeta\in \partial (\HH\setminus B_\eps).
\end{eqnarray}
Consequently we get \eqref{e:7.14} from \eqref{e:7.17}, \eqref{e:7.18} and
\eqref{e:7.6} with $D^j_\eps$ and $D^j$ in place of $D_\eps$ and $D$, respectively.
We arrive at \eqref{e:7.2}
by combining \eqref{e:7.6}, \eqref{e:7.9}, \eqref{e:7.12} and \eqref{e:7.14} with \eqref{e:7.10}.
\qed

\noindent
{\bf Proof of Theorem \ref{T:7.1}}.
The proof is essentially along the line of the proof in \cite[Lemma 6.24]{D},
but with some simplifications by
using BMD instead of ERBM.

Without loss of generality, we may assume that $B_1\cap \HH\subset D.$  We write $S=\partial B_1\cap \HH$ and take $t$ so small that $F_t\subset B_1\cap \HH.$  Along with the ABM $Z^\HH=(Z_t, \zeta, \P^\HH_z)$ on $\HH$, we consider BMD $Z^*=(Z_t^*, \zeta^*, \P_z^*)$ on $D^*=D \cup \{c_1^*, \cdots, c_N^*\}$ and define
\[
\begin{cases}
M_1(t)=\int_0^\pi E_{r_te^{i\theta}}\left[ \Im Z_{\sigma_{F_t}}; \sigma_{F_t}<\infty\right]\; \sin \theta\;d\theta,
\quad  M_1^*(t)=\int_0^\pi E^*_{r_te^{i\theta}}\left[ \Im Z^*_{\sigma_{F_t}}; \sigma_{F_t}<\infty\right]\; \sin \theta\;d\theta\\
M_2(t)=\int_0^\pi E_{r_te^{i\theta}}\left[ \Im Z_{\sigma_{F_t}}; \sigma_{F_t}<\sigma_S\right]\; \sin \theta\;d\theta,
\quad  M_2^*(t)=\int_0^\pi E^*_{r_te^{i\theta}}\left[ \Im Z^*_{\sigma_{F_t}}; \sigma_{F_t}<\sigma_S\right]\; \sin \theta\;d\theta.
\end{cases}
\]

It is known (see \cite[Theorem 1.6.6]{C}) that, if a real valued function $u(z)$ defined on a planar domain $E$ with $K\subset E\subset \HH$ is continuous on $E$, constant on each slit $C_j$, harmonic on $E\setminus K$
and its period around each slit vanishes, then $u$ is harmonic with respect to BMD on $(E\setminus K)\cup \{c_1^*,\cdots, c_N^*\}.$
Let $g_t(z)$ be the canonical map from $D\setminus F_t.$
Since the function $h_t(z)=\Im(z-g_t(z))$ enjoys all these properties, it is BMD-harmonic on $(D\setminus F_t)\cup \{c_1^*, \cdots, c_N^*\}.$  As $h_t(z)=\Im z$ on $F_t$ and $h_t$ vanishes on $\partial \HH$ and at $\infty$, we have by the maximum principle
\[
h_t(z)=\E_z^*\left[\Im Z_{\sigma_{F_t}}^*; \sigma_{F_t}<\infty\right].
\]

We fix $R>0$ so large that $\HH\setminus B_R\subset D.$
By (3.2), $a_t=\lim_{y\to\infty}iy(g_t(iy)-iy)=\lim_{y\to\infty} yh_t(iy).$
By the strong Markov property of $Z^*$
 and \eqref{e:7.a}, we have for $y\ge R,$
\[yh_t(iy)=y\E_{iy}^*\left[h_t(Z^*_{\sigma_{\partial B_R\cap\HH}})\right]=\frac{2R}{\pi}\int_0^\pi h_t(Re^{i\theta})\sin \theta d\theta\cdot (1+O(R/y)),\]
yielding an expression
\begin{equation}\label{e:7.19}
a_t=\frac{2R}{\pi}\int_0^\pi h_t(Re^{i\theta})\;\sin \theta\; d\theta.
\end{equation}

Define $K^*_D(\infty,0)=\lim_{y\uparrow \infty} y K^*_D(iy,0).$
Since $K^*_D(z,0)$ is $Z^\HH$-harmonic on $\HH\setminus B_R,$ we have
from \eqref{e:7.a}
\[K^*_D(z,0)=\frac{2R}{\pi}\frac{\Im z}{|z|^2}\left[\int_0^\pi K^*_D(Re^{i\theta},0)\sin \theta d\theta\right](1+O(R/|z|),\]
 which implies that
\begin{equation}\label{e:7.b}
K^*_D(\infty,0)=\frac{2R}{\pi}\int_0^\pi K^*_D(Re^{i\theta},0)\sin\theta d\theta,\quad K^*_D(z,0)=\frac{\Im z}{|z|^2}K^*_D(\infty,0)+O(1/|z|^2).
\end{equation}
Notice that \eqref{e:7.b} holds not only for a standard slit domain $D$ but also for a more general domain $D=\HH\setminus \bigcup_{j=1}^N A_j$ where $\{A_j\}$ are mutually disjoint compact continua contained in $\HH.$

It follows from \eqref{e:7.19}, the strong Markov property of $Z^*$, and Proposition \ref{P:7.2} for $\eps=r_t$ that
\begin{eqnarray*}
 a_t &=& \frac{2R}{\pi}\int_0^\pi\int_0^\pi \E_{Re^{i\theta_1}}^*\left[ h_t(Z^*_{\sigma_{\partial B_{r_t}\cap\HH}}); \sigma_{\partial B_{r_t}\cap\HH}<\infty\right]\sin \theta_1 d\theta_1\\
&=& \frac{2Rr_t}{\pi}\int_0^\pi \int_0^\pi K^*_{D_{r_t}}(Re^{i\theta_1}, r_te^{i\theta_2}) h_t(r_te^{i\theta_2})d\theta_2 \sin \theta_1d\theta_1\\
&=& \frac{2R}{\pi}\int_0^\pi
 K^*_D(Re^{i\theta_1},0)
\sin \theta_1d\theta_1\cdot 2r_t\int_0^\pi h_t(r_te^{i\theta_2})\sin \theta_2 d\theta_2 (1+O(r_t)) ,
\end{eqnarray*}
which combined with \eqref{e:7.b} gives
\begin{equation}\label{e:7.c}
a_t= 2r_t K^*_D(\infty,0) M_1^*(t) (1+O(r_t) ) .
\end{equation}

We claim that
\begin{equation}\label{e:7.d}
K^*_D(\infty,0)=\frac{1}{\pi}.
\end{equation}
To this end, consider the conformal map $f(z)=-\frac{1}{z}$ from $\HH$ onto $\HH$ and the image domain $\wh D=f(D)=\HH\setminus \bigcup_{j=1}^N f(C_j)$ of $D$.
Let $K^*_{\wh D}(z,\xi)$ be the BMD-Poisson kernel of $\wh D.$

We first show for $K^*_{\wh D}(\infty,0)=\lim_{y\to\infty} y K^*_{\wh D}(iy,0)$ that
\begin{equation}\label{e:7.e}
K^*_{\wh D}(\infty,0)=\frac{1}{\pi}.
\end{equation}
Let $\wh \Psi(z,\xi)$ be the BMD-complex Poisson kernel of $\wh D$: $\Im \wh\Psi(z,\xi)=K^*_{\wh D}(z,\xi),\ \lim_{z\to\infty} \wh\Psi(z,\xi)=0.$  Let $b$ be a half of the BMD-domain constant defined by (6.1) for $\wh D$:
$ b=\lim_{z\to 0} (\pi \wh \Psi(z,0)+\frac{1}{z}).$
Define $\varphi_D(z)=\pi \wh\Psi(-\frac{1}{z},0)-b.$  Then $\Im \varphi_D(z)=\pi
K^*_{\wh D}(f(z),0)$ is constant on each slit $C_j$ and $\lim_{z\to\infty}(\varphi_D(z)-z)=\lim_{w\to 0}(\pi\wh \Psi(w,0)+\frac{1}{w})-b=0.$  Therefore $\varphi_D$ is a canonical map from $D$, and consequently $z=\varphi_D(z),\; z\in D$ so that $y=\pi K^*_{\wh D}(i/y, 0).$  On the other hand, we see from \eqref{e:7.b} for $\wh D^*$ and $z=i/y$ that
$K^*_{\wh D}(i/y, 0)=y K^*_{\wh D}(\infty,0)+O(y^2),$
and accordingly $y=\pi y K^*_{\wh D}(\infty,0)+O(y^2),$ yielding \eqref{e:7.e}.

We next prove
\begin{equation}\label{e:7.f}
K^*_D(\infty,0)=K^*_{\wh D}(\infty, 0),
\end{equation}
which together with \eqref{e:7.e} gives \eqref{e:7.d}.
  Let $G^*_D (z,z')$ (resp. $G^*_{\wh D} (w,w')$) be the Green function ($0$-order resolvent density) of BMD on $D^*$ (resp. $\wh D^*$).  Then we have
$ K^*_D(z,0)=\lim_{\eps\downarrow 0}\frac{1}{2\eps} G^*_D (z,i\eps)$ and
$ K^*_{\wh D}(z,0)=\lim_{\eps\downarrow 0}\frac{1}{2\eps}  G^*_{\wh D}(z,i\eps).$
The conformal invariance of BMD (\cite[remark 7.8.2]{CF1}) readily implies
the identity
$G^*_{\wh D} (w,w')=G^*_D (f^{-1}(w), f^{-1}(w'))$ of BMD-Green functions for $f(z)=-\frac{1}{z}.$
Accordingly, using the symmetry of $G^*_D$, we get
\begin{eqnarray*}
  K^*_{\wh D}(\infty,0) &=& \lim_{y\to\infty} y K^*_{\wh D}(iy,0)
=\lim_{y\to\infty}\lim_{\eps\downarrow 0}\frac{y G^*_{\wh D} (iy,i\eps)}{2\eps}
=\lim_{y\to\infty}\lim_{\eps\downarrow 0}\frac{y G^*_D (f^{-1}(iy),f^{-1}(i\eps))}{2\eps}\\
&=&\lim_{y\to\infty}\lim_{\eps\downarrow 0}\frac{y G^*_D (i/y, i/\eps)}{2\eps}
= \lim_{\eps\downarrow 0}\lim_{y\to\infty}\frac{y G^*_D (i/\eps, i/y)}{2\eps}
= \lim_{\eps\downarrow 0}\frac{1}{\eps} K^*_D(i/\eps, 0)=K^*_D(\infty,0).
\end{eqnarray*}

From \eqref{e:7.c} and \eqref{e:7.d}, we finally arrive at
\begin{equation}\label{e:7.20}
a_t=\frac{2}{\pi} r_t M_1^*(t)[1+O(r_t)],\quad r_t\to 0.
\end{equation}
An analogous formula holds for $a_t^0$ (\cite[p 70]{L1}):
\begin{equation}\label{e:7.21}
a_t^0= \frac{2}{\pi} r_t M_1(t).
\end{equation}

We now use Proposition \ref{P:7.2} again to verify that
\begin{equation}\label{e:7.22}
\lim_{t\downarrow 0} \frac{r_tM_1^*(t)}{t}\quad\hbox {\rm exists if and only if}
\quad
\lim_{t\downarrow 0} \frac{r_tM_2^*(t)}{t}\quad {\rm exists},
\end{equation}
and, in this case, they are equal.
In fact, we have for
$h_t(z)=\E_z^*\left[\Im Z_{\sigma_{F_t}}^*; \sigma_{F_t}<\infty\right],$
\[ h_t(r_te^{ i \theta})=
\E_{r_te^{i\theta}}^*\left[\Im Z_{\sigma_{F_t}}^*; \sigma_{F_t}<\sigma_S\right]+\E_{r_te^{i\theta}}^\HH\left[h_t(Z_{\sigma_S}); \sigma_S<\infty\right]\]
and so
\begin{equation}\label{e:7.23}
M_1^*(t)=M_2^*(t)+
\int_0^\pi \E_{r_te^{i\theta}}^\HH\left[h_t(Z_{\sigma_S}); \sigma_S<\infty\right]\sin \theta d\theta.
\end{equation}
By substituting \eqref{e:7.2} into
$h_t(z)=\int_0^\pi K^*_{D_{r_t}}(z,r_te^{i\eta}) h_t(r_te^{i\eta}) r_t d\eta,\; z\in S$, we obtain
\begin{equation}\label{e:7.24}
h_t(z)
= 2r_t K^*_D(z,0) M_1^*(t) \left( 1+O(r_t) \right),
\quad z\in S.
\end{equation}

If
$\lim_{t\downarrow 0}\frac{r_t M_1^*(t)}{t}=\gamma$
exists, then $\frac{h_t(z)}{t}$ is uniformly bounded in $t>0$ and $z\in S$ by \eqref{e:7.24}, and we conclude that
$\lim_{t\downarrow 0}\frac{r_t M_2^*(t)}{t}=\gamma$ by \eqref{e:7.23}.
Conversely, suppose
$\lim_{t\downarrow 0}\frac{r_t M_2^*(t)}{t}=\gamma'$ exists.  Since
$M_1^*(t)-M_2^*(t)\le C r_tM_1^*(t)$ for some constant $C>0$ from \eqref{e:7.23} and \eqref{e:7.24}, we get $M_1^*(t)\le 2 M_2^*(t)$ for sufficiently small $t>0.$  Hence $\limsup_{t \downarrow 0}\frac{r_tM_1^*(t)}{t}<\infty$ and we conclude that
$\lim_{t\downarrow 0}\frac{r_t M_1^*(t)}{t}=\gamma'$ just as above.

In the same way, we can use \eqref{e:7.6} to verify that
\begin{equation}\label{e:7.25}
\lim_{t\downarrow 0} \frac{r_tM_1(t)}{t}\quad\hbox {\rm exists if and only if}
\quad
\lim_{t\downarrow 0} \frac{r_tM_2(t)}{t}\quad {\rm exists},
\end{equation}
and, in this case, they are equal.
As $M_2^*(t)=M_2(t),$ the desired statement of Theorem \ref{T:7.1} follows from \eqref{e:7.20}, \eqref{e:7.21}, \eqref{e:7.22} and \eqref{e:7.25}.
\qed

\bigskip

{\bf Acknowledgement.} We thank Wai Tong Fan to help us turning hand drawn graphs into digital ones.

\vskip 0.3truein

\noindent {\bf Zhen-Qing Chen}

\smallskip \noindent
Department of Mathematics, University of Washington, Seattle,
WA 98195, USA

\noindent
E-mail: \texttt{zqchen@uw.edu}

\medskip

\noindent {\bf Masatoshi Fukushima}:

\smallskip \noindent
 Branch of Mathematical Science,
Osaka University, Toyonaka, Osaka 560-0043, Japan.

\noindent Email: {\texttt fuku2@mx5.canvas.ne.jp}


{\small

\baselineskip9.0pt


\begin{thebibliography}{99}


 \bibitem[A]{A} L. V.  Ahlfors, {\it Complex Analysis},  McGraw-Hill, 1979

\bibitem[BF1]{BF1} R. O. Bauer and R. M. Friedrich, Stochastic Loewner evolution in multiply connected domains,
{\it C. R. Acad. Sci. Paris, Ser I  \bf  339} (2004),  579-584. 

\bibitem[BF2]{BF2} R. O. Bauer and R. M. Friedrich, On radial stochastic Loewner
 evolution in multiply connected domains, {\it J. Funct. Anal. \bf 237} 
(2006), 565-588.

\bibitem[BF3]{BF3} R. O. Bauer and R.M. Friedrich, On chordal and bilateral SLE in multiply connected domains, 
{\it Math. Z. \bf 258} (2008), 241-265.


\bibitem[C]{C} Z.-Q. Chen, Browniam Motion with Darning, Lecture notes for
talks given at RIMS, Kyoto University, 2012.

\bibitem[CF1]{CF1} Z.-Q. Chen and M. Fukushima, {\it Symmetric Markov Processes, Time Changes, and Boundary Theory}, Princeton University Press, 2012.

 \bibitem[CF2]{CF2}  Z.-Q. Chen and M. Fukushima, One-point reflection,
 {\it Stochastic Process Appl. \bf 125}(2015), 1368-1393.

\bibitem[CFR]{CFR} Z.-Q. Chen, M. Fukushima and S. Rhode, Chordal Komatu-Loewner equation and Brownian motion with darning in multply connected domains.
{\it Trans. Amer. Math. Soc. \bf 368}(2016), 4065-4114.


 \bibitem[CFS]{CFS} Z.-Q. Chen, M, Fukushima and H. Suzuki, Stochastic Komatu-Loewner evolutions and SLEs.
 Preprint. 

\bibitem[CL]{CL} E.A. Coddington and N. Levinson, {\it Theory of Ordinary Differential Equations}, Krieger, 1984.

\bibitem[D]{D} S. Drenning,  Excursion reflected Brownian Motions and Loewner equations in multiply connected domains,  arXiv:1112.4123, 2011.


\bibitem[Dy]{Dy} E. B. Dynkin,
{\it Markov Processes}, Vol. I, Springer, 1965.


\bibitem[FK]{FK} M. Fukushima and H. Kaneko, On Villat's kernels and BMD Schwarz kernels in Komatu-Loewner equations. in: {\it Stochastic Analysis and Applications 2014}, Springer Proc. in Math. and Stat. Vol.100 (Eds) D. Crisan, B. Hambly, T. Zariphopoulous, 2014, pp 327-348.

\bibitem[FOT]{FOT} M. Fukushima, Y.Oshima and M. Takeda,
{\it Dirichlet Forms and Symmetric Markov Processes}, De Gruyter, 2nd edition, 2011.

\bibitem[GM]{GM} J. B. Garnett and D. E. Marshall, {\it  Harmonic Measure},
Cambridge University Press, 2005.



\bibitem[H]{H} P.  Hartman, {\it Ordinary Differential Equations}, John Wiley, 1964.

\bibitem[IW]{IW} N. Ikeda and S. Watanabe, {\it Stochastic Differential Equations and Diffusion Processes}, North-Holland/Kodansha, 1981.


\bibitem[K]{K} Y. Komatu, On conformal slit mapping of multiply-connected domains, {\it Proc. Japan Acad. \bf 26} (1950), 26-31.

\bibitem[L1]{L1} G. F.  Lawler, {\it Conformally Invariant Processes in the Plane}, Mathematical Surveys and Monographs, AMS, 2005.


\bibitem[L2]{L2} G.  F.  Lawler, The Laplacian-$b$ random walk and the Schramm-Loewner evolution, {\it Illinois J. Math. \bf 50} (2006), 701-746 (Special volume in memory of Joseph Doob).




\bibitem[LSW1]{LSW1} G. Lawler, O. Schramm and W. Werner,  Values of Brownian intersection exponents, I: Half-plane exponents, 
{\it Acta Mathematica \bf 187}(2001), 237-273.

\bibitem[LSW2]{LSW2} G. Lawler, O. Schramm and W. Werner, Conformal restriction: the chordal case, {\it J. Amer.Math. Soc. \bf 16} 
(2003), 917-955.


\bibitem[RY]{RY} D. Revuz and M. Yor, {\it Continuous Martingales and Brownian Motion}, Springer, 1999.


\bibitem[RW]{RW} L. C. G. Rogers and D. Williams, {\it Diffusions, Markov Processes and Martingales}, Vol. 1, Cambridge University Press, 1979.

\bibitem[RS]{RS} S. Rohde and O. Schramm, Basic properties of SLE, {\it Ann. Math. \bf 161} (2005), 879-920

\bibitem[S]{S} O. Schramm, Scaling limits of loop-erased random walks and uniform spanning trees, {\it Israel J. Math. \bf 118} (2000), 221-288.


\bibitem[W]{W} W. Werner, {\it Random Planar Curves and Schramm-Loewner Evolutions}, Lecture Notes in Math. {\bf 1840}, Springer, 2004.


\end{thebibliography}
}
\end{document}